\documentclass{tglat2e}
\theoremstyle{plain}
\newtheorem{theo}{Theorem}[section]
\newtheorem{prop}[theo]{Proposition}
\newtheorem{lemme}[theo]{Lemma}
\newtheorem{cor}[theo]{Corollary}

\begin{document}
\title{A characterization of groups of parahoric type}
\authors{Fran\c cois Court\`es 
\address Universit\'e de Poitiers\\Laboratoire de Math\'ematiques et Applications\\UMR 7348 du CNRS\\T\'el\'eport 2\\Boulevard Marie et Pierre Curie\\86962 Futuroscope Chasseneuil Cedex 
\email courtes@math.univ-poitiers.fr}
\maketitle
\begin{abstract}
Let $F$ be a local henselian nonarchimedean field of residual field $k$, and let $G$ be the group of $F$-points of a connected reductive group defined over $F$. It is well-known that the quotient of any parahoric subgroup of $G$ by its pro-unipotent radical is isomorphic to the group of $k$-points of a reductive group defined over $k$, and conversely. In this paper, we generalize this result by studying a class of linear algebraic groups named groups of parahoric type; we prove that, under certain conditions, any such group is $k$-isomorphic to the quotient of a parahoric subgroup of some reductive group over $F$ by its $d$-th congruence subgroup for a suitable $d$.
\end{abstract}
\vskip 1cm
\font\teuf=eufm10
\font\seuf=eufm7
\font\sseuf=eufm6
\newfam\euffam
\textfont\euffam=\teuf
\scriptfont\euffam=\seuf
\scriptscriptfont\euffam=\sseuf
\def \got{\fam\euffam}
%\font\tmsy=msym10
%\font\smsy=msym7% 
%\font\ssmsy=msym5%  
%\newfam\msyfam
%\textfont\msyfam=\tmsy
%\scriptfont\msyfam=\smsy
%\scriptscriptfont\msyfam=\ssmsy
%\def \mth{\fam\msyfam}
\def \mth{\mathbb}

\section{Introduction}

Let $k$ be a perfect field, and let $F$ be a henselian local field of discrete valuation with residual field $k$. Let $G$ be the group of $F$-points of a connected reductive group defined over $F$, and let $K$ be a parahoric subgroup of $G$. It is well-known (see for example \cite[section 2.5]{mp0} or \cite[section 3.2]{mp}) that the quotient of $K$ by its pro-unipotent radical $K^0$ is the group ${\mth{G}}$ of $k$-points of a connected reductive group over $k$; conversely, any such group is isomorphic to the quotient of a parahoric subgroup of a suitable reductive $F$-group by its pro-unipotent radical. This is a very useful tool in the theory of representations of $G$: if $(\pi,V)$ is a smooth representation of $G$ such that the subspace $V^{K^0}$ of elements of $V$ fixed by $K^0$ is nontrivial, the restriction of $\pi$ to $K$ factors into a representation of ${\mth{G}}$ on $V^{K^0}$, and the theory of representations of ${\mth{G}}$, which is usually easier, in particular when $k$ is finite (see for example \cite{car}), gives us lots of important informations about our representation $\pi$. With $k$ ifinite, that possibility has been, for example, widely used by the author in \cite{cou0} to study invariant distributions in the case of level $0$; it is also an important ingredient in the theory of types, since the study of level $0$ representations of a $p$-adic group can be reduced to the study of representations of finite groups of Lie type, which already contain all the relevant data (see for example \cite{bk} for general facts about types, and works of Morris (\cite{mor1}, \cite{mor2}), or the section $6$ of \cite{mp}, for the level $0$ case).

But of course smooth representations of $G$ do not always have nontrivial $K^0$-fixed vectors. A natural question which arises then is: what do the quotients of $K$ by normal bounded open subgroups smaller than $K^0$ look like~? It is not hard to prove, directly in the equal characeristic case and using \cite{witt} in the mixed-characteristic case (see proposition \ref{kkd}) that these groups are the groups of $k$-points of connected linear algebraic groups defined over $k$, but such groups are not reductive in general; is it possible to generalize the definition of a reductive group to groups of that kind~? Unfortunately, finding a definition which is as simple and elegant as in the reductive case while encompassing all of our quotients and only them seems to be really hard; that explains our choice of giving a definition of groups of parahoric type (first introduced in \cite{cou}) largely inspired by the definition of Bruhat-Tits' valued root data (see \cite[Section I, 6.2]{bt}): that kind of definition does not really qualify as simple and elegant, but is useful enough for our purposes.

An important particular case of such normal bounded open subgroups of $K$ are the congruence subgroups introduced in \cite{scst} by Schneider and Stuhler (see subsection $2.2$ for the definition of these subgroups). They define a filtration of a parahoric subgroup $K$ of $G$ by subgroups $K^d$, $d\in{\mth{N}}$, $K^0$ being the pro-unipotent radical of $K$, and then use them to define the Schneider-Stuhler coefficient systems over the Bruhat-Tits building of $G$: for every smooth representation $(\pi,V)$ of $G$ and every nonnegative integer $d$, the Schneider-Stuhler coefficient system of depth $d$ of $\pi$ is the collection of the subspaces $V^{K^d}$, with $K$ being any parahoric subgroup of $G$ and $K^d$ being its $d$-th congruence subgroup, of $K^d$-fixed points of $V$ on which $K$ acts; these subspaces are nothing else than representations of the groups $K/K^d$, which are finite when $k$ is a finite field. These coefficient systems have turned out to be a useful tool in the recently developed homological/cohomological theory of smooth representations of $p$-adic groups; they allow us for example to construct projective resolutions of representations of $p$-adic groups, which are then useful to compute their Ext groups  ((see \cite{scst0} and \cite{scst}). You can also csee for example the works of Broussous (\cite{br1}, \cite{br2}) about the coefficient systems themselves,, and of Opdam-Solleveld (\cite{os1},\cite{os2}) about extension functors.  In \cite{ms1} and \cite{ms2}, Meyer and Solleveld also use them to prove in the modular case some results of Harich Chandra (see \cite{hc1} or \cite{hc2}) about the trace formula.  Such an interest in coefficient systems shos that a more detailed study of the representations they involve can be of some interest.

Note also that in the particular case where $K$ is a special maximal parahoric subgroup of $G$, the quotient $K/K^d$ is $k$-isomorphic to the group $G(\mathcal{O}_F/\mathfrak{p}_F^{d+1})$, where $\mathcal{O}_F$ is the unit ring of $F$ and $\mathfrak{p}_F$ its maximal ideal, and that the representations of these groups have already been studied in particular by Lusztig (see \cite{lu}) and Stasinski (see \cite{aops}, \cite{st1}, \cite{st2}, \cite{st3}); unfortunately, that special case is not enough to be able to deal with the Schneider-Stuhler coefficient systems, which happen to involve all possible quotients of parahoric subgroups of our $p$-adic group $G$ by their congruence subgroups of some given depth. In that context, a more general study of these quotients is definitely of some interest.

It is easy to check that, apart from a few degenerate cases, all quotients of parahorics by these congruence subgroups match our definition (since it was tailor-made for that). But what about the converse~? In other words, given a group of parahoric type $\underline{\mth{G}}$ of given depth $d$ defined over a perfect field $k$, we want to find a henselian local field of discrete valuation $F$ with residual field $k$, a reductive group $\underline{G}$ defined over $F$ and, if $F_{nr}$ is the maximal unramified extension of $F$, a parahoric group $\underline{K}$ of the group $G_{nr}$ of $F_{nr}$-points of $\underline{G}$, stable by the Galois group $Gal(F_{nr}/F)\simeq Gal(\overline{k}/k)$, such that $\underline{\mth{G}}$ is $k$-isomorphic to $\underline{K}/\underline{K}^d$, where $\underline{K}^d$ is the $d$-th congruence subgroup of $\underline{K}$.

In this paper we restrict ourselves to the case where $\underline{\mth{G}}$ is quasisplit over $k$ (for the definition of a quasisplit group of parahoric type, see the beginning of section $6$: that definition is simply a generalization of the definition of a quasisplit reductive group.) It turns out that at least in this case, when the characteristic of $k$ is either $0$ or $p$ not too small, as soon as we have a suitable field $F$, it is possible, at least when the center of $\underline{\mth{G}}$ satisfies a structural condition (see section $6$),  to find suitable $\underline{G}$ and $\underline{K}$ as well. More precisely, for every element $\alpha$ of the absolute root system of $\underline{\mth{G}}$, we can put a ring structure on the root subgroup $\underline{\mth{U}}_\alpha$, and that ring $\underline{R}_\alpha$ is the quotient of the ring of integers of some henselian local field $F_\alpha$ by its ideal of elements of valuation at least either $d+1$ or $d$, depending whether $\alpha$ lies inside or outside the reductive part of $\underline{\mth{G}}$. Our main result (theorem \ref{main}) states that as soon as all the $\underline{R}_\alpha$, with $\alpha$ lying inside (resp. outside) the reductive part of $\underline{\mth{G}}$ are isomorphic to each other, then $\underline{\mth{G}}$ is isomorphic to a quotient of a parahoric subgroup. Unfortunately, there are groups of parahoric type which do not satisfy that condition, but it holds at least for large classes of groups, as we will see in section $7$ of this paper. We conjecture that our result extends to non-quasisplit groups of parahoric type satisfying the required condition. By slightly generalizing the definition of a group of parahoric type, we also expect to be able to prove a similar result with the Schneider-Stuhler filtrations being replaced with the full Moy-Prasad filtrations (see \cite{mp} for their definition), which allow for a finer study of smooth representations oof a $p$-adic group.

This paper is organized as follows. In section $2$, we recall the definition of a group of parahoric type and prove a few results which will be needed in the sequel. We also observe that, apart from the aforementioned degenrate cases, the quotients of parahoric subgroups of reductive groups defined over a local field by their $d$-th congruence subgroup actually match our definition.

In section $3$, we introduce the notion of a truncated valuation ring. We also attach a ring $\underline{R}_\alpha$ (resp. $R_\alpha$) to every root subgroup of a group of parahoric type (resp. to the group of its $k$-points), and prove that in most cases, either the rings attached to two different roots or suitable quotients of them are isomorphic.

In section $4$, we prove our main result (with some additional hypotheses) in the particular case of groups with a root system of rank $1$.

In section $5$, we make a detailed study of the commutator relations, in order to prove that the structural constants involved in them are similar to the Chevalley constants of reductive groups (see \cite{chev}). An additional difficulty compared to the case of reductive groups is that we want our constants to lie in (quotients of) a truncated valuation ring $R$ of residual field $k$, and since the Lie algebra of the $k$-group ${\mth{G}}$ is not a $R$-module, we are forced to work with ${\mth{G}}$ directly, which makes things a bit messier. We also study the action of a Cartan subgroup on the root subgroups.

In section $6$, we prove our main result in the general case.

In section $7$, we give a few examples of cases in which the conditions of the theorem are fulfilled, as well as an explicit example of a case in which they are not.

\section{Definition and first properties}

\subsection{Definition of a group of parahoric type}

Let $\underline{\mth{G}}$ be a connected algebraic group defined over some perfect field $k$, and let $d$ be a positive integer. Let $\underline{\mth{T}}$ be a maximal torus of $\underline{\mth{G}}$, and let $\underline{\mth{H}}$ be the centralizer of $\underline{\mth{T}}$ in $\underline{\mth{G}}$; $\underline{\mth{H}}$ is a Cartan subgroup of $\underline{\mth{G}}$. Let $X^*(\underline{\mth{T}})$ (resp. $X_*(\underline{\mth{T}}$)) be the group of characters (resp. cocharacters) of $\underline{\mth{T}}$.

Since we are working with groups over very general fields, we will not define any topology on the groups over $k$ other than the Zariski topology. In this paper, topology-related terms (closed, open, dense subsets of $\underline{\mth{G}}$ for example), when applied to algebraic groups over $k$, always refer to the Zariski topology.

For every algebraic group $\underline{\mth{L}}$, we denote by $R_u(\underline{\mth{L}})$ the unipotent radical of $\underline{\mth{L}}$.

Let $\underline{\Phi}$ be the set of weights of $\underline{\mth{G}}$ with respect to $\underline{\mth{T}}$, and let $\underline{\Psi}$ be the root system of the reductive group $\underline{\mth{G}}/R_u(\underline{\mth{G}})$, also with respect to $\underline{\mth{T}}$, identified to a maximal torus of $\underline{\mth{G}}/R_u(\underline{\mth{G}})$; $\underline{\Psi}$ is then a root system included in $\underline{\Phi}$.

Let $Lie(\underline{\mth{G}})$ be the Lie algebra of $\underline{\mth{G}}$. For every $\chi\in X^*(\underline{\mth{T}})$, let $\underline{\mth{U}}_\chi$ be the unique closed connected subgroup of $\underline{\mth{G}}$ whose Lie algebra is the eigenspace of $\chi$ for the action of $\underline{\mth{T}}$. When $\chi\neq 0$ and $\underline{\mth{U}}_\chi$ is nontrivial, it will be called a {\em root subgroup} of $\underline{\mth{G}}$. Since $Lie(\underline{\mth{U}}_\chi)$ is then a nilpotent subalgebra of $Lie(\underline{\mth{G}})$, $\underline{\mth{U}}_\chi$ is unipotent. Note also that $\underline{\mth{U}}_0=\underline{\mth{H}}$.

The group $\underline{\mth{G}}$ is said to be {\em of parahoric type} of depth $d$ if it satisfies the following conditions:

\begin{itemize}
\item {\bf (PT1)} $\underline{\Phi}$ is a reduced root system, and if $\underline{\Phi}^\vee$ is the corresponding set of coroots in $X_*(\underline{\mth{T}})$, the quadruplet $(X^*(\underline{\mth{T}}),\underline{\Phi},X_*(\underline{\mth{T}}),\underline{\Phi}^\vee)$ is a root datum;
\item {\bf (PT2)} for every $\alpha\in\underline{\Phi}$, the intersection of $R_u(\underline{\mth{G}})$ with the root subgroup $\underline{\mth{U}}_\alpha$ of $\underline{\mth{G}}$ with respect to $\alpha$ is of dimension $d$;
\item {\bf (PT3)} the group $\underline{\mth{H}}$ is abelian, of dimension $(d+1)r'$, where $r'$ is the rank of $\underline{\mth{T}}$, and for every $\alpha\in\underline{\Phi}$, the intersection $\underline{\mth{H}}_{\alpha,1}$ of its unipotent radical $R_u(\underline{\mth{H}})$ with the subgroup of $\underline{\mth{G}}$ generated by $\underline{\mth{U}}_\alpha$ and $\underline{\mth{U}}_{-\alpha}$ is of dimension $d$. Moreover, for every $\alpha,\alpha'$ such that $\alpha\neq\pm\alpha'$, $\underline{\mth{H}}_{\alpha,1}\cap\underline{\mth{H}}_{\alpha',1}=\{1\}$,  and the product of the $\underline{\mth{H}}_{\alpha,1}$ when $\alpha$ runs over $\underline{\Phi}$ is of dimension $dr$, where $r$ is the rank of $\underline{\Phi}$;
\item {\bf (PT4)} there exists a concave function $f_0$ from $\underline{\Phi}$ to ${\mth{Z}}$ and, for every $\alpha\in\underline{\Phi}$ and every integer $i\geq f_0(\alpha)$, a connected subgroup $\underline{\mth{U}}_{\alpha,i}$ of $\underline{\mth{U}}_\alpha$ satisfying the following conditions:
\begin{itemize}
\item {\bf (PT4a)} $\underline{\mth{U}}_{\alpha,f_0(\alpha)}=\underline{\mth{U}}_\alpha$;
\item {\bf (PT4b)} for every $i$, $\underline{\mth{U}}_{\alpha,i+1}\subset\underline{\mth{U}}_{\alpha,i}$, and if $\underline{\mth{U}}_{\alpha,i}$ is nontrivial, $dim(\underline{\mth{U}}_{\alpha,i+1})=dim(\underline{\mth{U}}_{\alpha,i})-1$;
\item {\bf (PT4c)} the commutator relations: for every $\alpha,\beta,i,j$ such that $\alpha+\beta\in\underline{\Phi}$, $i\geq f_0(\alpha)$ and $j\geq f_0(\beta)$, we have $[Lie(\underline{\mth{U}}_{\alpha,i}),Lie(\underline{\mth{U}}_{\beta,j})]=Lie(\underline{\mth{U}}_{\alpha+\beta,i+j})$;
\item {\bf (PT4d)} for every $\alpha,i,j,k$ such that $i,k\geq f_0(\alpha)$ and $j\geq f_0(-\alpha)$, $[Lie(\underline{\mth{U}}_{\alpha,i}),Lie(\underline{\mth{U}}_{-\alpha,j})]$ is of dimension $d+1-i-j$ and we have $[[Lie(\underline{\mth{U}}_{\alpha,i}),Lie(\underline{\mth{U}}_{-\alpha,j})],Lie(\underline{\mth{U}}_{\alpha,k})]=Lie(\underline{\mth{U}}_{\alpha,i+j+k})$.
\end{itemize}
\end{itemize}

As for reductive groups, the group $\underline{\mth{G}}$ is defined over $k$ if $\underline{\mth{T}}$ is defined over $k$.

Note that the above definition of a group of parahoric type is slightly different from the definition used in \cite{cou}.

Since all maximal tori of $\underline{\mth{G}}$ are conjugates, these properties do not depend on the choice of $\underline{\mth{T}}$.

Remember that a (${\mth{Z}}$-valued) concave function $f_0$ on $\underline{\Phi}$ is a function from $\underline{\Phi}$ to ${\mth{Z}}$ satisfying the following properties:
\begin{itemize}
\item for every $\alpha\in\underline{\Phi}$, $f_0(\alpha)+f_0(-\alpha)\geq 0$;
\item for every $\alpha,\beta\in\underline{\Phi}$ such that $\alpha+\beta\in\underline{\Phi}$, $f_0(\alpha+\beta)\leq f_0(\alpha)+f_0(\beta)$.
\end{itemize}

The function $f_0$ in the above definition is not uniquely determined by $\underline{\mth{G}}$, but we can  easily check that it is entirely determined by its values on the elements of any arbitrarily chosen set of simple roots of $\underline{\Phi}$, and that those values can also be chosen arbitrarily. In particular we have:

\begin{prop}\label{fazero}
Let $\alpha$ be any arbitrarily chosen element of $\underline{\Phi}$; it is possible to choose $f_0$ in such a way that $f_0(\alpha)=0$.
\end{prop}

\begin{proof}
Since $\underline{\Phi}$ is reduced by {\bf (PT1)}, we deduce from \cite[\S 1, proposition 11]{bou} that there exists a set of simple roots of $\underline{\Phi}$ which contains $\alpha$; the result then follows from the above remark. $\Box$
\end{proof}

(Note that it is not always possible to have $f_0$ ${\mth{Z}}$-valued and $f_0(\alpha)=0$ for every $\alpha\in\underline{\Psi}$ at the same time: for example, when $\underline{\Phi}$ is of type $B_2$ and $\underline{\Psi}$ is the subsystem of long roots of $\underline{\Phi}$, if we set $f_0(\alpha)=0$ for every long $\alpha$, then we must have $f_0(\alpha)=\frac 12$ for every short $\alpha$. This is a consequence of propositions \ref{f0psi} and \ref{sfun} (see below) and of the fact that the Weyl group $W$ of $\underline{\mth{G}}$ acts transitively on both $\underline{\Psi}$ and $\underline{\Phi}-\underline{\Psi}$ in that particular case.)

\begin{prop}\label{crs}
Let $\underline{\Phi}'$ be a closed root subsystem of $\underline{\Phi}$, and let $\underline{\mth{G}}'$ be the subgroup of $\underline{\mth{G}}$ generated by $\underline{\mth{H}}$ and the $\underline{\mth{U}}_\alpha$, $\alpha\in\underline{\Phi}'$. Then $\underline{\mth{G}}'$ is a group of parahoric type, and is defined over $k$ as soon as $\underline{\mth{G}}$ is.
\end{prop}

\begin{proof}
The proof is straightforward. As the concave function $f'_0$ from $\underline{\Phi}'$ to ${\mth{Z}}$ associated to $\underline{\mth{G}}'$ in {\bf [PT4)}, we can simply choose the restriction to $\underline{\Phi}'$ of $f_0$. $\Box$
\end{proof}

Note that it is also often possible to choose as $f'_0$ a ${\mth{Z}}$-valued concave function which is not the restriction to $\underline{\Phi}'$ of any ${\mth{Z}}$-valued concave function on $\underline{\Phi}$. For example, in the case we were considering above, setting $\underline{\Phi}'=\underline{\Psi}$, we can simply choose $f'_0=0$.

Set $X=X^*(\underline{\mth{T}})$ and $X^\vee=X_*(\underline{\mth{T}})$. We say that $\underline{\mth{G}}$ is {\em of semisimple parahoric type} if the following equivalent conditions are satisfied:
\begin{itemize}
\item the rank of $X$ is equal to the rank of $\underline{\Phi}$;
\item the rank of $X^\vee$ is equal to the rank of $\underline{\Phi}^\vee$;
\item $\underline{\Phi}$ generates $X\otimes{\mth{Q}}$ as a ${\mth{Q}}$-vector space;
\item $\underline{\Phi}^\vee$ generates $X^\vee\otimes{\mth{Q}}$ as a ${\mth{Q}}$-vector space.
\end{itemize}

That notion is a generalization of the notion of semisimple groups in the same way as the notion of group of parahoric type is a generalization of the notion of reductive groups. Of course a group which is both reductive and of semisimple parahoric type is semisimple.

\begin{prop}
For every $\alpha\in\underline{\Phi}$, $\underline{\mth{U}}_\alpha$ is an abelian group.
\end{prop}

\begin{proof}
The group $\underline{\mth{T}}$ acts on the subspace $[Lie(\underline{\mth{U}}_\alpha),Lie(\underline{\mth{U}}_\alpha)]$ of $Lie(\underline{\mth{G}})$ via the weight $2\alpha$. Since $\underline{\Phi}$ is reduced, it cannot contain $2\alpha$, hence that subspace must be trivial. $\Box$
\end{proof}

\begin{lemme}\label{nilr}
Assume $j\geq f_0(\alpha)$, $k\geq f_0(-\alpha)$ and $j+k=1$. Then the subalgebta $[Lie(\underline{\mth{U}}_{\alpha,j}),Lie(\underline{\mth{U}}_{-\alpha,k})]$ is the nilpotent radical $L_{\alpha,1}$ of the intersection of $Lie(\underline{\mth{H}})$ with the subalgebra of $Lie(\underline{\mth{G}})$ generated by $Lie(\underline{\mth{U}}_{\pm\alpha})$.
\end{lemme}

\begin{proof}
The subalgebra $[Lie(\underline{\mth{U}}_{\alpha,j}),Lie(\underline{\mth{U}}_{-\alpha,k})]$ is a nilpotent subalgebra of $Lie(\underline{\mth{H}})$ of dimension $d$, hence must be $L_{\alpha,1}$.
\end{proof}

\begin{prop}\label{hndep}
Let $i$ be a positive integer; for every $\alpha\in\underline{\Phi}$, the subalgebra $[Lie(\underline{\mth{U}}_{\alpha,j}),Lie(\underline{\mth{U}}_{-\alpha,k})]$ of $Lie(\underline{\mth{H}})$, with $j\geq f_0(\alpha)$, $k\geq f_0(-\alpha)$ and $j+k=i$, depends only on $i$.
\end{prop}

\begin{proof}
We will proceed by induction on $i$. When $i=1$, the result is simply lemma \ref{nilr}. Assume now $i>1$ and write, for every $i'<i$, $L_{\alpha,i'}=[Lie(\underline{\mth{U}}_{\alpha,j'}),Lie(\underline{\mth{U}}_{-\alpha,k'})]$ for suitable $j',k'$ such that $j'+k'=i'$; by induction hypothesis $L_{\alpha,i'}$ depends only on $i'$. Let $j,k,l$ be such that $j\geq f_0(\alpha)$, $k\geq f_0(-\alpha)$ and $j+k+l=i$; we now prove that we have:
\[[Lie(\underline{\mth{U}}_{\alpha,j+l}),Lie(\underline{\mth{U}}_{-\alpha,k})]=[Lie(\underline{\mth{U}}_{\alpha,j}),Lie(\underline{\mth{U}}_{-\alpha,k+l})],\] 
and the result follows immediately.

By {\bf (PT4d)}, we have:
\[[L_{\alpha,l},Lie(\underline{\mth{U}}_{\alpha,j})]=Lie(\underline{\mth{U}}_{\alpha,j+l});\]
\[[L_{\alpha,l},Lie(\underline{\mth{U}}_{-\alpha,k})]=Lie(\underline{\mth{U}}_{-\alpha,k+l}).\]
On the other hand, since $\underline{\mth{H}}$ is abelian, we have:
\[[Lie(\underline{\mth{U}}_{\alpha,j}),Lie(\underline{\mth{U}}_{-\alpha,k})],L_{\alpha,l}]=0.\]
The result then follows immediately from the Jacobi identity. $\Box$
\end{proof}

For every $\alpha\in\underline{\Phi}$ and every $u\in\underline{\mth{U}}_\alpha$ (resp. every $y\in Lie(\underline{\mth{U}}_\alpha)$), we call valuation of $u$ in $\underline{\mth{U}}_\alpha$ (resp. valuation of $y$ on $Lie(\underline{\mth{U}}_\alpha)$) and denote by $v(u)$ (resp. $v(y)$) the largest integer $v$ such that $u\in\underline{\mth{U}}_{\alpha,v}$ (resp $y\in Lie(\underline{\mth{U}}_{\alpha,v})$). By convention the valuation of the identity element is infinite.

\subsection{Parahoric groups and groups of parahoric type}

In this section, we give the definition of the $d$-th congruence subgroup (in the sense of Schneider-Stuhler) of a parahoric group of some algebraic group over a nonarchimedean local field with discrete valuation. These filtrations of parahoric subgroups were first introduced in \cite{scst0} for $\underline{G}=GL_n$, then in \cite{pr} for $\underline{G}$ quasi-simple and simply-connected, and finally generalized to any reductive group in \cite{scst}.  We now check that the quotients of parahoric subgroups by these subgroups are actually groups of parahoric type.

Let $\underline{G}$ be a connected reductive algebraic group defined over a henselian local field $F$ with discrete valuation and perfect residual field $k$, and split over the unramified closure $F_{nr}$ of $F$, and let $G_{nr}$ be the group of $F_{nr}$-points of $\underline{G}$. In the sequel, the topology we will be using on groups over $F$ is the usual analytic topology.

Let $\mathcal{B}_{nr}$ be the Bruhat-Tits building of $G_{nr}$. Let $\underline{T}$ be a $F_{nr}$-split maximal torus of $\underline{G}$, and let $\mathcal{A}$ be the apartment of $\mathcal{B}_{nr}$ associated to $\underline{T}$; $\mathcal{A}$ is isomorphic as an affine space to $(X_*(T_{nr})/X_*(Z_{nr}))\otimes{\mth{R}}$, where $Z_{nr}$ is the group of $F_{nr}$-points of the neutral component of the center $\underline{Z}$ of $\underline{G}$, and that isomorphism is canonical up to translation. We can thus identify $\mathcal{A}$ with $(X_*(T_{nr})/X_*(Z_{nr}))\otimes{\mth{R}}$ by setting the origin at some arbitrarily chosen point $x_0$ of $\mathcal{A}$; moreover, if we take as $x_0$ a special vertex of $\mathcal{B}_{nr}$; the walls of $\mathcal{A}$ are precisely the hyperplanes defined by equations of the form $\alpha(x)=i$ for some $\alpha\in\underline{\Phi}$ and some $i\in{\mth{Z}}$. Let $Q$ be the subgroup of $X^*(T_{nr})$ generated by $\underline{\Phi}$; for every $\alpha\in\underline{\Phi}$ and every $x\in\mathcal{A}$, we can set $\alpha(x)=<\alpha,x>$, where $<.,.>$ is the usual pairing between $X^*(T_{nr})$ and $X_*(T_{nr})$, whose restriction to $Q\times X_*(T_{nr})$ factors through $Q\times X_*(Z_{nr})$, extended to $Q\times\mathcal{A}$.

For every $\alpha\in\Phi$, let $U_\alpha$ be the root subgroup of $G_{nr}$ associated to $\alpha$, and for every $i\in{\mth{Z}}$, we denote by $U_{\alpha,i}$ the subgroup of $U_\alpha$ fixing every element $x$ of $\mathcal{A}$ such that $\alpha(x)\leq i$. We have $U_{\alpha,i}\supsetneq U_{\alpha,i+1}$ for every $i$, and the family $(U_{\alpha,i})_{i\in{\mth{Z}}}$ is a fundamental system of open neighborhoods of unity in $U_\alpha$ for the $p$-adic topology.

Let $G'$ be the derived subgroup of $G_{nr}$, and let $K_T$ be the maximal bounded subgroup of $T_{nr}$. The group $G_0=G'K_T$ is an open normal subgroup of $G_{nr}$, and we deduce from \cite[corollaries 2.19 and 2.21]{im} that $G$ is the semidirect product of $G_0$ by a discrete group (the group denoted by $\Omega$ in \cite{im}). For every part $S$ of $\mathcal{B}$, the {\em connected fixator} $K_S$ of $S$ in $G_{nr}$ is the subgroup of elements of $G_0$ which fix $S$ pointwise;  when $S$ is bounded, $K_S$ contains open subgroups of $T$ and the $U_\alpha$, hence is itself an open subgroup of $G_0$, hence also of $G$.

Let now $A$ be any facet of $\mathcal{A}$, and for every $\alpha\in\underline{\Phi}$, set:
\[f_A(\alpha)=Sup(\{\alpha(x)|x\in A\}).\]
By \cite[I. 6.4.3]{bt}, $f_A$ is a concave function; moreover, since every wall of $\mathcal{A}$ either contains $A$ or does not meet it, we have $f_A(\alpha)+f_A(-\alpha)\leq 1$ for every $\alpha\in\underline{\Phi}$. Let $K$ be the connected fixator of $A$ in $G_{nr}$; $K$ is a parahoric subgroup of $G$. Let now $f_A^0$ be the function defined by:
\[f_A^0:\alpha\in\underline{\Phi}\longmapsto 1-f_A(-\alpha).\] 
We will check that $f_A^0$ is also concave. For every $\alpha\in\underline{\Phi}$, we have:
\[f_A^0(\alpha)+f_A^0(-\alpha)=2-f_A(\alpha)-f_A(-\alpha)\geq 1\geq 0\]
since $f_A(\alpha)+f_A(-\alpha)\geq 1$. Now let $\alpha,\beta$ be elements of $\underline{\Phi}$ such that $\alpha+\beta\in\underline{\Phi}$; we have:
\[f_A^0(\alpha+\beta)=1-f_A(-\alpha-\beta)\leq 1-f_A(-\alpha)+f_A(\beta)=f_A^0\alpha)+f_A(\beta);\]
we are using here the fact that $f(-\alpha-\beta)+f(\beta)\leq f(-\alpha)$ by concavity of $f_A$. Moreover, since $f_A(\beta)+f_A(-\beta)=1$, $f_A(\beta)\leq f_A^0(\beta)$; hence $f_A^0(\alpha+\beta)\leq f_A^0(\alpha)+f_A^0(\beta)$ and $f_A^0$ is concave.

Let $K^0$ be the subgroup of $G_{nr}$ defined by:
\[K^0=Z^0\prod_{\alpha\in\underline{\Phi}}U_{\alpha,f_A(\alpha)},\]
with $Z^0$ being the pro-unipotent radical of $Z_{nr}$; $K^0$ is the pro-unipotent radical of $K$. Similarly, for every $i\geq 0$, set:
\[Z^i=\prod_{\xi\in X_*(Z)}\xi(1+\mathfrak{p}_{nr}^i),\]
where $X_*(Z)$ is the group of cocharacters of $Z$, and:
\[K^i=Z^i\prod_{\alpha\in\underline{\Phi}}U_{\alpha,f_A(\alpha)+i};\]
by \cite[I.2]{scst}, $K^i$ is an open normal subgroup of $K$, which we call the {\em $i$-th congruence subgroup} of $G$, and the $K^i$ form a fundamental system of open neighborhoods of unity in $G_{nr}$.

\begin{prop}\label{gfparah}
Assume the residual characteristic $p$ of $F$ and the root system $\underline{\Phi}$ of $\underline{G}$ satisfy one of the following conditions:
\begin{itemize}
\item $p>3$;
\item $p=3$ and $\underline{\Phi}$ has no irreducible component of type $G_2$;
\item $p=2$ and every irreducible component of $\underline{\Phi}$ is of type $A_n$ for some $n$.
\end{itemize}
Let $K$ be the parahoric subgroup of $G_{nr}$ associated to the facet $A$; assume $K$ is stable by the action of $Gal(F_{nr}/F)$ over $G_{nr}$. For every integer $d\geq 0$, let $K^d$ be the $d$-th congruence subgroup of $K$. Then $\underline{\mth{G}}=K/K^d$ is an algebraic group of parahoric type of depth $d$ defined over the residual field $k$ of $F$.
\end{prop}

\begin{proof}
The fact that $\underline{\mth{G}}$ satisfies {\bf (PT1)} to {\bf (PT3)} is an easy consequence of the definitions. Moreover, the concave function $f_0$ of {\bf (PT4)} is the function $f_A$ defined above, and for every $\alpha$ and every $i\geq f_0(\alpha)$, we obviously have $\underline{\mth{U}}_{\alpha,i}=U_{\alpha,i}/(U_{\alpha,i}\cap K^d)=U_{\alpha,i}/U_{\alpha,Sup(i,d+1-f_A(-\alpha)})$, hence {\bf (PT4a)} and {\bf (PT4b)} hold. Finally, when $p$ satisfies the required conditions, the commutator relations of {\bf (PT4c)} come from \cite[Theorem 1]{chev}, and the relations of {\bf (PT4d)} are easy to check directly. Hence {\bf (PT4)} holds as well. $\Box$
\end{proof}

Note that when $\underline{\Phi}$ and $p$ do not satisfy the condition of the above proposition, we can also deduce from \cite[Theorem 1]{chev} that {\bf (PT4c)} does not hold. Hence in most of these cases, $\underline{\mth{G}}$ cannot be a group of parahoric type although it comes from a quotient of a parahoric subgroup by some congruence subgroup. The author has chosen to keep these degenerate cases out of the definition, in order not to make it even more complicated than it already is.

\subsection{More properties of groups of parahoric type}

\begin{lemme}\label{ruconn}
The groups $R_u(\underline{\mth{H}})$ and, for every $\alpha\in\Phi$, $R_u(\underline{\mth{G}})\cap\underline{\mth{U}}_\alpha$, are connected.
\end{lemme}

\begin{proof}
Since the torus $\underline{\mth{T}}$ acts on $R_u(\underline{\mth{G}})$, that group is the product of the weight subgroups for that actions, which are $R_u(\underline{\mth{H}})$ and the intersections of $R_u(\underline{\mth{G}})$ with the $\underline{\mth{U}}_\alpha$, $\alpha\in\Phi$. Since $R_u(\underline{\mth{G}})$ is connected, these intersections must then be connected too. $\Box$
\end{proof}

\begin{prop}\label{d0red}
Assume $d=0$. Then $\underline{\mth{G}}$ is simply a reductive group.
\end{prop}

\begin{proof}
By lemma \ref{ruconn}, the intersections of $R_u(\underline{\mth{G}})$ with $\underline{\mth{H}}$ and the $\underline{\mth{U}}_\alpha$, $\alpha\in\Phi$,  are connected. On the other hand, for every $\alpha$, $R_u(\underline{\mth{G}})\cap\underline{\mth{U}}_\alpha$ is of dimension $d=0$ by {\bf (PT2)}, hence trivial, and since $dim(\underline{\mth{H}})=r=dim(\underline{\mth{T}})$, $R_u(\underline{\mth{H}})$ must be trivial too. The proposition follows. $\Box$
\end{proof}

\begin{prop}\label{f0psi}
Assume the characteristic $p$ of $k$ is not $2$. Let $\alpha$ be an element of $\underline{\Psi}$; then $f_0(\alpha)+f_0(-\alpha)=0$, and $\underline{\mth{U}}_\alpha\cap R_u(\underline{\mth{G}})=\underline{\mth{U}}_{\alpha,f_0(\alpha)+1}$.
\end{prop}

\begin{proof}
For every subspace $L$ of $Lie(\underline{\mth{G}})$, set $C_\alpha(L)=[[L,\underline{\mth{U}}_{-\alpha}],\underline{\mth{U}}_\alpha]$, and for every $i\geq 0$, define inductively $C^i_\alpha(L)$ by $C^0_\alpha(L)=L$ and $C^i_\alpha(L)=C_\alpha(C^{i-1}_\alpha(L))$ when $i>0$. By ({\bf PT4d)} and an obvious induction, we must have $C_\alpha^i(\underline{\mth{U}}_\alpha)=\underline{\mth{U}}_{\alpha,f_0(\alpha)+i(f_0(\alpha)+f_0(-\alpha))}$ for every $i$, hence if $f_0(\alpha)+f_0(-\alpha)>0$, $C_\alpha^i(\underline{\mth{U}}_\alpha)=0$ for $i$ large enough. On the other hand, since $\alpha\in\underline{\Psi}$, $\underline{\mth{U}}_\alpha$ has a nontrivial image in $Lie(\underline{\mth{G}}/R_u(\underline{\mth{G}}))$, and it is easy to check by induction that the image of $C_\alpha^i(\underline{\mth{U}}_\alpha)$ is equal to the image of $\underline{\mth{U}}_\alpha$ for every $i$, which leads to a contradiction if $f_0(\alpha)+f_0(-\alpha)>0$. Hence $f_0(\alpha)+f_0(-\alpha)$ must be zero.

Now we prove that $\underline{\mth{U}}_\alpha\cap R_u(\underline{\mth{G}})=\underline{\mth{U}}_{\alpha,f_0(\alpha)+1}$. Since these groups are both connected (by respectively lemma \ref{ruconn} and {\bf (PT4)}) and have the same dimension it is enough to prove that the latter is contained in the former; we will in fact prove the inclusion for their respective Lie algebras, which is equivalent since they are connected. Assume some element $u$ of $Lie(\underline{\mth{U}}_{\alpha,f_0(\alpha)+1})$ has a nontrivial image in $Lie(\underline{\mth{G}}/R_u(\underline{\mth{G}}))$; a simple computation in the Lie algebra $\mathfrak{sl}_2$ shows that there exists then $u'\in Lie(\underline{\mth{U}}_{-\alpha})$ such that $u_1=[[u,u'],u]$ also has a nontrivial image in $Lie(\underline{\mth{G}}/R_u(\underline{\mth{G}}))$. On the other hand, we have $u_1\in Lie(\underline{\mth{U}}_{f_0(\alpha)+2})$ by {\bf (PT4d)}, and by an easy induction (consider $u_2=[[u_1,u'],u]$ and so on), applying ({\bf PT4d)} again at each step, we find an element of the trivial space $Lie(\underline{\mth{U}}_{f_0(\alpha)+d+1})$ with a nontrivial image in $Lie(\underline{\mth{G}}/R_u(\underline{\mth{G}}))$ which is of course impossible. Hence the result. $\Box$
\end{proof}

\begin{prop}\label{sfun}
Assume $\alpha$ is an element of $\underline{\Phi}$ which does not belong to $\underline{\Psi}$. Then $f_0(\alpha)+f_0(-\alpha)=1$.
\end{prop}

\begin{proof}
If $d=0$, then $\underline{\mth{G}}$ is reductive, hence $\underline{\Psi}=\underline{\Phi}$ and there is nothing to prove; assume then $d>0$. Let $\alpha$ be any element of $\underline{\Phi}-\underline{\Psi}$; the dimension of $[Lie(\underline{\mth{U}}_\alpha),Lie(\underline{\mth{U}}_{-\alpha})]$, which is contained in the Lie algebra of $R_u(\underline{\mth{H}})$, is equal to $d$ by {\bf (PT3)}; on the other hand, it is also equal to $Sup(d+1-f_0(\alpha)-f_0(-\alpha),0)$ by {\bf (PT4d)}. Hence we must have $d+1-f_0(\alpha)-f_0(-\alpha)=d>0$ and the result follows immediately. $\Box$
\end{proof}

\begin{prop}\label{f0ab}
Let $\alpha,\beta$ be two elements of $\underline{\Phi}$ such that either $\alpha$ or $\beta$ lies in $\underline{\Psi}$. Assume $\alpha+\beta$ is also a root. Then $f_0(\alpha+\beta)=f_0(\alpha)+f_0(\beta)$.
\end{prop}

\begin{proof}
Assume for example $\alpha\in\underline{\Psi}$. We have $f_0(\alpha+\beta)\leq f_0(\alpha)+f_0(\beta)$ and $f_0(\beta)\leq f_0(\alpha+\beta)+f_0(-\alpha)$. On the other hand, since $\alpha\in\underline{\Psi}$, we have $f_0(\alpha)+f_0(-\alpha)=0$, hence:
\[f_0(\beta)\leq f_0(\alpha)+f_0(-\alpha)+f_0(\beta)=f_0(\beta).\]
Hence the first two inequalities are equalities and the proposition is proved. $\Box$
\end{proof}

\begin{prop}\label{f0ab2}
Let $\alpha,\beta$ be two elements of $\underline{\Phi}$ such that $\alpha+\beta\in\underline{\Phi}$. Assume in addition that none of the three lies in $\underline{\Psi}$.
\begin{itemize}
\item Assume $f_0(\alpha+\beta)=f_0(\alpha)+f_0(\beta)$. Then we also have $f_0(-\alpha)=f_0(\beta)+f_0(-\alpha-\beta)$ and $f_0(-\beta)=f_0(-\alpha-\beta)+f_0(\alpha)$.
\item Assume $f_0(\alpha+\beta)<f_0(\alpha)+f_0(\beta)$. Then we also have $f_0(-\alpha)<f_0(\beta)+f_0(-\alpha-\beta)$ and $f_0(-\beta)<f_0(-\alpha-\beta)+f_0(\alpha)$.
\end{itemize}
\end{prop}

\begin{proof}
First we prove the first assertion. We will only prove the first equality, the proof of the second one being similar. We have:
\[f_0(-\alpha)=1-f_0(\alpha)=1+f_0(\beta)-f_0(\alpha+\beta)=f_0(\beta)+f_0(-\alpha-\beta).\]
The second assertion is an immediate consequence of the first one. $\Box$
\end{proof}

\begin{prop}\label{f0ab3}
Let $\alpha,\beta$ be two elements of $\underline{\Phi}$ such that $f_0(\alpha+\beta)<f_0(\alpha)+f_0(\beta)$. Then $f_0(\alpha)=f_0(\alpha+\beta)+f_0(-\beta)$.

\end{prop}

\begin{proof}
We have:
\[f_0(\alpha)>f_0(\alpha+\beta)-f_0(\beta)\geq f_0(\alpha+\beta)+f_0(-\beta)-1;\]
since $f_0$ is concave and ${\mth{Z}}$-valued, we must then have $f_0(\alpha)=f_0(\alpha+\beta)+f_0(-\beta)$ and the proposition is proved.  $\Box$
\end{proof}

Now we look at the commutator relations. We can easily check in a similar way as in proposition 8.2.3 of \cite{spr} that for every $\alpha,\beta\in\underline{\Phi}$ such that $\alpha+\beta\in\underline{\Phi}$, the set of commutators $[\underline{\mth{U}}_\alpha,\underline{\mth{U}}_\beta]$ is contained in the product of the $\underline{\mth{U}}_{i\alpha+j\beta}$, with $i$ and $j$ being positive integers. Moreover, for every $u\in \underline{\mth{U}}_\alpha$ and every $u'\in\underline{\mth{U}}_\beta$, if we write $[u,u']=\prod_{i,j>0}u_{ij}$, with $u_{ij}\in\underline{\mth{U}}_{i\alpha+j\beta}$, $u_{11}$ does not depend on the choice of the order in which the product is taken. Hence the map $(u,u')\mapsto u_{11}$ from $\underline{\mth{U}}_\alpha\times\underline{\mth{U}}_\beta$ to $\underline{\mth{U}}_{\alpha+\beta}$ is a morphism of algebraic groups. We can then rewrite ${\bf (PT4c)}$ as follows:

\begin{itemize}
\item {\bf (PT4c')} for every $\alpha,\beta,i,j$ such that $\alpha+\beta\in\underline{\Phi}$, $i\geq f_0(\alpha)$ and $j\geq f_0(\beta)$, the projection of $[\underline{\mth{U}}_{\alpha,i},\underline{\mth{U}}_{\beta,j}]$ on $\underline{\mth{U}}_{\alpha+\beta,i+j}$ is surjective.
\end{itemize}

Now we consider the subsets of the form $\prod_{\chi\in\{\alpha,0,-\alpha\}}\underline{\mth{U}}_\chi$ of $\underline{\mth{G}}$, the product being taken in some given order; we first remark that they are of the same dimension as the subgroup $\underline{\mth{G}}_\alpha$ of $\underline{\mth{G}}$ they generate, hence open and dense in $\underline{\mth{G}}_\alpha$. Since they are obviously normalized by $\underline{\mth{H}}$, there exists an open dense subset $S$ of $\underline{\mth{U}}_\alpha\times\underline{\mth{U}}_{-\alpha}$ such that if $(u,u')$ belongs to $S$, we can write $[u,u']=\prod_{\chi\in\{\alpha,0,-\alpha\}}u_\chi$, with $u_\chi\in\underline{\mth{U}}_\chi$ for every $\chi$; moreover, since $\underline{\mth{H}}$ normalizes both $\underline{\mth{U}}_\alpha$ and $\underline{\mth{U}}_\alpha$, $S$ depends only on the relative position of $\alpha$ and $-\alpha$ for the order we have taken, and for every $(u,u')\in S$, the component $u_0$ also depends only on the relative position of $\alpha$ and $-\alpha$. On the other hand, if we write $[u,u']=u_0u_\alpha u_{-\alpha}$ and (assuming it is possible) $[u_\alpha,u_{-\alpha}^{-1}]=u'_0u'_\alpha u'_{-\alpha}$, we obtain:
\[[u,u']=u_0u_{-\alpha}^{-1}[u_\alpha,u_{-\alpha}^{-1}]^{-1}u_\alpha\]
\[=u_0u_{-\alpha}^{-1}u'_{-\alpha}{}^{-1}u'_\alpha{}^{-1}u'_0{}^{-1}u_\alpha.\]
from which we deduce that swapping $\alpha$ and $-\alpha$ in the definition of our order replaces $u_0$ by $u_0u'_0{}^{-1}$. Hence the subgroup $\underline{\mth{H}}_<$ of $\underline{\mth{H}}$ generated by the $u_0$ when considering an order such that $\alpha<-\alpha$ contains an open dense subset which is included in the subgroup $\underline{\mth{H}}_>$ we obtain when considering an order such that $\alpha>-\alpha$; since any open dense subset of $\underline{\mth{H}}_<$ generates $\underline{\mth{H}}_<$, we obtain that $\underline{\mth{H}}_<\subset\underline{\mth{H}}_>$. By symmetry, these subgroups must then be equal, which implies that our subgroup does not depend on the choice of the order at all. We obtain a similar result by replacing $\underline{\mth{U}}_\alpha$ (resp. $\underline{\mth{U}}_{-\alpha}$) by $\underline{\mth{U}}_{\alpha,i}$ for any $i$ (resp. $\underline{\mth{U}}_{-\alpha,j}$ for any $j$).

We can then rewrite {\bf (PT4d)} as follows:

\begin{itemize}
\item {\bf (PT4d')} for every $\alpha,i,j,k$ such that $i,k\geq f_0(\alpha)$ and $j\geq f_0(-\alpha)$, the projection of $[\underline{\mth{U}}_{\alpha,i},\underline{\mth{U}}_{-\alpha,j}]$ on $\underline{\mth{H}}$ generates a subgroup $\underline{\mth{L}}$ of $\underline{\mth{H}}$ of dimension $Sup(h-i-j,0)$, and is equal to $\underline{\mth{L}}$ when $i+j>0$; moreover, we have $[\underline{\mth{L}},\underline{\mth{U}}_{\alpha,k}]=\underline{\mth{U}}_{\alpha,i+j+k}$.
\end{itemize}

%Of course, when $\underline{\mth{G}}$ is defined over $k$, that property still holds with {\bf (PT4c}) and/or {\bf (PT4d)} being replaced with {\bf (PT4c')} and/or {\bf (PT4d')}.

We also obtain the following corollary to proposition \ref{hndep}:

\begin{cor}
Let $i$ be a positive integer; for every $\alpha\in\underline{\Phi}$, the subgroup of $\underline{\mth{H}}$ which is the projection of $[\underline{\mth{U}}_{\alpha,j},\underline{\mth{U}}_{-\alpha,k}]$, with $j\geq f_0(\alpha)$, $k\geq f_0(-\alpha)$ and $j+k=i$, depends only on $i$.
\end{cor}

For every $\alpha\in\underline{\Phi}$ and every $i>0$, we denote by $\underline{\mth{H}}_{\alpha,i}$ the subgroup given by the above corollary. We obviously have $\underline{\mth{H}}_{\alpha,i}=\underline{\mth{H}}_{-\alpha,i}$. Moreover, when $i=1$, we deduce immediately from lemma \ref{nilr} that $\underline{\mth{H}}_{\alpha,1}$ is the unipotent radical of $\underline{\mth{H}}_\alpha$, which matches the definition of $\underline{\mth{H}}_{\alpha,1}$ given in the statement of {\bf (PT3)}.

For every $\alpha\in\underline\Phi$, if $\alpha^\vee$ is the corresponding coroot in $\underline\Phi^\vee$ and if $\underline{\mth{T}}_\alpha$ is the image of $\alpha^\vee$ in $\underline{\mth{T}}$, let $\underline{\mth{H}}_\alpha$ be the subgroup of $\underline{\mth{H}}$ generated by $\underline{\mth{T}}_\alpha$ and $\underline{\mth{H}}_{\alpha,1}$. Since for every $x\in k^*$, $(-\alpha^\vee)(x)=\alpha^\vee(x^{-1})$, we also have $\underline{\mth{T}}_{-\alpha}=\underline{\mth{T}}_\alpha$, hence $\underline{\mth{H}}_{-\alpha}=\underline{\mth{H}}_\alpha$.

\subsection{The group of $k$-points of $\underline{\mth{G}}$}

Now assume $\underline{\mth{G}}$ is defined over $k$. In this subsection, we prove some useful results about the group of $k$-points ${\mth{G}}$ of $\underline{\mth{G}}$. Note that the corresponding results for $\underline{\mth{G}}$ are simply particular cases of the results below (with $k$ algebraically closed). Remember that we are assuming that $k$ is perfect.

Let $\underline{\mth{L}}$ be any linear algebraic group defined over $k$, and let ${\mth{L}}$ be the group of its $k$-points. Since $k$ is perfect, by \cite[14.4.5]{spr}, $R_u(\underline{\mth{L}})$ is defined over $k$; in the sequel, by a slight abuse of notation, we will denote by $R_u({\mth{L}})$ the group of $k$-points of $R_u(\underline{\mth{L}})$. Obviously, we have $R_u({\mth{L}})=R_u(\underline{\mth{L}})\cap{\mth{L}}$.

Let ${\mth{G}}$ (resp. ${\mth{T}}$, ${\mth{H}}$) be the group of $k$-points of $\underline{\mth{G}}$ (resp. $\underline{\mth{T}}$, $\underline{\mth{H}}$), let $\Phi$ (resp. $\Psi$) be the relative root system of ${\mth{G}}$ (resp. ${\mth{G}}/R_u({\mth{G}})$), and for every $\alpha\in\Phi$, let ${\mth{U}}_\alpha$ be the corresponding root subgroup of ${\mth{G}}$.

We say that $\underline{\mth{G}}$ is {\em split} over $k$, or $k$-split, if it contains a $k$-split maximal torus. It implies in particular that $\Phi=\underline{\Phi}$, $\Psi=\underline{\Psi}$ and for every $\alpha\in\Phi$, ${\mth{U}}_\alpha$ is the group of $k$-points of $\underline{\mth{U}}_\alpha$.

Now we prove that for solvable groups, our definition of a split group is consistent with the usual definition of a split solvable group:

\begin{prop}
Assume $\underline{\mth{G}}$ is solvable. Then $\underline{\mth{G}}$ is split according to the above definition if and only if it  is split according to \cite[15.1]{bor}.
\end{prop}

\begin{proof}
For every $i\leq d$, consider the subgroup $\underline{\mth{G}}_i$ of $\underline{\mth{G}}$ defined by:
\[\underline{\mth{G}}_i=\underline{\mth{H}}_i\prod_{\alpha\in\overline{\Phi}}\underline{\mth{U}}_{f_0(-\alpha)-i-1},\]
with ${\mth{H}}_i$ being the product of the ${\mth{H}}_{\alpha,i}$, $\alpha\in\underline{\Phi}$. Then $\underline{\mth{G}}_d$ is trivial, $\underline{\mth{G}}/\underline{\mth{G}}_0$ is the product of a central subgroup and of the reductive quotient of $\underline{\mth{G}}$, which is a torus since $\underline{\mth{G}}$ is solvable, and isomorphic to $\underline{\mth{T}}$ hence split, and for every $i\in\{0,\dots,d-1\}$, $\underline{\mth{G}}_i/\underline{\mth{G}}_{i+1}$ is a direct product of one-dimensional abelian groups which are all $k$-isomorphic to either $\underline{\mth{G}}_a$ or $\underline{\mth{G}}_m$. 

Conversely, assume there exists subgroups $\underline{\mth{G}}=\underline{\mth{L}}_0\supset\underline{\mth{L}}_1\supset\dots\supset\underline{\mth{L}}_r=\{0\}$ of $\underline{\mth{G}}$ suchthat for every $i$, $\underline{\mth{L}}_i/\underline{\mth{L}}_{i+1}$ is $k$-isomorphic to either $\underline{\mth{G}}_a$ or $\underline{\mth{G}}_m$. Then the reductive quotient of $\underline{\mth{G}}$ must be solvable and split too, hence is a split torus, which proves the proposition. $\Box$
\end{proof}

From now on until part $5$ included, we will make the following assumptions:
\begin{itemize}
\item $\underline{\mth{G}}$ and $\underline{\mth{T}}$ are split over $k$;
\item the characteristic $p$ of $k$ is not $2$, and if $\Phi$ has any irreducible component of type $G_2$, $p$ is not $3$ either;
\item the map $x\mapsto x^2$ is surjective on $k$ (which means, since $p\neq 2$, that every element of $k$ has two opposite square roots);
\item if $P^\vee$ is the subgroup of the elements $\xi$ of $X^*(\underline{\mth{T}}^\vee)\otimes{\mth{Q}}$ such that $<\alpha,\xi>$ is an integer for every $\alpha\in\underline{\Phi}$ and $Q^\vee$ is the subgroup of $X_*(\underline{\mth{T}})$ generated by $\underline{\Phi}^\vee$, the inverse image of every element of $k$ by the map $x\mapsto x^{[P^\vee:Q^\vee]}$ is also contained in $k$.
\end{itemize}

It is easy to check that in most cases, the third condition implies the fourth one. In fact, when $\underline{\Phi}$ is irreducible, the only cases where it does not are the cases $A_n$, $n+1$ not being a power of $2$, and $E_6$.

\begin{lemme}\label{kinf}
With the above hypotheses, $k$ is infinite.
\end{lemme}

\begin{proof}
Asssume $k$ is finite. Since $p\neq 2$, we have $-1\neq 1$ in $k$, hence $1$ has two distinct square roots. The map $x\mapsto x^2$ is then not injective, hence not surjective either, which leads to a contradiction. Hence $k$ must be infinite. $\Box$
\end{proof}

Since $k$ is infinite, by \cite[proposition 18.3]{bor}, ${\mth{G}}$ is dense in $\underline{\mth{G}}$ and it is not hard to check that the following properties respectively imply {\bf (PT4c')} and {\bf (PT4d')}:

\begin{itemize}
\item {\bf (PT4c'')} for every $\alpha,\beta,i,j$ such that $\alpha+\beta\in\underline{\Phi}$, $i\geq f_0(\alpha)$ and $j\geq f_0(\beta)$, the projection of $[{\mth{U}}_{\alpha,i},{\mth{U}}_{\beta,j}]$ on ${\mth{U}}_{\alpha+\beta,i+j}$ is surjective,
\item {\bf (PT4d'')} for every $\alpha,i,j,k$ such that $i,k\geq f_0(\alpha)$ and $j\geq f_0(-\alpha)$, the projection of ${\mth{U}}_{\alpha,i},{\mth{U}}_{-\alpha,j}$ on ${\mth{H}}$ generates a subgroup ${\mth{L}}$ of ${\mth{H}}$ of dimension $Sup(h-i-j,0)$, and is equal to ${\mth{L}}$ when $i+j>0$; moreover, we have $[{\mth{L}},{\mth{U}}_{\alpha,k}]={\mth{U}}_{\alpha,i+j+k}$.
\end{itemize}

We will prove later (corollary \ref{priisec} and proposition \ref{prisec2}) that these implications are in fact equivalences. For the moment, we simply prove our statements for a group satisfying {\bf (PT4c'')} and {\bf (PT4d'')} instead of {\bf (PT4c')} and {\bf (PT4d')}

For every $\alpha\in\Phi$, let ${\mth{H}}_\alpha$ (resp. ${\mth{T}}_\alpha$) be the group of $k$-points of $\underline{\mth{H}}_\alpha$ (resp. $\underline{\mth{T}}_\alpha$). For every positive integer $i$, we also denote by ${\mth{H}}_{\alpha,i}$ the group of $k$-points of $\underline{\mth{H}}_{\alpha,i}$.

\begin{lemme}\label{dpta}
With the same hypotheses on $k$ as above, let $S$ be a generating subset of $\Phi$ whose cardinality is the rank of $\Phi$. Then ${\mth{T}}$ is the direct product of the ${\mth{T}}_\alpha$, $\alpha\in S$.
\end{lemme}

\begin{proof}
Let $t$ be any element of ${\mth{T}}$; we have a unique decomposition $t=\prod_{\alpha\in S}t_\alpha$, where for every $\alpha\in S$, $t_\alpha\in\underline{\mth{T}}_\alpha$, and we thus only have to prove that $t_\alpha\in{\mth{T}}_\alpha$ for every $\alpha$.

Since $X_*(\underline{\mth{T}}))$ is contained in $P^\vee$, $[P^\vee:Q^\vee]X_*(\underline{\mth{T}}))$ is contained in $[P^\vee:Q^\vee]P^\vee\subset Q^\vee\subset X_*(\underline{\mth{T}})$. We deduce from this that since $\prod_{\alpha\in S}t_\alpha$ is an element of ${\mth{T}}$, we have $t_\alpha^{[P^\vee:Q^\vee]}\in{\mth{T}}$ for every $\alpha$.

Hence for every $\alpha\in S$, there exists $x_\alpha\in k$ such that $t_\alpha^{[P^\vee:Q^\vee]}=\alpha^\vee(x_\alpha)$, hence $t_\alpha=\alpha^\vee(y_\alpha)$, with $y_\alpha$ being a $[P^\vee:Q^\vee]$-th root of $x_\alpha$. By the hypotheses we have made, we then have $y_\alpha\in k$, which implies $t_\alpha\in{\mth{T}}$. The lemma is now proved. $\Box$
\end{proof}

\begin{lemme}\label{hsq}
Let $\alpha$ be an element of $\Phi$. Then every element of ${\mth{H}}_\alpha$ has two square roots in ${\mth{H}}_\alpha$.
\end{lemme}

\begin{proof}
Consider the decomposition ${\mth{H}}_\alpha={\mth{T}}_\alpha{\mth{H}}_{\alpha,1}$; for every $h\in{\mth{H}}_\alpha$, the decomposition $h=tu$, with $t\in{\mth{T}}_\alpha$ and $u\in{\mth{H}}_{\alpha,1}$, is simply the Jordan decomposition of $h$. This implies in particular that if $h'=t'u'$ is such that $h'^2=h$, then we have $t'^2=t$ and $u'^2=u$, and conversely.

The group ${\mth{T}}_\alpha$ is isomorphic to $k^*$, which implies that $t$ admits two square roots in ${\mth{T}}_\alpha$. On the other hand, for every $i\geq 1$, ${\mth{H}}_{\alpha,i}/{\mth{H}}_{\alpha,i+1}$ is isomorphic to the additive group $k$; since the characteristic of $k$ is not $2$, the map $x\mapsto x^2$ is an automorphism of ${\mth{H}}_{\alpha,i}/{\mth{H}}_{\alpha,i+1}$ for every $i$, hence also an automorphism of ${\mth{H}}_{\alpha,1}$. Therefore, $u$ admits exactly one square root in ${\mth{H}}_{\alpha,1}$. The result follows.  $\Box$
\end{proof}

\begin{prop}\label{htrans2}
For every $\alpha\in\Phi$ and every integer $i$, the action of ${\mth{H}}_\alpha$ by conjugation on the set of elements of ${\mth{U}}_\alpha$ of valuation $i$ is transitive.
\end{prop}

\begin{proof}
We can of course assume that $i$ is such that ${\mth{U}}_\alpha$ actually contains elements of valuation $i$. To simplify the notations, using proposition \ref{fazero}, we will also assume $f_0(\alpha)=0$.

Assume first $i=0$. Let $u_1,u_2$ be two elements of ${\mth{U}}_\alpha$ of valuation $0$; we want to prove that there exists $h\in{\mth{H}}_\alpha$ such that $hu_1h^{-1}=u_2$. Since ${\mth{T}}$ acts transitively by conjugation on the non-identity elements of the one-dimensional quotient group ${\mth{U}}_\alpha/{\mth{U}}_{\alpha,1}$, by conjugating $u_1$ by an element of ${\mth{T}}$, we can assume that $u_1\in{\mth{U}}_{\alpha,1}u_2$.

Consider the application $(h,u)\mapsto[h,u]$ from ${\mth{H}}_{\alpha,1}\times{\mth{U}}_{\alpha,0}$ to ${\mth{U}}_{\alpha,1}$. Since by {\bf (PT4c'')} and {\bf (PT4d'')} the image of ${\mth{H}}_{\alpha,2}\times{\mth{U}}_{\alpha,0}$ is contained in ${\mth{U}}_{\alpha,2}$, this application induces an application from ${\mth{H}}_{\alpha,1}/{\mth{H}}_{\alpha,2}\times{\mth{U}}_{\alpha,0}/{\mth{U}}_{\alpha,1}$ to ${\mth{U}}_{\alpha,1}/{\mth{U}}_{\alpha,2}$, which we will denote by $(\overline{h},\overline{u})\mapsto[\overline{h},\overline{u}]$.

Now we identify both ${\mth{H}}_{\alpha,1}/{\mth{H}}_{\alpha,2}$ and ${\mth{U}}_{\alpha,0}/{\mth{U}}_{\alpha,1}$ to the group ${\mth{G}}_a$ of $k$-points of $\underline{\mth{G}}_a$, the additive group of dimension $1$, via suitable isomorphisms, and we consider the application $(\overline{h},\overline{u})\mapsto[\overline{h},\overline{u}]$ as an algebraic map from ${\mth{G}}_a^2$ to ${\mth{G}}_a$.

\begin{lemme}\label{huhu}
There exists a nonzero constant $c$ such that we have $[\overline{h},\overline{u}]=c\overline{h}\overline{u}$ for every $\overline{h},\overline{u}\in{\mth{G}}_a$.
\end{lemme}

\begin{proof}
For every $h,h'\in{\mth{H}}_{\alpha,1}$ and every $u\in{\mth{U}}_\alpha$, we have:
\[[hh',u]=hh'uh'^{-1}h^{-1}u^{-1}=h[h',u]h^{-1}[h,u]=[h,[h',u]][h',u][h,u],\]
and since $[h,[h',u]]\in{\mth{U}}_{\alpha,2}$ and ${\mth{U}}_\alpha$ is abelian, in the quotient group identified to the additive group ${\mth{G}}_a$ we obtain $[\overline{h}\overline{h'},\overline{u}]=[\overline{h},\overline{u}]+[\overline{h'},\overline{u}]$. Similarly, for every $h,u,u'$, we have $[h,uu']=[h,u]u[h,u']u^{-1}$, from which we obtain in the quotient group, identified to ${\mth{G}}_a$  again:
\[[\overline{h},\overline{u}\overline{u'}]=[\overline{h},\overline{u}]+[\overline{h},\overline{u'}].\]
Finally, since $t$ and $h$ commute, we have for every $t\in{\mth{T}}$:
\[[h,Ad(t)u]=Ad(t)[h,u].\]
Since both $u$ and $[h,u]$ belong to ${\mth{U}}_\alpha$, $Ad(t)$ acts on ${\mth{U}}_{\alpha,0}/{\mth{U}}_{\alpha,1}$ identified to ${\mth{G}}_a$ by multiplication by $\alpha(t)$; since for every $a\in k^*$, there exists $t\in{\mth{T}}$ such that $\alpha(t)=a$, we have shown that the map $(\overline{h},\overline{u})\mapsto[\overline{h},\overline{u}]$ is $k$-linear in $\overline{u}$ and ${\mth{Z}}$-linear in $\overline{h}$; that map must then be of the form $(\overline{h},\overline{u})\mapsto cP(\overline{h})\overline{u}$, where $P(X)=X$ when $char(k)=0$ and $P$ is a polynomial whose nonzero terms are all of degree a power of $p$ when $char(k)=p>0$, and $c$ is an element of $k$ which is nonzero because {\bf (PT4d'')} implies that $(\overline{h},\overline{u})\mapsto[\overline{h},\overline{u}]$ is surjective. The lemma is now proved when $char(k)=0$; we will thus assume $char(k)=p>0$ in the rest of the proof.

Consider now the group ${\mth{U}}_{-\alpha,1}/{\mth{U}}_{-\alpha,2}$, also identified with ${\mth{G}}_a$; by the same reasoning as above, the map which associates to $u^-\in{\mth{U}}_{-\alpha,1}$ and $u\in{\mth{U}}_\alpha$ the component in ${\mth{H}}_{\alpha,1}$ of $[u^-,u]$ (which is an element of the group ${\mth{H}}_{\alpha,1}{\mth{U}}_{\alpha,1}{\mth{U}}_{-\alpha,2}$) induces in the quotient groups a map $(\overline{u}^-,\overline{u})\mapsto[\overline{u}^-,\overline{u}]$ which is ${\mth{Z}}$-bilinear and such that for every $\overline{u}^-,\overline{u}$ and every nonzero constant $\lambda$, we have:
\[[\lambda\overline{u}^-,\frac 1\lambda\overline{u}]=[\overline{u}^-,\overline{u}].\]
Hence the map is of the form $(\overline{u}^-,\overline{u})\mapsto Q(\overline{u}^-\overline{u})$, where $Q$ is a polynomial whose nonzero terms are all of degree a power of $p$. Finally we have, for every $u^-,u,u'$:
\[[[u^-,u],u']=u^-u(u^-)^{-1}u'u^-u^{-1}(u^-)^{-1}u'^{-1}=Ad(u^-u)[[(u^-)^{-1},u'],u^{-1}],\]
hence in the (additive) quotient group:
\[[[\overline{u}^-,\overline{u}],\overline{u'}]=[[-\overline{u}^-,\overline{u'}],-\overline{u}],\]
which can be rewritten as:
\[P(Q(\overline{u}^-\overline{u}))\overline{u'}=-P(Q(-\overline{u}^-\overline{u'}))\overline{u}.\]
From {\bf (PT4d'')} we know that both $P$ and $Q$ are nonzero. Since the right-hand side of the equation is of degree $1$ in $\overline{u}$, in the left-hand side, the polynomial $P\circ Q$ must also be of degree $1$. Hence both $P$ and $Q$ are monomials of degree $1$ and the lemma is proved. $\Box$
\end{proof}

We deduce from lemma \ref{huhu} that the map $\overline{h}\mapsto[\overline{h},\overline{u}]$ is surjective for every nonzero $\overline{u}$, which implies that there exists $h\in{\mth{H}}_{\alpha,1}$ such that $[h,u_1]\in u_2u_1^{-1}{\mth{U}}_{\alpha,2}$, hence $hu_1h^{-1}\in u_2{\mth{U}}_{\alpha,2}$.

We can iterate the reasoning, replacing $u_1$ by $hu_1h^{-1}$ and ${\mth{H}}_{\alpha,1}$ by ${\mth{H}}_{\alpha,2}$, then by ${\mth{H}}_{\alpha,3}$... and we finally obtain that $u_1$ is a conjugate of an element of ${\mth{U}}_{\alpha,d+1}u_2$ by an element of ${\mth{H}}_\alpha$. Since the group ${\mth{U}}_{\alpha,d+1}$ is trivial, the proposition is proved for $i=0$.

Now assume $i>0$. Let $h_0$ be an element of ${\mth{H}}_{\alpha,i}-{\mth{H}}_{\alpha,i+1}$; from the above discussion we deduce easily that the map $u\mapsto[h_0,u]$ from ${\mth{U}}_\alpha$ to ${\mth{U}}_{\alpha,i}$ is surjective and that the inverse image of ${\mth{U}}_{\alpha,i+1}$ by that map is ${\mth{U}}_{\alpha,1}$, hence the image of the set of elements of valuation $0$ is the set of elements of valuation $i$. Moreover, we have, for every $h\in{\mth{H}}_\alpha$:
\[Ad(h)[h_0,u]=[h_0,Ad(h)u].\]
Hence ${\mth{H}}_\alpha$ acts transitively on the set of elements of ${\mth{U}}_\alpha$ of valuation $i$ and the proposition is proved. $\Box$
\end{proof}

For every group $L$ and every pair of subgroups $A,B$ of $L$, we denote by $Z_A(B)$ the intersection with $A$ of the centralizer of $B$ in $L$.

\begin{cor}\label{priisec}
{\bf (PT4c')} implies {\bf (PT4c'')}).
\end{cor}

\begin{proof}
If $k$ is algebraically closed, ({\bf (PT4c'))} and ({\bf (PT4c'')} a re identical and there is nothing to prove; hence we can replace {\bf (PT4c')}) with ({\bf (PT4c'')} in the proof of propocition \ref{htrans2} in that case. Assume now we are in the general case. By proposition \ref{htrans2} applied to $\underline{\mth{H}}$ and the $\underline{\mth{U}}_{\alpha,i}$, $\alpha\in\Phi$ and $i\geq f_0(\alpha)$, we obtain that for every $\alpha,\beta,i,j$ such that $\underline{\mth{U}}_{\alpha,i}$ is nontrivial, we have:
\[p_{\alpha+\beta}([\underline{\mth{U}}_{\alpha,i},\underline{\mth{U}}_{\beta,j}])=p_{\alpha+\beta}([u,\underline{\mth{U}}_{\beta,j}]),\]
where $p_{\alpha+\beta}$ designates the projection on $\underline{\mth{U}}_{\alpha+\beta}$, and where $u$ is any element of $\underline{\mth{U}}_\alpha$ of valuation exactly $i$. By choosing $u$ inside  ${\mth{U}}_\alpha$, we obtain the result. $\Box$
\end{proof}

\begin{prop}\label{centru}
Let $\alpha$ be any element of $\Phi$ and let $u$ be an element of minimal valuation $v$ of ${\mth{U}}_\alpha$. Then $Z_{\mth{H}}(u)=Z_{\mth{H}}({\mth{U}}_\alpha)$.
\end{prop}

\begin{proof}
Let $h$ be any element of $Z_{\mth{H}}(u)$. Since ${\mth{H}}$ is abelian, $h$ also commutes with every element of the form $h'uh'^{-1}$, $h'\in{\mth{H}}$; by proposition \ref{htrans2}, $h$ then commutes with every element of valuation $v$ of ${\mth{U}}_\alpha$, hence with every element of ${\mth{U}}_\alpha$ by density. Hence $Z_{\mth{H}}(u)\subset Z_{\mth{H}}({\mth{U}}_\alpha)$; the other inclusion being obvious, the proposition is proved. $\Box$
\end{proof}

\begin{prop}\label{sumlie}
The sum of the $k$-vector spaces $Lie({\mth{H}}_{\alpha,1}))$, $\alpha\in\Phi$, is of dimension $dr$.
\end{prop}

\begin{proof}
If $k$ is algebraically closed, this is an immediate consequence of ${\bf (PT3)}$. In the general case, the sum of the $Lie({\mth{H}}_{\alpha,1})$, $\alpha\in\Phi$, generates as a $\overline{k}$-vectoor space the sum of the $Lie({\underline{\mth{H}}_{\alpha,1}})$, which is of dimension $dr$ , hence must be of that same fimension as a $k$-vector space. $\Box$
\end{proof}

Let $\alpha,\beta$ be two elements of $\Phi$; we say $\alpha$ and $\beta$ are {\em strongly orthogonal} if they are orthogonal and $\alpha+\beta$ is not a root.

Let $S$ be a subset of $\Phi$. We say $S$ {\em generates} $\Phi$ if every element of $\Phi$ is a linear combination with integer coefficients of the elements of $S$, or equivalently if the only closed root subsystem of $\Phi$ containing $S$ is $\Phi$ itself.

\begin{prop}\label{basprod}
Assume $\underline{\mth{G}}$ is of semisimple parahoric type, and let $S$ be any subset of $\Phi$ which generates it and whose cardinality is the rank of $\Phi$. Then ${\mth{H}}$ is the almost direct product of the ${\mth{H}}_\alpha$, $\alpha\in S$.
\end{prop}

\begin{proof}
We already know by lemma \ref{dpta} that ${\mth{T}}$ is the direct product of the ${\mth{T}}_\alpha$, $\alpha\in S$. It remains to prove that $R_u({\mth{H}})$ is the almost direct product of the ${\mth{H}}_{\alpha,1}$, $\alpha\in S$. We will in fact prove the equivalent result that $Lie(R_u({\mth{H}}))$ is the direct sum of the $Lie({\mth{H}}_{\alpha,1})$, $\alpha\in S$; by equality of dimensions it is enough to prove that the $Lie({\mth{H}}_{\alpha,1})$ generate $Lie(R_u({\mth{H}}))$.

First we consider the following Jacobi identity: let $\alpha,\beta$ be two elements of $\Phi$ such that $\alpha+\beta$ is also a root, and let $y_\alpha$ (resp. $y_\beta$, $y_{-\alpha-\beta}$) be an element of $Lie({\mth{U}}_{\alpha, f_0(\alpha)+1}$ (resp. $Lie({\mth{U}}_{\beta})$,  $Lie({\mth{U}}_{-\alpha-\beta,f_0(-\alpha-\beta)}$); we have:
\[[[y_\alpha,y_\beta],y_{-\alpha-\beta}]+[[y_\beta,y_{-\alpha-\beta}],y_\alpha]+[[y_{-\alpha-\beta},y_\alpha],y_\beta]=0.\]
The three terms of the left-hand side belong respectively to $Lie({\mth{H}}_{-\alpha-\beta,1})$, $Lie({\mth{H}}_{\alpha,1})$ and $Lie({\mth{H}}_{\beta,1})$, and, assuming that each one of these three subspaces of $Lie({\mth{H}})$ is generated by elements of the form of the corresponding term when $y_\alpha$ (resp. $y_\beta$, $y_{-\alpha-\beta}$) runs over $Lie({\mth{U}}_{\alpha,f_0(\alpha)+1})$ (resp. $Lie({\mth{U}}_\beta$, $Lie({\mth{U}}_{-\alpha-\beta})$), the sum of our three subspaces is generated by any two of them. By {\bf (PT4c)} and {\bf(PT4d)}, that condition is fulfilled as soon as $f_0(\alpha+\beta)=f_0(\alpha)+f_0(\beta)$; if this is not the case, then by proposition \ref{f0ab3} we must have $f_0(\alpha)=f_0(\alpha+\beta)+f_0(-\beta)$, and since $Lie({\mth{H}}_{-\gamma})=Lie({\mth{H}}_\gamma)$ for every $\gamma\in\Phi$, we obtain the same result. Hence if $\alpha,\beta$ are two elements of $S$ such that $\alpha+\beta\in\Phi$, the assertion of the proposition for $S$ is equivalent to the same assertion with $S$ replaced by the set $S'$ obtained by replacing either $\alpha$ or $\beta$ by $\alpha+\beta$. Note also that, again since $Lie({\mth{H}}_{-\gamma})=Lie({\mth{H}}_\gamma)$ for every $\gamma$, we can also replace any element of $S$ with its opposite. We will call such replacements {\em elementary replacements} in the rest of the proof. 

First assume $S$ is a basis of $\Phi$. By proposition \ref{sumlie},  the sum of the $Lie({\mth{H}}_{\alpha,1})$, $\alpha\in\Phi$, is $Lie({\mth{H}}_1)$; all we now have to prove is that for every $\alpha\in\Phi$, $Lie({\mth{H}}_{\alpha,1})$ is contained in the sum of the $Lie({\mth{H}}_{\beta,1})$, $\beta\in S$. Let then $\alpha$ be an element of $\Phi$ which is not in $S$; by eventually replacing some elements of $S$ with their opposites, we can assume that $\alpha$ is contained in the set of positive roots $\Phi^+_S$ associated to $S$. Since $\alpha$ does not belong to $S$, by \cite[\S 1, proposition 19]{bou}, there exists $\beta\in S$ such that $\alpha-\beta$ is a root, and by the above remarks, $Lie({\mth{H}}_{\alpha,1})$ is contained in $Lie({\mth{H}}_{\beta,1})+Lie({\mth{H}}_{\alpha-\beta,1})$. 

Let $l_S(\alpha)$ be the number of elements of $S$ (counted with multiplicities) $\alpha$ is the sum of; we will prove the desired assertion by induction on $l_S(\alpha)$. If $l_S(\alpha)=2$, then $\alpha-\beta\in S$ and there is nothing to prove. Assume now $l_S(\alpha)>2$; then $l_S(\alpha-\beta)=l_S(\alpha)-1$ and the result follows from the induction hypothesis.

%Now we consider the general case; we will prove that, by using elementary replacements, we can reduce ourselves to the previous case. Since $S$ generates $\Phi$ and since its cardinality is the rank of $\Phi$, it is also a subset of linearly independent roots of $\Phi$; by \cite[\S 1, proposition 22]{bou} applied to the subset $P$ of elements of $\Phi$ which are linear combinations with nonnegative integer coefficients of the elements of $S$, there exists then a set of positive roots $\Phi^+$ of $\Phi$ containing $S$. Let $\Delta$ be the set of simple roots of $\Phi^+$, and for every $\alpha\in S$, let $l_\Delta(\alpha)$ be defined as $l_S(\alpha)$, but relative to $\Delta$ this time. If $\sum_{\alpha\in S}l_\Delta(\alpha)=r$, then $S=\Delta$ and we are done; assume then $\sum_{\alpha\in S}l_\Delta(\alpha)>r$.

Now we prove the general case. Assume first $\Phi$ is irreducible. Let $\beta$ be an element of $\Delta$ which does not belong to $S$; we then have an expression of the form:
\[\beta=\sum_{\alpha\in S}\lambda_\alpha\varepsilon_\alpha\alpha,\]
with the $\lambda_\alpha$ being nonnegative integers and the $\varepsilon_\alpha$ being $\pm 1$; moreover, since $S$ is a generating subset of $\Phi^+$ which is not its basis and $\beta$ is one of its simple roots, the map $\alpha\in S\mapsto\varepsilon_\alpha$ is not constant.

We now prove the following lemma:

\begin{lemme}
There exist $\alpha,\alpha'\in S$ such that $\varepsilon_\alpha=-1$, $\varepsilon_{\alpha'}=1$ and $\alpha-\alpha'$ is a root.
\end{lemme}

\begin{proof}
Let $\alpha_1,\dots,\alpha_r$ be the elements of $S$, each of the $\alpha\in S$ showung up exactly $\lambda_\alpha$ times. By \cite[\S 1, proposition 19]{bou}, it is possible to choose the labelling in such a way that $\alpha_1+\dots+\alpha_i$ is a root for every $i\in\{1,\dots,r\}$, and since the map $\alpha\mapsto\varepsilon_\alpha$ is not constant, there exists $i\geq 2$ be such that $\varepsilon_1=\varepsilon_2=\dots=\varepsilon_{i-1}=-\varepsilon_i$; it is not hard to check that there is exactly one such $i$, and we have:
\[(\alpha_1+\dots+\alpha_{i-1})-\alpha_i\in\Phi.\]
Hence we must have $\alpha_j-\alpha_i\in\Phi$ for some $j<i$, and $\varepsilon_j=-\varepsilon_i$. The lemma is proved. $\Box$
\end{proof}

%We deduce immediately from the above equality the following one:
%\[\beta+\sum_{\alpha,\varepsilon_\alpha=-1}\lambda_\alpha\alpha=\sum_{\alpha',\varepsilon_{\alpha'}=1}\lambda_{\alpha'}\alpha'.\]
%Since both sides are obviously nonzero, there exists an $\alpha'$ such that we arein one of the following cases, $(.,.)$ being a $W$-invariant scalar product on $X^*({\mth{T}})$:
%\begin{itemize}
%\item there exists an $\alpha$ such that $(\alpha,\alpha')>0$;
%\item $(\beta,\alpha')>0$.
%\end{itemize}
%In the first case, we only have to apply \cite[\S 1, theorem 1]{bou} and we are done. In the second case, assume the first one does not hold as well. Then $\beta$ being a simple root, there exists an $\alpha'$ such that $\alpha'-\beta$ is a root, and we must then have $(\alpha,\alpha'-\beta)>0$ for some $\alpha$. 
%
%On the other hand, since $S$ generates $\Phi$ and $\Phi$ is irreducible, $S$ cannot be divided into two strongly orthogonal subsets, hence there must exist $\alpha,\alpha'\in S$ such that $\varepsilon_\alpha=-\varepsilon_{\alpha'}$ and we are in one of the following cases:
%\begin{itemize}
%\item $(\alpha,\alpha')>0$;
%\itemù $\alpha$ and $\alpha'$ are orthogonal but not strongly orthogonal.
%\end{itemize}
%
%In both cases, $\alpha-\alpha'$ is an element of $\Phi$: in the first case, this is a consequence of \cite[\S 1, Theorem 1]{bou}, in the second case, since $\alpha$ and $-\alpha'$ are also orthogonal but not strongly orthogonal, $\alpha-\alpha'$ must also be a root.

Now we go back to the proof of proposition \ref{basprod}. Let $\alpha,\alpha'$ be geiven by the previous lemma; assume for example $\alpha-\alpha'$ is positive, the other case being symmetrical. Then by elementary replacements we can replace $\alpha$ by $\alpha-\alpha'$ in $S$ (since we have $\alpha=\alpha'+(\alpha-\alpha')$, $(S-\{\alpha\})\cup\{\alpha-\alpha'\}$ still generates $\Phi$, and is still contained in our system of positive roots since we have assumed $\alpha-\alpha'$ is positive); we obviously have $l_\Delta(\alpha-\alpha')<l_\Delta(\alpha)$, and it is not hard to check that $S-\{\alpha\}\cup\{\alpha-\alpha'\}$ also generates $\Phi$.
By iterating the process, after a finite number of steps we reach an $S$ satisfying $\sum_{\alpha\in S}l_\Delta(\alpha)=r$, or in other words $S=\Delta$, as desired. If now $\Phi$ is reducible, we can do the same induction componentwise. $\Box$
\end{proof}

Note that as an obvious consequence of {(\bf PT3)}, when $\underline{\mth{G}}$ is of semisimple rank $2$, the almost direct product is in fact direct.

We deduce immediately from the proposition the following corollaries:

\begin{cor}
For every subset $S'$ of $S$ of cardinality $r'$, the product of the ${\mth{H}}_{\alpha,1}$, $\alpha\in S'$, is almost direct, and its dimension is $dr'$.
\end{cor}

Let $\Phi'$ be a root subsystem of $\Phi$. We say $\Phi'$ is a {\em Levi subsystem} if there exists a set of simple roots $\Delta$ of $\Phi$ such that $\Phi'$ is generated by some subset of $\Delta$.

\begin{cor}
Let $\Phi'$ be a Levi subsystem of $\Phi$ of dimension $r'$; the dimension of the product ${\mth{H}}_{\Phi,1}$ of the ${\mth{H}}_{\alpha,1}$, $\alpha\in\Phi'$, is $dr'$.
\end{cor}

\begin{proof}
Let $S'$ be a basis of $\Phi'$; we prove that the product of the ${\mth{H}}_{\alpha,1}$, $\alpha\in S'$, is ${\mth{H}}_{\Phi,1}$ the exact same way as proposition \ref{basprod}, and we conclude using the previous corollary. $\Box$
\end{proof}

\begin{cor}
The group ${\mth{H}}$ is the product of the center of ${\mth{G}}$ and the subgroups ${\mth{H}}_\alpha$, $\alpha\in S$.
\end{cor}

Note that the above product is not necessarily direct. For every $\alpha,\beta\in\Phi$ such that $\beta\not\in\Psi$, we deduce from {\bf (PT2)} and {\bf (PT4d'')} that the groups ${\mth{H}}_{\alpha,d}$ and ${\mth{U}}_\beta$ commute. Hence it may happen (for example when $\Psi$ is empty) that ${\mth{H}}_{\alpha,d}$ commutes with every ${\mth{U}}_\beta$, $\beta\in\Phi$, which then implies, since ${\mth{H}}$ is abelian, that ${\mth{H}}_{\alpha,d}$ is central in ${\mth{G}}$.

\section{Truncated valuation rings}

\subsection{Generalities}

Let $R$ be a commutative ring with unit, and let $d$ be a nonnegative integer. A {\em truncated valuation} (of depth $d$) on $R$ is a surjective application $v$ from $R$ to the set $V=\{0,1,\dots,d,+\infty\}$ satisfying the following conditions:
\begin{itemize}
\item for every $x\in R$, $v(x)=+\infty$ if and only if $x=0$;
\item for every $x,y\in R$, $v(xy)=v(x)+v(y)$ if $v(x)+v(y)\in V$, and $+\infty$ otherwise;
\item for every $x,y\in R$, $v(x+y)\geq Inf(v(x),v(y))$.
\end{itemize}

Note that the first two conditions imply $xy=0$ as soon as $v(x)+v(y)> d$. In particular, if $d\geq 1$, $R$ is not an integral domain.

As for local rings, a {\em uniformizer} of $R$ is an element of $R$ of valuation $1$.

The quotient $R/R_1$ will be called the {\em residual ring} (resp. the {\em residual field} when $R$ is local) of $R$. 

Since the product of two elements of $R$ of valuation $0$ is also of valuation $0$, $R/R_1$ is an integral domain, hence its characteristic $p$ is either $0$ or a prime number. Let $e$ be the valuation of $p$ in $R$; $e$ will be called the {\em absolute ramification index} of $R$. Note that a truncated valuation ring of characteristic $0$ is always of infinite absolute ramification index.

\begin{prop}\label{ideal}
The invertible elements of $R$ are all of valuation $0$, and the subsets $R_i=\{x\in R|v(x)\geq i\}$, $i\in\{0,\dots,d+1\}$ are ideals of $R$.

Conversely, assume
%that $R$ is noetherian,
that $R_1$ is the only prime ideal of $R$ and that the map $x\mapsto\varpi x$ from $R$ to $R_1$, where $\varpi$ is some given uniformizer of $R$, is surjective. Then $x\in R$ is invertible if and only if $v(x)=0$, and the $R_i$ are the only ideals of $R$.
\end{prop}

\begin{proof}
The first assertion is clear from the definitions. Now assume the hypotheses of the second one are satisfied. Let $x\in R$ be such that $v(x)=0$ and $x$ is not invertible. Then $xR+R_1$ is an ideal of $R$ which strictly contains $R_1$ but does not contain the identity element of $R$. By induction we can define an increasing set $I$ of ideals of $R$, all containing $R_1$ but not $R$, such that the union of every nonempty subset of $I$ is also in $I$. Applying Zorn's lemma, we obtain that this set admits a maximal element, which must be a maximal, hence also prime, ideal of $R$ strictly greater than $R_1$. On the other hand, such an ideal cannot exist by the hypotheses we have made; hence every element of valuation $0$ of $R$ is invertible.

Now let $R'$ be any ideal of $R$. Assume first $R'$ is not contained in $R_1$; it then contains an element of valuation $0$, hence invertible, which implies that we must have $R'=R=R_0$. Now assume $R'$ is contained in $R_1$, and let $i$ be the greatest integer such that $R'\subset R_i$; $R'$ then contains elements of valuation exactly $i$. Let $R''$ be the inverse image of $R'$ by the map $x\mapsto\varpi^ix$; it is not hard to see that $R''$ is also an ideal of $R$; on the other hand, we have the following lemma:

\begin{lemme}\label{rri}
The map $x\mapsto\varpi^ix$ from $R$ to $R_i$ is surjective.
\end{lemme}

\begin{proof}
We prove this lemma by induction on $i$, the case $i=1$ being one of the hypotheses of proposition \ref{ideal}. Let $y$ be an element of $R_i$; by that same hypothesis, there exists $z\in R$ such that $\varpi z=y$, and we then have $v(z)=v(y)-1\geq i-1$, hence $z\in R_{i-1}$; by induction hypothesis there exists $x\in R$ such that $\varpi^{i-1}x=z$, hence $\varpi^ix=y$, as desired. $\Box$
\end{proof}

By this lemma, since $R'$ contains an element $y$ of $R_i$ of valuation $i$, $R''$ contains an element $x$ of $R$ such that $\varpi^ix=y$, hence of valuation $0$; by the previous case, we must then have $R''=R$, hence by lemma \ref{rri} again, $R'=R_i$. The proposition is now proved. $\Box$
\end{proof}

From now on we assume that $R$ satisfies the hypotheses of the second assertion of proposition \ref{ideal}. This is in particular the case when $R$ is local.

 We have the following result:

\begin{prop}
Let $\varpi$ be a uniformizer of $R$, and let $S$ be a system of representatives of $R/R_1$ in $R$. Then any element $x\in R$ has a unique decomposition of the form:
\[x=\sum_{i=0}^ds_i\varpi^i\]
where the $s_i$ are elements of $S$.
\end{prop}

\begin{proof}
The proof is similar to the corresponding proof for valued rings (see \cite[II, proposition 5]{ser}). $\Box$
\end{proof}

\begin{prop}
Let $R$ be a valued ring and let the $R_i$, $i\in{\mth{N}}$, be defined as in proposition \ref{ideal}. For every $i\geq 1$, the valuation on $R$ canonically induces a truncated valuation on $R/R_i$.
\end{prop}

\begin{proof}
This is obvious from the definitions. $\Box$
\end{proof}

Conversely, we have the following result:

\begin{prop}\label{exfield}
Let $R$ be a local ring with a truncated valuation of depth $d$ and whose residual ring $k$ is a perfect field. There exists a henselian local field $F$ such that $R$ is isomorphic as a filtered ring to the quotient $\mathcal{O}/\mathfrak{p}^{d+1}$, where $\mathcal{O}$ is the ring of integers of $F$ and $\mathfrak{p}$ is the maximal ideal of $\mathcal{O}$.
\end{prop}

\begin{proof}
First assume that $R$ is either of valuation $0$ or of infinite absolute ramification index; the rings $R$ and $k$ then have the same characteristic, and we deduce from \cite[II.4, prop. 6 and 8]{ser} that there exists a subring $R^0$ of $R$ isomorphic to $k$ and such that the restriction to $R^0$ of the canonical projection $R\rightarrow k$ is an isomorphism; $R$ is then a finite-dimensional $k$-algebra. If $\varpi$ is any uniformizer of $R$, consider the morphism $k[[X]]\rightarrow R$ sending $k$ to $R^0$ and $X$ to $\varpi$; it is obviously surjective and its kernel is the ideal $X^{d+1}k[[X]]$. Hence we can take $F=k((X))$ in this case.

Now assume that $R$ is of positive characteristic and of finite absolute ramification index. By \cite[II.5, theorem 3]{ser}, there exists a unique local field $F_0$ which is complete (hence henselian), of charateristic $0$, absolutely unramified, and whose residual field is $k$. Let $\mathcal{O}_0$ be its ring of integers; by \cite[II.5, proposition 10]{ser}, there exists a ring morphism $\phi$ from $\mathcal{O}_0$ to $R$ such that, if $\pi$ (resp. $\pi'$)) is the canonical projection from $\mathcal{O}_0$ (resp. $R$) to $k$, $\pi'\circ\phi=\pi$. Let $R_\#$ be the image of $\phi$ in $R$; $R_\#$ is then an absolutely unramified local ring with a truncated valuation, and $R$ is of the form $R_\#[\varpi]$, where $\varpi$ is any uniformizer of $R$.

Let $\varpi$ be such a uniformizer, and let $P$ be a polynomial of degree $e$ in $F_0[X]$ such that (with a slight abuse of notation) $P(\varpi)=0$; $P$ is then an Eisenstein polynomial, hence irreducible. Let $F$ be a splitting field for $P$ over $F_0$, and let $\mathcal{O}$ be its ring of integers; the morphism $\phi$ then extends to a surjective morphism from $\mathcal{O}$ to $R$, and its kernel must then be $\mathfrak{p}^{d+1}$, where $\mathfrak{p}$ is the maximal ideal of $\mathcal{O}$, which proves the proposition. $\Box$
\end{proof}

Note that the field $F$ is not unique. In the case of a ring of infinite absolute ramification index, we can take as $F$ any henselian local field of absolute ramification index strictly greater than $d$ and of residual field $k$. Even in the case of finite ramification index, and even if we restrict ourselves to complete fields, there may be several possibilities for $F$ (for example, with $R={\mth{Z}}_2[\sqrt{2}]/2\sqrt{2}{\mth{Z}}_2[\sqrt{2}]$, we can of course take $F={\mth{Q}}_2[\sqrt{2}]$, but $F={\mth{Q}}_2[\sqrt{3}]$ works as well.)

\begin{cor}
The truncated valuation ring $R$ is a henselian ring.
\end{cor}

\begin{proof}
This is an immediate consequence of the henselianity of $F$. (Note that it can also be proved directly; the proof is similar to the proof of Hensel's lemma for complete local fields.) $\Box$
\end{proof}

\begin{cor}\label{tralgg}
The truncated valuation ring $R$ is, as an additive group, the group of $k$-points of a finite-dimensional unipotent abelian algebraic group of dimension $d+1$ defined over $k$.
\end{cor}

\begin{proof}
When $R$ is of infinite absolute ramification index, it is isomorphic as an additive group to the direct sum of $d+1$ copies of $k$ and the result follows immediately. Assume now $R$ is of finite absolute ramification index, and let $F_0$ be the complete absolutely unramified local field defined in the proof of proposition \ref{exfield}; by \cite[Satz 10]{witt}, the ring of integers  $\mathcal{O}_0$ of $F_0$ is isomorphic to a ring of Witt vectors over $k$ (see \cite{witt} or \cite{ser} for a detailed definition of a ring of Witt vectors). which is itself the group of $k$-points of an infinite-dimensional abelian unipotent algebraic group defined over $k$. Hence the ring $R_\#$ defined as in the proof of proposition \ref{exfield}, being a finite-dimensional quotient of $\mathcal{O}_0$, is itself the group of $k$-points of a finite-dimensional abelian unipotent algebraic group defined over $k$. The corollary is then proved when $R=R_\#$, or in other words when $R$ is absolutely unramified. In the general case, let $R_b$ be the quotient of $R_\#$ by the ideal generated by one of its elements of maximal finite valuation; if $e$ is the absolute ramification index of $R$, we have, as additive groups, if $\varpi$ is any uniformizer of $R$:
\[R=R_\#\oplus\varpi R_\#\oplus\dots\oplus\varpi^iR_\#\oplus\varpi^{i+1}R_b\oplus\dots\oplus\varpi^{e-1}R_b,\]
where $i$ is such that $e\,dim(R_\#)+i=d+1$. The corollary follows immediately. $\Box$
\end{proof}

We will generalize the above corollary later (proposition \ref{kkd}).

\begin{lemme}
Let $R$ be a truncated valuation ring of residual field $k$. Assume every element of $k$ has two opposite square roots in $k$. Then every element of $R$ of valuation $0$ has two opposite square roots in $R$.
\end{lemme}

\begin{proof}
Since the square roots are distinct, $k$ is not of characcteristic $2$; the result is then simply a particular case of Hensel's lemma. $\Box$
\end{proof}

\subsection{The group ${\mth{U}}_\alpha$ as a truncated valuation ring}

Let $\alpha$ be an element of $\Phi$; set ${\mth{U}}={\mth{U}}_\alpha$. In this section, we treat the group ${\mth{U}}$ as an additive group on which ${\mth{H}}$ acts by adjunction, and write its composition law as $+$. For the rest of the section we will assume $f_0(\alpha)=0$ to simplify the notations; it is possible by proposition \ref{fazero}. (In the general case, an element of ${\mth{U}}$ will have two distinct valuations: its valuation as an element of ${\mth{U}}_\alpha$, and its valuation as an element of the ring $R_\alpha$ defined below, which differ by $f_0(\alpha)$. We thus will have to be careful about which valuation we are talking about.)

We fix arbitrarily an element $u_1$ of ${\mth{U}}$ of valuation $0$; by proposition \ref{htrans2} and proposition \ref{centru}, for every $u\in{\mth{U}}$ of valuation $0$, there exists an element $h_u$ of ${\mth{H}}$, unique up to the centralizer of ${\mth{U}}$ in ${\mth{H}}$, such that $Ad(h_u)u_1=u$. In particular, $Ad(h_u)$ is entirely determined by $u$. We can define the following multiplication on ${\mth{U}}$: let $u,u'\in{\mth{U}}$. If $v(u)=0$, we set:
\[u.u'=Ad(h_u)u'.\]
If now $v(u)>0$, then we must have $v(u+u_1)=0$, and we can thus set:
\[u.u'=(u+u_1).u'-u'.\]
Of course the multiplication depends on the choice of $u_1$.

\begin{lemme}\label{adhuu}
Let $u,u'$ be two elements of valuation $0$ of ${\mth{U}}$. We have $Ad(h_{u.u'})=Ad(h_u)Ad(h_{u'})$.
\end{lemme}

\begin{proof}
We have $u.u'=Ad(h_{u.u'})(u_1)$; on the other hand, we also have $u.u'=Ad(h_u)(u')=Ad(h_uh_{u'})(u_1)$. The lemma then comes from the unicity of $Ad(h_{u.u'})$. $\Box$
\end{proof}

\begin{cor}\label{adcom}
We have $Ad(h_{u.u'})=Ad(h_{u'.u})$.
\end{cor}

\begin{proof}
This is an immediate consequence of the previous lemma and the commutativity of ${\mth{H}}$. $\Box$
\end{proof}

\begin{lemme}\label{vuhu}
For every $u\in{\mth{U}}$ such that $v(u)>0$ and every choice of $h_{u+u_1}\in{\mth{H}}$, we have $v_\alpha(h_{u+u_1})=v(u)$.
\end{lemme}

\begin{proof}
We have $[h_{u+u_1},u_1]=u+u_1-u_1=u$, which proves the assertion since $v(u_1)=0$. $\Box$
\end{proof}

\begin{prop}
The multiplication operation defined in this section makes the group $({\mth{U}},+)$ into a commutative local ring $R=R_\alpha$, and the valuation $v$ is a truncated valuation on that ring.
\end{prop}

\begin{proof}
First we check that the multiplication is commutative. We immediately deduce from corollary \ref{adcom} that $u.u'=u'.u$ when $u$ and $u'$ are both of valuation $0$. Now assume $v(u)=0$ and $v(u')>0$; we have:
\[u.u'=Ad(h_u)(u')=Ad(h_u)(u'+u_1)-Ad(h_u)(u_1)\]
\[=u.(u'+u_1)-u=(u'+u_1).u-u=u'.u.\]
Finally assume $u$ and $u'$ are both of positive valuation; we have:
\[u.u'=(u+u_1).u'-u'=u'.(u+u_1)-u'=(u'+u_1)(u+u_1)-u-u_1-u',\]
and by a similar reasoning:
\[u'.u=(u+u_1)(u'+u_1)-u'-u_1-u.\]
Since $u+u_1$ and $u'+u_1$ commute, $u$ and $u'$ also commute. Hence the multiplication is commutative.

We deduce in particular from this that for every two elements $u,u'$ of ${\mth{U}}$ such that $v(u')>0$, we have:
\[u.u'=u.(u'+u_1)-u.\]

We now check that $u_1$ is the identity element for the multiplication. Since $Ad(h_{u_1})(u_1)=u_1$, by unicity $Ad(h_{u_1})$ must be the identity, hence $u_1.u=u$ for every $u\in{\mth{U}}$; we also obtain $u.u_1=u$ by commutativity.

We now check the distributivity. Since we already know that the multiplication is commutative we only have to check its left distributivity; let then $u,u',u''$ be three elements of ${\mth{U}}$, we prove now that $u.(u'+u'')=u.u'+u.u''$. When $v(u)=0$, that simply comes from the distributivity of $Ad(h_u)$. Assume now $v(u)>0$; we have, using the fact that $v(u+u_1)=0$:
\[u.(u'+u'')=(u+u_1).(u'+u'')-u'-u''=(u+u_1).u'-u'+(u+u_1).u''-u''=u.u'+u.u'',\]
which proves the assertion.

We now check the associativity. Let $u,u',u''$ be three elements of ${\mth{U}}$; we prove that $(u.u').u''=u.(u'.u'')$. First assume $v(u)=v(u')=0$; we have $(u.u').u''=Ad(h_{u.u'})(u'')$ and $u.(u'.u'')=Ad(h_u)(Ad(h_{u'})(u''))$, and lemma \ref{adhuu} then implies the result. Next assume $v(u)>0$ and $v(u')=0$; we have:
\[(u.u').u''=((u+u_1).u'-u').u''=((u+u_1).u').u''-u'.u''\]
\[=(u+u_1).(u'.u'')-u'.u''=u.(u'.u'').\]
Finally assume $v(u')>0$; we then have $v(u.u')>0$, and, using the previous case:
\[(u.u').u''=(u.(u'+u_1)-u).u''=(u.(u'+u_1)).u''-u.u''\]
\[=u.((u'+u_1).u'')-u.u''=u.((u'+u_1).u''-u'')=u.(u'.u'').\]

We thus have proved that $({\mth{U}},+,.)$ is a commutative ring. Moreover, we deduce immediately from the definition of the multiplication that for every $i\geq 0$, the subgroup ${\mth{U}}_i$ of the $u\in{\mth{U}}$ such that $v(u)\geq i$ is an ideal of that ring. To prove that ${\mth{U}}$ is local, we only have to check that every element of ${\mth{U}}$ of valuation $0$ is invertible. Let $u$ be such an element, and set $u'=Ad(h_u^{-1})(u_1)$; we then have $u'.u=u_1$, hence $u$ is invertible.

We finally have to check that $v$ is a truncated valuation on the ring ${\mth{U}}$. Since $v$ is already a truncated valuation on the root subgroup ${\mth{U}}$ of ${\mth{G}}$, we only have to check that for every $u,u'\in{\mth{U}}$ such that $u.u'\neq 0$, we have $v(u.u')=v(u)+v(u')$.

For every $h\in{\mth{H}}$, let $v_\alpha(h)$ be the unique integer such that $[h,{\mth{U}}]={\mth{U}}_{v_\alpha(h)}$; obviously, $v_\alpha(h)$ depends only on the class of $h$ modulo the centralizer of ${\mth{U}}$ in ${\mth{H}}$, hence $v_\alpha(h_u)$ depends only on $u\in{\mth{U}}$. 

If $v(u)=0$, the equality $v(u.u')=v(u)+v(u')$ comes immediately from the definition of the multiplication; assume now $v(u)>0$. We then have $u.u'=(u+u_1).u'-u'=[h_{u+u_1},u']$ by lemma \ref{vuhu}, hence when $u.u'\neq 0$, $v(u.u')=v_\alpha(h_{u+u_1})+v(u')=v(u)+v(u')$, as required.  $\Box$
\end{proof}

Let $R_\alpha$ be the commutative ring defined in the previous proposition. By proposition \ref{exfield},  we know that there exists a nonarchimedean local field $F$ such that $R_\alpha$ is isomorphic to the quotient ring $R_F=\mathcal{O}_F/\mathfrak{p}_F^\delta$, where $\mathcal{O}_F$ is the ring of integers of $F$, $\mathfrak{p}$ is the maximal ideal of $\mathcal{O}$ and $\delta$ is either $d+1$ or $d$, depending whether $\alpha$ lies in $\Psi$ or not. Moreover, by corollary \ref{tralgg}, $R_F$ is the group of $k$-points of a unipotent algebraic group of dimension $d+1$ defined over $k$. We also define the ring $R'_\alpha$ as the quotient of $R_\alpha$ by the ideal generated by $\varpi^{d-1}$; note that when $\alpha\not\in\Psi$, we have $R'_\alpha=R_\alpha$.

Note that the above construction of $R_\alpha$ (resp. $R'_\alpha$) works in particular when $k$ is algebraically closed. We will thus denote, for every $\alpha\in\underline{\Phi}$, by $\underline{R}_\alpha$ (resp. $\underline{R'}_\alpha$) the ring associated to the subgroup $\underline{\mth{U}}_\alpha$ (resp. the quotient of such a ring by the ideal of its elements of valuation $d-1$ or greater). 

Now we prove the following result:

\begin{prop}
There exists a ring isomorphism between $R_F$ and $R=R_\alpha$ which is also an isomorphism of $k$-algebraic groups.
\end{prop}

\begin{proof}
Let $\phi$ be a ring isomorphism between $R_F$ and $R$; since $\phi$ sends invertible (resp. non-invertible) eleents of $R_F$ on invertible (resp. non-invertible) elements of $R$ we must have $\phi(R_{F_1})=R_1 $, and $\phi$ thus induces a field isomorphism $\phi_0$ between $R_F/R_{F,1}$ and $R/R_1$, and we can always choose $\phi$ in such a way that $\phi_0$ is $k$-algebraic. Then $\phi$ is a bijective morphism of $k$-varieties, hence also of $k$-algebraic groups. According to \cite[5.3.3]{spr}, to prove that $\phi$ is an isomorphism of algebraic groups, we now only have to prove that the tangent map $d\phi$ between the Lie algebras $Lie(R_F)$ and $Lie(R)$ is surjective.

For every $i$, we can check with the help of an easy induction that $\phi(R_{F,i})=R_i$. Let then $\phi_i$ be the lift to $R_{F,i}/R_{F,i+1}$ of the restriction of $\phi$ to $R_{F,i}$; $\phi_i$ is then a bijective morphism of algebraic groups between $R_{F,i}/R_{F,i+1}$ and $R_i/R_{i+1}$.

Now for every element $u$ of ${\mth{U}}_\alpha$ of valuation $1$, consider the element $h_{u+u_1}$ of ${\mth{H}}_\alpha$; it belongs to ${\mth{H}}_{\alpha,1}$, and we know from lemma \ref{huhu} and the remarks folloowing it that for every integer $i\leq d-1$, the composition map:
\[R_{F,i}/R_{F,i+1}\rightarrow R_i/R_{i+1}\rightarrow R_{i+1}/R_{i+2}\rightarrow R_{F,i+1}/R_{F,i+2},\]
where the first and third maps are respectively $\phi_i$ and $\phi_{i+1}^{-1}$ and the second one is induced by $u'\mapsto[h_{u+u_1},u']$, is of the form $x\mapsto c_ix$, where $c_i$ is a nonzero constant which does not depend on the choice of $u$; hence if $\phi_i$ is an isomorphism of $k$-algebraic groups, then so is $\phi_{i+1}$. Since we already know that $\phi_0$ is an isomorphism of $k$-algebraic groups, by induction every $\phi_i$, $i>0$, is also an isomorphism of $k$-algebraic groups, hence for every $i$, $d\phi_i$ is surjective. Hence $d\phi$ must be surjective as well, which proves the proposition. $\Box$
\end{proof}

We now use this proposition to prove the following one:

\begin{prop}\label{kkd}
Let $G$ be a reductive algebraic group defined over a nonarchimedian local field $F$ of residual field $k$ and split over the maximal unramified extension $F_{nr}$ of $F$;, let $K$ be a $Gal(F_{nr}/F)$-stable parahoric subgroup of $G$ and let $K'$ be a $Gal(F_{nr}/F)$-stable normal subgroup of $K$ contained in the pro-unipotent radical $K^0$ of $K$. Then the quotient $K/K'$ can be viewed as an algabraic group defined over $k$.
\end{prop}

\begin{proof}
It is already well-know, that $K/K^0$ is a reductive group defined over $k$; we thus only have to consider the group $K^0/K'$. SInce this group is normalized by any maximal $k$-torus $\underline{\mth{T}}$ of $\underline{\mth{G}}$, it is a product of subquotients of root subgroups relative to $\underline{\mth{T}}$ and of a subquotient of the maximal torus $T$ of $G$ corresponding to ${\mth{T}}$; the first kind of subgroups are isomorphic to rings of the form described by proposition \ref{tralgg}, hence have a $k$-structure by that proposition, and it is not hard to check that the latter group is a subgroup of a product of groups of unit elements of such rings, hence also has a $k$-structure. Hence $K/K'$ also has a $k$-structure and the proposition is proved. $\Box$
\end{proof}

We now check that up to isomorphism, the ring $R$ does not depend on the choice of $u_1$.

\begin{prop}
Let $u'_1$ be another element of ${\mth{U}}$ of valuation $0$, and let $R_*$ be the ring constructed over the group ${\mth{U}}$ the same way as $R=({\mth{U}},+,.)$, but taking $u'_1$ instead of $u_1$ as its unit element. Then the map $\phi=Ad(h_{u'_1})$ is a ring isomorphism from $R$ to $R_*$.
\end{prop}

\begin{proof}
We already know that $\phi$ is a group automorphism of ${\mth{U}}$; moreover, by definition of $h_{u'_1}$, $\phi(u_1)=u'_1$. We thus only have to check that for every $u,u'\in{\mth{U}}$, $\phi(u.u')=\phi(u)*\phi(u')$, where $*$ designates the multiplication in $R_*'$.

Assume first $v(u)=0$; we have:
\[\phi(u.u')=Ad(h_{u'_1})(Ad(h_u)u')=Ad(h_{u'_1}h_uh_{u'_1}^{-1})(\phi(u'))\]
\[=Ad(h_{u'_1}h_uh_{u'_1}^{-1})(u'_1)*\phi(u')=Ad(h_{u'_1})(u)*\phi(u')=\phi(u)*\phi(u').\]
Assume now $v(u)>0$; we have:
\[\phi(u.u')=\phi((u+u_1).u'-u')=(\phi(u)+u'_1)*\phi(u')-\phi(u')=\phi(u)*\phi(u'),\]
which proves the result. $\Box$
\end{proof}

\section{The case of rank $1$}

Now we investigate what happens when $\Phi$ is of rank $1$. We will consider the two possible cases separately: the case $\Psi=\Phi$ (which implies ${\mth{G}}$ is not solvable) and the case $\Psi=\emptyset$ (which implies ${\mth{G}}$ is solvable). We will also assume in this section that $\underline{\mth{G}}$ is of semisimple parahoric type; the general case of rank $1$ groups will be covered by the main theorem in section $6$.

\subsection{The nonsolvable case}

In this subsection, we assume ${\mth{G}}$ is not solvable. Let $\pm\alpha$ be the elements of $\Phi=\Psi$; we can also assume $f_0(\alpha)=f_0(-\alpha)=0$.

Since $\Psi=\Phi$, $\Psi$ has a nontrivial Weyl group $W$. Moreover, since the group $W$ normalizes ${\mth{T}}$, it also normalizes ${\mth{H}}$, and since we have:
\[W=N_{\mth{G}}({\mth{T}})/Z_{\mth{G}}({\mth{T}})=N_{\mth{G}}({\mth{T}})/{\mth{H}},\]
the group $W$ embeds canonically into $N_{\mth{G}}({\mth{H}})/{\mth{H}}$. Moreover, every element of ${\mth{G}}$ which normalizes ${\mth{H}}$ also normalizes its only maximal torus ${\mth{T}}$, hence the above embedding is an isomorphism.

In the sequel, we will simply identify $W$ with the group $N_{\mth{G}}({\mth{H}})/{\mth{H}}$. Let $w$ be its unique nontrivial element; we then have $w={\mth{H}}n=n{\mth{H}}$, with $n$ being any element of $w$. We also have $w(\alpha)=-\alpha$, which implies that if $n$ is any representant of $w$, $n{\mth{U}}_\alpha n^{-1}={\mth{U}}_{-\alpha}$, or in other words that $Ad(n$) induces an isomorphism between $R_\alpha$ and $R_{-\alpha}$ (resp. $R'_\alpha$ and $R'_{-\alpha}$), which of course depends on $n$.

We start by proving the following results. Remember that we are assuming $p\neq 2$.

\begin{lemme}\label{unord2}
There exists a unique element of order $2$ in ${\mth{H}}$, which belongs to ${\mth{T}}$.
\end{lemme}

\begin{proof}
By proposition \ref{hsq}, every element of ${\mth{H}}={\mth{H}}_\alpha$, and in particular the identity, has exactly two square roots in ${\mth{H}}$. Since one of them is the identity itself, ${\mth{H}}$ contains exactly one element $h$ of order $2$. Moreover, since $p\neq 2$, $h$ is semisimple, hence belongs to ${\mth{T}}$. $\Box$
\end{proof}

\begin{lemme}\label{centg}
The center of ${\mth{G}}$ is trivial if the semisiple group $\underline{\mth{G}}/R_u(\underline{\mth{G}})$ is adjoint, and of order $2$ if $\underline{\mth{G}}/R_u(\underline{\mth{G}})$ is simply-connected.
\end{lemme}

\begin{proof}
Since the corresponding result for ${\mth{G}}/R_u({\mth{G}})$ is well-known (it comes for example from the fact that $G$  is isomorphic to either $SL_2$ or $PGL_2$, see \cite[proposition 7.2.4]{spr}), all we have to do to prove the first assertion is to check that $R_u({\mth{G}})$ does not contain any nontrivial central element. We have $R_u({\mth{G}})={\mth{U}}_{\alpha,1}{\mth{H}}_{\alpha,1}{\mth{U}}_{-\alpha,1}$; let then $g=uhu'$ be a central element of ${\mth{G}}$ contained in $R_u({\mth{G}})$. By conjugating $g$ by a suitable element of ${\mth{T}}$ we see that $u$ and $u'$ must be trivial, hence $g\in{\mth{H}}$. On the other hand, we see with the help of {\bf (PT4d'')} and of lemma \ref{huhu} that every nontrivial element of ${\mth{H}}_{\alpha,1}$ acts nontrivially on the root subgroups; hence $g$ must be trivial.

Assume now $\underline{\mth{G}}/R_u(\underline{\mth{G}})$ is simply-connected (hence isomorphic to $SL_2$). Since the above reasoning works just as well in this case, the center of ${\mth{G}}$ is canonically isomorphic to a subgroup of the center of ${\mth{G}}/R_u({\mth{G}})$; to finish the proof of the lemma, we thus only have to check that it contains a nontrivial element. Let then $g$ be the unique element of order $2$ of ${\mth{T}}$; it acts on the rings $R_\alpha$ and $R_{-\alpha}$ by multiplication by either $1$ or $-1$, and since a simple computation in $SL_2$ shows that its action on the projections of these two rings on ${\mth{G}}/R_u({\mth{G}})$ is trivial, its action on the rings themselves must be trivial too. Hence $g$ commutes with ${\mth{U}}_\alpha$ and ${\mth{U}}_{-\alpha}$, and also with ${\mth{H}}$ by definition of ${\mth{H}}$; it is then central in ${\mth{G}}$ and the lemma is proved. $\Box$
\end{proof}

\begin{prop}\label{isoo}
Assume $k$ is any infinite field. Let ${\mth{L}},{\mth{L}}'$ be the groups of $k$-points of two connected algebraic groups defined over $k$. Assume there exist a dense open subset $\Omega$ (resp. $\Omega'$) of ${\mth{L}}$ (resp. ${\mth{L}}'$), a dense open subset $\Omega''$ of ${\mth{L}}\times{\mth{L}}$ such that the map $(h,l)\mapsto hl$ from $\Omega''$ to ${\mth{L}}$ contains $\Omega$ in its image, and an isomorphism of algebraic varieties $\phi: \Omega\rightarrow\Omega'$ such that for every $h,l\in\Omega$ such that $(h,l)\in\Omega''$ and $hl\in\Omega$, $\phi(hl)=\phi(h)\phi(l)$. Then $\phi$ can be extended into an isomorphism of algebraic groups between ${\mth{L}}$ and ${\mth{L}}'$.
\end{prop}

\begin{proof}

We first prove the following lemmas:

\begin{lemme}
Let $l$ be any element of $\Omega$. The subset $\Omega_l$ of the elements $h\in\Omega$ such that $(h^{-1},hl)\in\Omega''$ is dense in ${\mth{L}}$.
\end{lemme}

Set $V''={\mth{L}}\times {\mth{L}}-\Omega''$; $V''$ is a closed subvariety of ${\mth{L}}\times {\mth{L}}$, hence its intersection with the closed connected subvariety $\{(h^{-1},hl)|h\in {\mth{L}}\}$ is also closed; by definition of $\Omega''$,  it is strictly contained in that subvariety, hence of dimension strictly smaller. The result follows. $\Box$

\begin{lemme}
Let $h$ be an element of ${\mth{L}}$ which does not belong to $\Omega$; there exist $h_1,h_2\in\Omega$ such that $h=h_1h_2$.
\end{lemme}

\begin{proof}
Since $\Omega$ is open and dense, $V={\mth{L}}-\Omega$ is a closed subset of ${\mth{L}}$ of strictly smaller dimension. Hence for every $h\in {\mth{L}}$, the set $V\cup hV^{-1}$, where $V^{-1}$  is the set of inverses of elements of $V$, is strictly contained in ${\mth{L}}$. Let $h_1$ be an element of ${\mth{L}}$ which does not belong to $V\cup hV^{-1}$; we then have $h_1\in\Omega$, and since $h_1\not\in hV^{-1}$, $h_2=h_1^{-1}h$ also belongs to $\Omega$, which proves the lemma. $\Box$
\end{proof}

Now we prove proposition \ref{isoo}. For every $h=h_1h_2\in {\mth{L}}$, with $h_1$ and $h_2$ being elements of $\Omega$, set $\phi(h)=\phi(h_1)\phi(h_2)$; we first have to prove that our definition is consistent.

First assume $h\in\Omega$. By definition of $\phi$, we have $\phi(h)=\phi(h_1)\phi((h_2)$ as soon as $(h_1,h_2)\in\Omega''$. Since $h\in\Omega$, the set of $h_1\in\Omega$ such that, setting $h_2=h_1^{-1}h$,  $(h_1,h_2)\in\Omega''$ is nonempty, hence dense since it is open. The map:
\[h_1\mapsto\phi(h_1)\phi(h_1^{-1}h)\]
is then of constant value $\phi(h)$ on an open dense subset of ${\mth{L}}$. Since it is obviously continuous on its whole domain of definition, that map must then take that same constant value everywhere. Hence the identity $\phi(h)=\phi(h_1)\phi(h_2)$ holds regardless of whether $(h_1,h_2)$ lies in $\Omega''$ or not.

Now assume $h$ does not belong to $\Omega$. Write $h=h_1h_2=h_3h_4$, and let $l$ be an element of $\Omega$ such that $lh_1$ and $lh_3$ also belong to $\Omega$; such an $l$ exists because $\Omega\cap\Omega h_1^{-1}\cap\Omega h_3^{-1}$ is open and dense in ${\mth{L}}$. We then have, by the previous case:
\[\phi(lh)=\phi(lh_1)\phi(h_2)=\phi(lh_3)\phi(h_4);\]
\[\phi(lh_1)=\phi(l)\phi(h_1);\]
\[\phi(lh_3)=\phi(l)\phi(h_3).\]
We deduce from these three equalities that $\phi(h_1)\phi(h_2)=\phi(h_3)\phi(h_4)$, as desired. Hence the definition of $\phi$ is consistent.

Now let $h,l$ be any two elements of ${\mth{L}}$; we will check that $\phi(hl)=\phi(h)\phi(l)$. If $h,l\in\Omega$ and $(h,l)\in\Omega''$ this is simply the definition of $\phi$, and we have previously seen that it also holds when $(h,l)\not\in\Omega''$. Assume now $h\in V$ and $l\in\Omega$, and write $h=h_1h_2$, with $h_1,h_2\in\Omega$ being such that $h_2l\in\Omega$; we then have:
\[\phi(hl)=\phi(h_1)\phi(h_2l)=\phi(h_1)\phi(h_2)\phi(l)=\phi(h)\phi(l).\]
The case $h\in\Omega$ and $l\in V$ is symmetrical, and the case $h,l\in V$ is treated similarly, using that last case.

We now have to prove that $\phi$ is an isomorphism of abstract groups. Let $h$ be an element of its kernel; for every $l\in\Omega$ such that $hl\in\Omega$, we then have $\phi(l)=\phi(hl)$; since $\phi$ is injective on $\Omega$, we must then have $h=1$. Let now $h'$ be an element of ${\mth{L}}'$ which does not belong to the image of $\phi$; then $h'\Omega'$ is a nonempty open subset of ${\mth{L}}$ disjoint from $\Omega'$, which is impossible since $\Omega'$ is dense. Hence $\phi$ is an isomorphism.

It only remains to check that $\phi$ is an isomorphism of algebraic varieties. For every $h\in{\mth{L}}$, the application $l\mapsto\phi(h)\phi(h^{-1}l)$ from $h\Omega$ to $\phi(h)\Omega'$ is an isomorphism of algebraic varieties which coincides with $\phi$ on $\Omega\cap h\Omega$. Moreover, since $\Omega$ is open and dense, it is easy to check that ${\mth{L}}$ is the union of a finite number of subsets of the form $h\Omega$; by glueing the corresponding applications we obtain that $\phi$ is algebraic on the whole group and the proposition is proved. $\Box$
\end{proof}

Before proving the main result of this section, we need some more preliminary results. Remember that a pseudo-$k$-parabolic subgroup (or simply pseudo-parabolic when there is no ambiguity) of a nonreductive group $\underline{\mth{G}}$ is the group generated by the centralizer $\underline{\mth{H}}$ of some maximal $k$-split torus $\underline{\mth{T}}$ and the root subgroups $\underline{\mth{U}}_\alpha$, $\alpha\in X$, where $X$ contains some set of positive weights of $\underline{\mth{G}}$ relative to $\underline{\mth{T}}$. It is not hard to check that this definition of a pseudo-parabolic subgoroup is equivalent to the definition of \cite[15.1]{spr}. When $\underline{\mth{H}}$ is a Cartan subgroup of $\underline{\mth{G}}$ and $X$ is a set of positive weights itself we will speak of a pseudo-Borel subgroup. We will also say ${\mth{B}}$ is a pseudo-Borel subgroup of ${\mth{G}}$ if it is the group of $k$-points of some pseudo-Borel subgroup of $\underline{\mth{G}}$ defined over $k$.

Let ${\mth{B}}$ be the pseudo-Borel subgroup of ${\mth{G}}$ generated by ${\mth{H}}$ and ${\mth{U}}_\alpha$, and let $w$ be the nontrivial element of the Weyl group $W$ of ${\mth{G}}$ with respect to ${\mth{T}}$; the Bruhat decomposition of ${\mth{G}}$ (lifted from the Bruhat decoposition of its reductive quotient) can be written the following way:
\[{\mth{G}}={\mth{B}}{\mth{U}}_{-\alpha,1}\sqcup{\mth{B}}w{\mth{B}}.\]
Here we are using the fact that ${\mth{B}}{\mth{U}}_{-\alpha,1}w{\mth{B}}{\mth{U}}_{-\alpha,1}={\mth{B}}w{\mth{B}}$, whose proof is straightforward. A proof of a similar result for $GL_n$ can be found in \cite{hill} for example.

Since ${\mth{B}}{\mth{U}}_{-\alpha,1}$ is a closed subgroup of ${\mth{G}}$ of strictly smaller dimension, the subset ${\mth{B}}w{\mth{B}}$ is open and dense in ${\mth{G}}$. According to proposition \ref{isoo} we then only have to find a bijection $\phi$ between ${\mth{B}}w{\mth{B}}$ and the corresponding subset of ${\mth{G}}_0$ such that for every $g,g'\in{\mth{B}}w{\mth{B}}$ such that $gg'\in{\mth{B}}w{\mth{B}}$, $\phi(gg')=\phi(g)\phi(g')$. In this case, the subset $\Omega''$ of ${\mth{B}}w{\mth{B}}\times{\mth{B}}w{\mth{B}}$ we use for the proposition is simply ${\mth{B}}w{\mth{B}}\times{\mth{B}}w{\mth{B}}$ itself, which satisfies the required condition by the following lemma:

\begin{lemme}
We have ${\mth{B}}w{\mth{B}}\subset {\mth{B}}w{\mth{B}}w{\mth{B}}$.
\end{lemme}

\begin{proof}
Since both members of the above equality are unions of cosets modulo $R_u({\mth{G}})$, the lemma is an immediate consequence of the corresponding result for reductive groups (see for example \cite[lemma 8.3.7]{spr}). $\Box$
\end{proof}

Since ${\mth{G}}$ is of semisimple parahoric type and not solvable, we deduce from {\bf (PT4d'')} that the group generated by ${\mth{U}}_\alpha$ and ${\mth{U}}_{-\alpha}$ contains ${\mth{H}}$, hence is the whole group ${\mth{G}}$; we then have the following proposition:

\begin{prop}\label{nunu}
Let $n$ be any representative of $w$ in $N_{\mth{G}}({\mth{H}})$; there exist $u\in{\mth{U}}_\alpha$ and $u'\in{\mth{U}}_{-\alpha}$ such that $n=uu'u=u'uu'$; moreover, we have $u'=n^{-1}un=nun^{-1}$.
\end{prop}

\begin{proof}
Let $u'$ be any element of ${\mth{U}}_{-\alpha}$ of valuation $0$. Since it does not belong to ${\mth{B}}{\mth{U}}_{-\alpha,1}$, it must belong to ${\mth{B}}w{\mth{B}}={\mth{B}}w{\mth{U}}_\alpha$. Moreover, since ${\mth{B}}={\mth{U}}_\alpha{\mth{H}}$, we have ${\mth{B}}w{\mth{U}}_{\alpha}={\mth{U}}_\alpha{\mth{H}}w{\mth{U}}_\alpha={\mth{U}}_\alpha w{\mth{U}}_\alpha$ (remember that $w$, as an element of the Weyl group, is a class mod ${\mth{H}}$, hence the expression ${\mth{U}}_\alpha w{\mth{U}}_\alpha$ makes sense); there exists then a representative $n_{u'}$ of $w$ such that $u'$ belongs to ${\mth{U}}_\alpha n_{u'}{\mth{U}}_\alpha$. Write $n_{u'}=u_1u'u_2$, $u_1,u_2\in{\mth{U}}_\alpha$; to prove the first equality of the lemma for $n_{u'}$, we only have to check that $u_1=u_2$.

We now prove the following lemmas:

\begin{lemme}\label{nhn}
For every $h\in{\mth{H}}$, we have $nhn^{-1}=h^{-1}$.
\end{lemme}

\begin{proof}
Let $x,y\in R_\alpha$ be such that $Ad(h)u_{\alpha,1}=u_{\alpha,x}$ and $Ad(n)u_{\alpha,1}=u_{-\alpha,y}$; we have:
\[Ad(nhn^{-1})u_{-\alpha,1}=Ad(nh)u_{\alpha,y^{-1}}=Ad(n)u_{\alpha,xy^{-1}}=u_{-\alpha,x},\]
hence $nhn^{-1}$ and $h^{-1}$ act the same way on ${\mth{U}}_{-\alpha}$. Hence the map $h\mapsto nhn^{-1}h$ is a continuous group morphism from ${\mth{H}}$ to the center of ${\mth{G}}$, which by lemma \ref{centg} is a finite group; that morphism must then be trivial and the lemma is proved. $\Box$
\end{proof}

For every $x\in R_\alpha$, we will denote by $u_{\alpha,x}$ the element $x$ viewed as an element of ${\mth{U}}_\alpha$.

\begin{lemme}
We have $N_{\mth{G}}({\mth{H}})\cap Z_{\mth{G}}(u_{\alpha,1})=Z({\mth{G}}$).
\end{lemme}

\begin{proof}
We will prove that $N_{\mth{G}}({\mth{H}})\cap Z_{\mth{G}}(u_{\alpha,1})\subset Z({\mth{G}})$, the other inclusion being trivial. Let $h$ be an element of $N_{\mth{G}}({\mth{H}})\cap Z_{\mth{G}}(u_{\alpha,1})$. Since $h$ commutes with a nontrivial element of ${\mth{U}}_\alpha$ it must be contained in ${\mth{H}}$, which implies that it commutes with ${\mth{H}}$ since ${\mth{H}}$ is abelian. It also commutes with every element of the form $hu_{\alpha,1}h^{-1}=u_{\alpha,x}$, with $x$ being of valuation $0$, hence with every element of ${\mth{U}}_\alpha$ by density. Finally, by lemma \ref{nhn}, we have $h=nh'^{-1}n^{-1}$, which commutes with $n{\mth{U}}_\alpha n^{-1}={\mth{U}}_{-\alpha}$, hence $h$ commutes with ${\mth{U}}_{-\alpha}$ as well. Since the subgroups ${\mth{H}}$, ${\mth{U}}_\alpha$ and ${\mth{U}}_{-\alpha}$ generate ${\mth{G}}$, the lemma is proved. $\Box$
\end{proof}

\begin{lemme}\label{hnh}
For every $n\in w$, there exists $h=h_n\in{\mth{H}}$ such that $n=hn_{u'}h^{-1}$.
\end{lemme}

\begin{proof}
Let $x_0$ (resp. $x$) be the element of $R_{-\alpha}$ such that $Ad(n_{u'})u_{\alpha,1}=u_{-\alpha,x_0}$ (resp. $Ad(n)u_{\alpha,1}=u_{-\alpha,x}$); $x_0$ and $x$ are both of valuation $0$. Let $y$ be a square root of $x_0^{-1}x$ in $R_{-\alpha}$, and let $h$ be an element of ${\mth{H}}_{-\alpha}$ such that $Ad(h)u_{\alpha,1}=u_{\alpha,y}$; we have:
\[Ad(hn_{u'}h^{-1})u_{\alpha,1}=Ad(hn_{u'})u_{\alpha,y^{-1}}=Ad(h)u_{-\alpha,x_0y}=u_{-\alpha,x_0y^2}=u_{-\alpha,x}.\]
Hence we have $n=zhn_{u'}h^{-1}$, with $z$ being central in ${\mth{G}}$. If $z$ is trivial we are done; if $z$ is nontrivial, then by lemma \ref{centg} it is of order $2$. Let $y'=-y$ be the other square root of $x_0^{-1}x$; we have:
\[u_{\alpha,y'}=Ad(hz')u_{\alpha,1},\]
where $z'$ is an element of ${\mth{H}}$ such that $Ad(z')u_{\alpha,1}=u_{\alpha,-1}$. On the other hand, since $z'\neq z$, by lemma \ref{unord2}, $z'$ cannot be of order $2$; we must then have $z'^2=z$. By lemma \ref{nhn}, we have $z'n_{u'}=n_{u'}z'^{-1}$; we thus obtain:
\[n=zhn_{u'}h^{-1}=z'^2hn_{u'}h^{-1}=(z'h)n_{u'}(z'h)^{-1}.\]
The lemma is now proved. $\Box$
\end{proof}

Note that $h_n$ is not unique. When ${\mth{G}}$ has a nontrivial center, $h_n$ is defined up to an element of that center. When ${\mth{G}}$ has a trivial center, a simple computation in $PGL_2$ shows that the unique element of order $2$ of ${\mth{H}}_\alpha$, which acts on both ${\mth{U}}_\alpha$ and ${\mth{U}}_{-\alpha}$ by $x\mapsto -x$, commutes with $n_{u'}$,  hence $h_n$ is defined up to that element. On the other hand, it is easy to check that two different candidates for $h_n$ always differ by an element of order $2$ of ${\mth{H}}$, hence by lemma \ref{unord2} $h_n$ can take exactly two different values.

We deduce from lemma \ref{hnh} that for $h=h_n$, we have $n=(hu_1h^{-1})(hu'h^{-1})(hu_2h^{-1})$. This is true in particular for $n=n_{u'}^{-1}$.

\begin{lemme}
The element $h=h_{n_{u'}^{-1}}$ can be chosen in such a way that the automorphism $Ad(h_{n_{u'}^{-1}})$ acts on both ${\mth{U}}_\alpha$ and ${\mth{U}}_{-\alpha}$ by $x\mapsto -x$.
\end{lemme}

\begin{proof}
Let $h$ be one of the two possible choices for $h_{n_{u'}^{-1}}$.
%We have $n_{u'}^{-1}=hn_{u'}h^{-1}$ and $n_{u'}=(n_{u'}^{-1})^{-1}=hn_{u'}^{-1}h^{-1}$; hence $h^2$ and $n_{u'}$ commute. We deduce from this that the action of $h^2$ by adjunction is the same on both root subgroups, hence that it is either trivial or $x\mapsto -x$.
When ${\mth{G}}$ has a trivial center, we can easily check that $n_{u'}=n_{u'}^{-1}$, hence $h$ and $n_{u'}$ commute; by the remark following lemma \ref{hnh}, it is then possible to choose $h$ in such a way that $Ad(h)$ acts on the root subgroups by $x\mapsto -x$. Now assume that ${\mth{G}}$ has a nontrivial center. Since $h^2$ and $n_{u'}$ commute, by the remark following lemma \ref{hnh}, $h^2$ belongs to the center of ${\mth{G}}$, hence $Ad(h)$ must act on the root subgroups either trivially or by $x\mapsto -x$. Assume that the action of $Ad(h)$ on the root subgroups is trivial. Then $h$ commutes with ${\mth{H}}_\alpha$, ${\mth{U}}_\alpha$ and ${\mth{U}}_{-\alpha}$, hence is central in ${\mth{G}}$. On the other hand, in the proof of lemma \ref{hnh}, we have already exhibited an element $z'$ of ${\mth{H}}$ which does not commute with $n_{u'}$, and this leads to a contradiction. Hence $Ad(h)$ acts on the root subgroups as $x\mapsto -x$ and the lemma is proved. $\Box$
\end{proof}

From now on until the end of the proof of proposition \ref{nunu}, we assume $h=h_{n_{u'}^{-1}}$ actually satisfies the condition of the above lemma. We then obtain:
\[n_{u'}=(n_{u'}^{-1})^{-1}=(hu_1h^{-1})(hu'h^{-1})(hu_2h^{-1})^{-1}\]
\[=(u_1^{-1}u'^{-1}u_2^{-1})^{-1}=u_2u'u_1,\]
hence $u'=u_2^{-1}n_{u'}u_1^{-1}$. On the other hand, since ${\mth{U}}_\alpha$ and $n_{u'}{\mth{U}}_\alpha n_{u'}^{-1}={\mth{U}}_{-\alpha}$ have a trivial intersection, any element of ${\mth{U}}_\alpha n_{u'}{\mth{U}}_\alpha$ can be written in a unique way in the form $un_{u'}v$; hence $u_1=u_2$. For every $n\in w$, we then have:
\[n=(h_nu_1h_n^{-1})(h_nu'h_n^{-1})(h_nu_1h_n^{-1}),\]
with $h_n$ defined as in lemma \ref{hnh}, which completes the proof of the first equality of proposition \ref{nunu}.

To prove the next two, we observe that we have:
\[u'^{-1}u^{-1}n=u,\]
which can be rewritten as:
\[u'^{-1}n(n^{-1}u^{-1}n)=u,\]
hence:
\[n=u'u(n^{-1}un)\]
By the same reasoning as above (with $\alpha$ and $-\alpha$ switched), we must then have $n^{-1}un=u'$, hence $n=u'uu'$. The equality $u'=nun^{-1}$ is obtained in a similar way.  $\Box$
\end{proof}

Now we prove the main result of the paper in the particular case of a group $\underline{\mth{G}}$ of semisimple rank $1$, not solvable and of semisimple parahoric type. Actually, we prove the corresponding result for the group ${\mth{G}}$ of $k$-points of $\underline{\mth{G}}$, assuming $k$ satisfies the same hypotheses as in the previous section; the assertion for $\underline{\mth{G}}$ is simply the particular case $k=\overline{k}$.
%since, by lemma \ref{kinf}, $k$ is then infinite, by \cite[proposition 18.3]{bor}, ${\mth{G}}$ is dense in $\underline{\mth{G}}$, hence the result for ${\mth{G}}$ implies the result for $\underline{\mth{G}}$ as well (and conversely).

\begin{prop}\label{a1nr}
Assume $\Phi$ is of rank $1$ and ${\mth{G}}$ is of semisimple parahoric type and not solvable. There exists a nonarchimedean henselian local field $F$ of residual field $k$ and an algebraic group $\underline{G}$ of type $A_1$ defined and split over $F$ such that ${\mth{G}}$ is isomorphic to the quotient ${\mth{G}}_0$ of a maximal parahoric subgroup of the group $G$ of $F$-points of $\underline{G}$ by its $d$-th congruence subgroup.
\end{prop}

For the definition of the congruence subgroups of the parahoric subgroups of $G$, see subsection $2.2$.

\begin{proof}
By proposition \ref{exfield}, we know that there exists a nonarchimedean henselian local field $F$ such that the quotient of its ring of integers $\mathcal{O}_F$ by $\mathfrak{p}_F^{d+1}$, where $\mathfrak{p}_F$ is the maximal ideal of $\mathcal{O}_F$, is isomorphic to $R_\alpha$ and $R_{-\alpha}$. Let $\underline{G}$ be an algebraic group defined and split over $F$ whose root datum is $(X^*({\mth{T}}),\Phi,X_*({\mth{T}}),\Phi^\vee)$; such a group exists by theorem 10.1.1 of \cite{spr}. Let $G$ be the group of $F$-points of $\underline{G}$, let $K$ be any maximal parahoric subgroup of $G$ and set ${\mth{G}}_0=K/K^d$, where $K^d$ is the $d$-th congruence subgroup of $K$ (see eection $2.2$); ${\mth{G}}_0$ is then the group of $k$-points of an algebraic group of parahoric type defined and split over $k$ which has the same root datum as ${\mth{G}}$; in particular, since $\Psi=\Phi$ is nontrivial, ${\mth{G}}_0$ is not solvable. Let ${\mth{B}}_0$ be any pseudo-Borel subgroup of ${\mth{G}}_0$; consider a pseudo-Levi decomposition ${\mth{B}}_0={\mth{H}}_0{\mth{U}}_{0,\alpha_0}$ of ${\mth{B}}_0$, where $\alpha_0$ is the root of ${\mth{G}}_0/{\mth{T}}_0$ corresponding to ${\mth{B}}_0$, ${\mth{T}}_0$ being the unique maximal torus of ${\mth{H}}_0$. Let $\phi$ be a group isomorphism between ${\mth{U}}_\alpha$ and ${\mth{U}}_{0,\alpha_0}$ coming from a ring isomorphism between $R_\alpha$ and $\mathcal{O}_F/\mathfrak{p}_F^{d+1}$; $\phi$ can be extended to an isomorphism from ${\mth{B}}$ to ${\mth{B}}_0$ by considering the isomorphism between ${\mth{T}}$ and ${\mth{T}}_0$ and setting, for every $h\in{\mth{H}}_{\alpha,1}$, $\phi(h)=h_0$ where $h_0$ is the unique element of the group ${\mth{H}}_{\alpha,0}$ corresponding to ${\mth{H}}_\alpha$ in ${\mth{G}}_0$ such that $Ad(h_0)\phi(u)=\phi(Ad(h)u)$ for every $u\in{\mth{U}}$. Let $w_0$ be the nontrivial element of the Weyl group of ${\mth{G}}_0$ relative to ${\mth{T}}_0$; we can extend $\phi$ to a bijection between ${\mth{B}}w{\mth{B}}$ and ${\mth{B}}_0w_0{\mth{B}_0}$ by choosing an arbitrary element $n_0$ (resp. $n$) of $w_0$ (resp. $w$) and setting for every element $g=unb$ of ${\mth{U}}w{\mth{B}}={\mth{B}}w{\mth{B}}$, $\phi(unb)=\phi(u)n_0\phi(b)$.

First we remark that if we write $n=uu^-u$, with $u\in{\mth{U}}_\alpha$ and $u^-\in{\mth{U}}_{-\alpha}$, a necessary condition for $\phi$ to satisfy $\phi(gg')=\phi(g)\phi(g')$ for every $g,g'$ is that $n_0$ satisfies $n_0=\phi(u)\phi(u')\phi(u)$, hence $\phi(u)^{-1}n_0\phi(u)^{-1}\in{\mth{U}}_{0,-\alpha_0}$. Hence we must take as $n_0$ an element of $w_0$ satisfying that condition; the existence of such an element follows from proposition \ref{nunu},

Now we prove the following lemma:

\begin{lemme}
For every $h,h'\in{\mth{H}}$ and every $g=unb\in{\mth{U}}w{\mth{B}}$, we have $\phi(hgh')=\phi(h)\phi(g)\phi(h')$.
\end{lemme}

\begin{proof}
As an immediate consequence of the definition of $\phi$ we can assume $h'=1$. We have:
\[\phi(hg)=\phi(hunb)=\phi(huh^{-1})n_0\phi(n^{-1}hn)\phi(b).\]
By lemma \ref{nhn}, we have $n^{-1}hn=h^{-1}$ and $n_0^{-1}\phi(h)^{-1}n_0=\phi(h)$, and we obtain: 
\[\phi(n^{-1}hn)=\phi(h^{-1})=\phi(h)^{-1}=n_0^{-1}\phi(h)n_0,\]
hence:
\[\phi(hg)=\phi(hunb)=\phi(huh^{-1})n_0n_0^{-1}\phi(h)n_0\phi(b)\]
\[=\phi(huh^{-1})\phi(h)n_0\phi(b)=\phi(h)\phi(u)n_0\phi(b)=\phi(h)\phi(unb)=\phi(h)\phi(g).\] $\Box$
\end{proof}

We now have to check that $\phi(gg')=\phi(g)\phi(g')$ for every $g,g'\in{\mth{B}}w{\mth{B}}$ such that $gg'\in{\mth{B}}w{\mth{B}}$. We see immediately from the definitions that we can assume $g\in w{\mth{B}}$ and $g'\in{\mth{U}}n^{-1}$; moreover, writing $g=nhu$ and $g'=u'n^{-1}$, with $h\in{\mth{H}}$ and $u,u'\in{\mth{U}}_\alpha$, the condition $gg'\in{\mth{B}}w{\mth{B}}$ is equivalent to $v(uu')=0$. Moreover, we have $\phi(g)=\phi(nhn^{-1})\phi(n)\phi(u)$, $\phi(g')=\phi(u')\phi(n^{-1})$ and $\phi(uu')=\phi(u)\phi(u')$, which proves that, by replacing $g$ by $(nhn^{-1})^{-1}gu'=nuu'$ and $g'$ by $u'^{-1}g'=n^{-1}$, we can assume $h=u'=1$; we thus only have to prove that for every $u\in{\mth{U}}_\alpha$ of valuation $0$, $\phi(nu)\phi(n^{-1})=\phi(nun^{-1})$.

First we assume $u$ is such that $n$ is of the form $n=uu^-u=u^-uu^-$, with $u^-\in{\mth{U}}_{-\alpha}$; such a $u$ exists by proposition \ref{nunu}. We then have, by definition of $\phi$:
\[\phi(u^-)=\phi(u)^{-1}\phi(n)\phi(u)^{-1},\]
hence:
\[\phi(n)=\phi(u)\phi(u^-)\phi(u)=\phi(u^-)\phi(u)\phi(u^-),\]
hence:
\[\phi(nu)\phi(n^{-1})=\phi(n)\phi(u)\phi(n^{-1})\]
\[=\phi(u^-)\phi(u)\phi(u^-)\phi(u)\phi(n)^{-1}\]
\[=\phi(u^-)\phi(n)\phi(n)^{-1}=\phi(u^-)=\phi(nun^{-1}).\]
which proves the assertion.

Now we prove the general case. Let $u'$ be the element of ${\mth{U}}_\alpha$ satisfying $n=u'u^-u'$ for some $u^-\in{\mth{U}}_{-\alpha}$ and let $h\in{\mth{H}}_\alpha$ be such that $hu'h^{-1}=u$; we have, using the fact that $nhn^{-1}h^{-1}\in{\mth{H}}$:
\[\phi(nu)=\phi(nhu'h^{-1})=\phi((nhn^{-1}h^{-1})(hnh^{-1})(hu'h^{-1}))=\phi(nhn^{-1}h^{-1})\phi((hnh^{-1})(huh^{-1}));\]
\[\phi(n^{-1})=\phi((hn^{-1}h^{-1})(hnh^{-1}n^{-1}))=\phi(hn^{-1}h^{-1})\phi(hnh^{-1}n^{-1});\]
\[\phi(nun^{-1})=\phi(nhu'h^{-1}n^{-1})=\phi((nhn^{-1}h^{-1})(hnh^{-1})(hu'h^{-1})(hn^{-1}h^{-1})(hnh^{-1}n^{-1}))\]
\[=\phi(nhn^{-1}h^{-1})\phi((hnh^{-1})(hu'h^{-1})(hn^{-1}h^{-1}))\phi(hnh^{-1}n^{-1}).\]
Since $hnh^{-1}=(hu'h^{-1})(hu^-h^{-1})(hu'h^{-1})$, we are reduced to the previous case. $\Box$
\end{proof}

Now we prove a relation between some elements of ${\mth{G}}$ which will be useful in the sequel. The group $\underline{G}$ is semisimple, hence isomorphic to either $SL_2$ or $PGL_2$, and can thus be written as (a quotient of) a group of $2\times 2$ matrices, and it is then also the case for ${\mth{G}}_0$. We can then assume that ${\mth{H}}_0$ is the Cartan subgroup of diagonal elements of ${\mth{G}}_0$, and that ${\mth{U}}_{0,\alpha_0}$ (resp. the image of ${\mth{U}}_{-\alpha}$ in ${\mth{G}}_0$) is the subgroup of upper (resp. lower) triangular elements of ${\mth{G}}$ with both diagonal terms being $1$. (Note that we cannot just speak about unipotent upper (lower) triangular elements because an upper (lower) triangular element of ${\mth{G}}$ may be unipotent while having diagonal terms different from $1$.)

For every $u\in{\mth{U}}_{\alpha}$ (resp. ${\mth{U}}_{-\alpha}$), if $x$ is the corresponding element in $R_\alpha$ (resp. $R_{-\alpha}$), we have set $u_{\alpha,x}=u$ (resp. $u_{-\alpha,x}=u$). Similarly, for every $x\in 1+\mathfrak{p}_F$, if $h$ is the inverse image of $\left(\begin{array}{cc}x&0\\0&x^{-1}\end{array}\right)$ (resp. $\left(\begin{array}{cc}x^{-1}&0\\0&x\end{array}\right)$) in ${\mth{H}}_\alpha$ we write $h=h_{\alpha^\vee,x}$ (resp. $h=h_{-\alpha^\vee,x}$). Note that for every $x$, we have $h_{-\alpha^\vee,x}=h_{\alpha^\vee,x^{-1}}$.

Note that, although the map $x\in R_\alpha^*\mapsto h_{-\alpha^\vee,x}$ is always surjective, it is not necessarily injective: when ${\mth{G}}$ has a trivial center, its kernel is $\{\pm 1\}$. 

\begin{prop}\label{ualphac}
For suitable choices of the unit elements in ${\mth{U}}_\alpha$ and ${\mth{U}}_{-\alpha}$, we have, for every $\lambda\in R_\alpha$, $\lambda\neq -1$:
\[u_{-\alpha,\lambda}u_{\alpha,1}=u_{\alpha,\frac 1{1+\lambda}}h_{\alpha^\vee,\frac 1{1+\lambda}}u_{-\alpha,\frac\lambda{1+\lambda}}.\]
\end{prop}

\begin{proof}
We can see by an easy matrix computation that the result holds for the corresponding elements of ${\mth{G}}_0$. By the isomorphism obtained in proposition \ref{a1nr}, it holds in ${\mth{G}}$ as well. $\Box$
\end{proof}

\subsection{The solvable case}

Assume now ${\mth{G}}$ is solvable, and let $\pm\alpha$ be the elements of $\Phi$ once again. We can put a ring structure on ${\mth{U}}_\alpha$ the same way as  in section $3$; if $R_\alpha$ is that ring,, for every $x\in R_\alpha$, we will write $u_{\alpha,x}$ instead of $x$ for clarity; same for ${\mth{U}}_{-\alpha}$, which defines a ring $R_{-\alpha}$. Moreover, ${\mth{H}}_\alpha$ acts on both ${\mth{U}}_\alpha$ and ${\mth{U}}_{-\alpha}$ by adjunction, which allows us to define the elements $h_{\alpha^\vee,x}$, $x\in {\mth{H}}_\alpha^*$, the following way: let $x$ be an element of $\mathcal{O}_F$, and write it $x=x_tx_u$, with $x_t\in k^*$  and $x_u\in 1+\mathfrak{p}_F$. Then $h_{\alpha^\vee,x}=h_{\alpha^\vee,x_t}h_{\alpha^\vee,x_u}$, where $h_{\alpha^\vee,x_t}=\alpha^\vee(x_t)$ and $h_{\alpha^\vee,x_u}$ is a representant of the only element of ${\mth{H}}_{\alpha,1}/{\mth{H}}_{\alpha,d}$ acting on $U_\alpha$ by adjunction by $x^2$. (it is easy to check that it is possible to choose these representants in such a way that the map $x_u\mapsto h_{\alpha^\vee,x_u}$ is a group morphism whose kernel is $1+\mathfrak{p}^{d+1}$). We then have $Ad(h_{\alpha^\vee,x})(u_{\alpha,y})=u_{\alpha,x^2y}$ and $Ad(h_{\alpha^\vee,x})(u_{-\alpha,y})=u_{-\alpha,\frac y{x^2}}$ for every $x,y$; note that, as in the solvable case, when ${\mth{G}}$ has a trivial center, the map $x\in R_\alpha^*\mapsto h_{\alpha^\vee,x}$ is not injective.

We first have to prove that $R_\alpha$ and $R_{-\alpha}$ are isomorphic; we then establish a relation between elements of ${\mth{H}}$, ${\mth{U}}_\alpha$ and ${\mth{U}}_{-\alpha}$ which is analogous to the one of proposition \ref{ualphac} for nonsolvable groups (relation (\ref{scom})), and we are then ready to prove the main result (proposition \ref{soca}).

Since ${\mth{G}}$ is solvable, we have the following Iwahori decompositions:
\[{\mth{G}}={\mth{U}}_\alpha{\mth{H}}{\mth{U}}_{-\alpha}={\mth{U}}_{-\alpha}{\mth{H}}{\mth{U}}_\alpha.\]
We deduce from that decomposition that there exist functions $a$, $b$, $c$ from $R_{-\alpha}$ to respectively $R_\alpha$, ${\mth{H}}_\alpha$, $R_{-\alpha}$ such that we have, for every $\lambda\in R_{-\alpha}$:
\begin{equation}\label{abc}u_{-\alpha,\lambda}u_{\alpha,1}=u_{\alpha,a(\lambda)}h_{\alpha^\vee,b(\lambda)}u_{-\alpha,c(\lambda)}.\end{equation}
Note that we have already remarked that when ${\mth{G}}$ has a trivial center, there are two possible choices for $b(\lambda)$. On the other hand, we deduce from {\bf (PT4c'')} that $h_{\alpha^\vee,b(\lambda)}\in{\mth{H}}_{\alpha,1}$; hence $b(\lambda)$ is a square root of some element of $1+R_{\alpha,1}$, and we set it to be the one belonging to $1+R_{\alpha,1}$.

Similarly, there exist functions $a'$, $b'$, $c'$ from $R_\alpha$ to respectively $R_{-\alpha}$, ${\mth{H}}_\alpha$, $R_\alpha$ such that for every $\lambda'\in R_\alpha$, we have:
\begin{equation}\label{abc1}u_{\alpha,\lambda'}u_{-\alpha,1}=u_{-\alpha,a'(\lambda')}h_{\alpha^\vee,b'(\lambda')}u_{\alpha,c'(\lambda')}.\end{equation}
These six functions are morphisms of $k$-varieties. (Note that since $1+R_{\alpha,1}$ is a unipotent group, it is easy to check (details are left to the reader) that its automorphism $x\mapsto\sqrt{x}$ is $k$-algebraic, hence the presence of a square root in the definitions of the function $b$ and $b'$ (see above) does not cause any algebraicity issues.)

Let $\phi$ be the group isomorphism between $R_\alpha^*$ and $R_{-\alpha}^*$ such that for every $h\in{\mth{H}}$, $\phi(Ad(h)u_{\alpha,1})=Ad(h)^{-1}u_{-\alpha,1}$. Choose a uniformizer $\varpi$ of $R_{-\alpha}$ and a uniformizer $\varpi'$ of $R_\alpha$; $\phi$ extends to a semigroup isomorphism between $R_\alpha$ and $R_{-\alpha}$ (viewed as multiplicative semigroups), which we will also denote by $\phi$, by setting $\phi(\varpi')=\varpi$. To prove that $\phi$ is a ring isomorphism between $R_\alpha$ and $R_{-\alpha}$, we only have to check that it preserves addition; this is not true for every choice of $\varpi$ and $\varpi'$, and we will have to find suitable ones. For the moment, though, we choose them arbitrarily.

First we remark that since $-1$ is the only element of order $2$ in $R_\alpha^*$ (resp. $R_{-\alpha}^*$), we must have $\phi(-1)=-1$. By multiplicativity we obtain that for every $x\in R_\alpha$, $\phi(-x)=-\phi(x)$.

Let $\lambda$ be an element of $R_{-\alpha}$; assume first $v(\lambda)=0$. Let $\mu$ be a square root of $-\lambda$ in $R_{-\alpha}$; by conjugating both sides of (\ref{abc}) by $h_{\alpha^\vee,\mu}$, we obtain:
\[u_{-\alpha,-1}u_{\alpha,\phi^{-1}(-\lambda)}=u_{\alpha,\phi^{-1}(-\lambda)a(\lambda)}h_{\alpha^\vee,b(\lambda)}u_{-\alpha,-\frac{c(\lambda)}\lambda},\]
from which we deduce, using (\ref{abc1}), that for every $\lambda$ such that $v(\lambda)=0$, we have:
\[\lambda a'(\phi^{-1}(\lambda))=c(\lambda),\]
\[b'(\phi^{-1}(\lambda))=b(\lambda)^{-1};\]
\[c'(\phi^{-1}(\lambda))=\lambda a(\lambda).\]

When $v(\lambda)>0$, the left-hand sides of these three equalities depend on the choice of $\phi$; we will find a suitable $\phi$ later on. We can already remark that all three of them are trivially true when $\lambda=0$.

Now we prove a few unicity lemmas:

\begin{lemme}
Let $\lambda,\lambda'$ be two elements of $R_{-\alpha}$. Assume that $c(\lambda)=c(\lambda')$. Then $\lambda=\lambda '$.
\end{lemme}

\begin{proof}
We have, using (\ref{abc}):
\[u_{-\alpha,\lambda-\lambda'}=(u_{\alpha,1}u_{-\alpha,\lambda})^{-1}(u_{\alpha,1}u_{-\alpha,\lambda'})\]
\[=u_{\alpha,a(\lambda)}h_{\alpha^\vee,b(\lambda)}u_{-\alpha,c(\lambda)}u_{-\alpha,-c(\lambda')}h_{\alpha^\vee,-b(\lambda')^{-1}}u_{\alpha,-a(\lambda')}\]
\[=u_{\alpha,a(\lambda)}h_{\alpha^\vee,b(\lambda)b(\lambda')^{-1}}u_{\alpha,-a(\lambda')}.\]
The left-hand side belongs to ${\mth{U}}_{-\alpha}$ and the right-hand side to ${\mth{B}}={\mth{H}}{\mth{U}}_\alpha$; since the intersection of these two groups is trivial, the assertion follows. $\Box$
\end{proof}

\begin{lemme}\label{unicb}
Let $\lambda,\lambda'$ be two elements of $R_{-\alpha}$. Assume that $b(\lambda)=b(\lambda')$. Then $\lambda=\lambda '$.
\end{lemme}

\begin{proof}
We have, using (\ref{abc}) again:
\[[u_{\alpha,-1},u_{-\alpha,\lambda}]=u_{\alpha,a(\lambda)-1}h_{\alpha^\vee,b(\lambda)}u_{-\alpha,c(\lambda)-\lambda}\]
and a similar expression with $\lambda'$ instead of $\lambda$. The result then follows immediately from the following lemma:

\begin{lemme}
The application $\psi$ which associates to every $u\in{\mth{U}}_{-\alpha}$ the projection on ${\mth{H}}$ of $[u_{\alpha,-1},u]$ is a bijective morphism of $k$-varieties between ${\mth{U}}_{-\alpha}$ and ${\mth{H}}_{\alpha,1}$.
\end{lemme}

\begin{proof}
From {\bf (PT4d'')}, we know that every element of ${\mth{H}}_{\alpha,d}$ is the projection on ${\mth{H}}$ of $[u,u']$ for some $u\in{\mth{U}}_{-\alpha,d-1}$ and some $u'\in{\mth{U}}_{\alpha}$ of valuation $0$. By conjugating the above commutator by an appropriate element of ${\mth{H}}$ we can even assume $u'=u_{\alpha,-1}$. The restriction of $\psi$ to ${\mth{U}}_{-\alpha,d-1}$ is then a surjective morphism of $k$-varieties.

On the other hand, ${\mth{H}}_{\alpha,d}$ lies in the center of ${\mth{G}}$, and a simple computation shows that the restriction $\psi'$ of $\psi$ to ${\mth{U}}_{-\alpha,d-1}$ is a morphism of algebraic groups between ${\mth{U}}_{-\alpha,d-1}$ and ${\mth{H}}_{\alpha,d}$. Since both groups are isomorphic to ${\mth{G}}_a$ and since $\psi'$ is onto ${\mth{H}}_{\alpha,d}$, it must be bijective.

For every $i\leq d$, by replacing ${\mth{G}}$ by its quotient by its $i$-th congruence subgroup, we see in a similar way that the restriction of $\psi$ to ${\mth{U}}_{-\alpha,i-2}$ factors into a bijective morphism of $k$-varieties between ${\mth{U}}_{-\alpha,i-2}/{\mth{U}}_{-\alpha,i-1}$ and ${\mth{H}}_{\alpha,i-1}/{\mth{H}}_{\alpha,i}$. By combining all these morphisms we see that $\psi$ is a bijective morphism as well. $\Box$
\end{proof}
\end{proof}

\begin{cor}\label{condpi}
For every uniformizer $\varpi'$ of $R_\alpha$, there exists a unique uniformizer $\varpi$ of $R_{-\alpha}$ such that $b(\varpi)=b'(\varpi')^{-1}$.
\end{cor}

\begin{proof}
Unicity is an immediate consequence of lemma \ref{unicb}. For existence, we remark that $Im(b)$ is the projection on ${\mth{H}}$ (for the order $\alpha<0<-\alpha$) of $[{\mth{U}}_{-\alpha},u_{-\alpha,1}]$, hence also of $[{\mth{U}}_{-\alpha},u_{-\alpha,x}]$ for every $x\in R_{-\alpha}$ of valuation $0$ since $u_{-\alpha,x}=Ad(h_{-\alpha^\vee,\sqrt{x}})u_{\alpha_1}$ for any choice of $\sqrt{x}$. On the other hand, for every uniformizer $\varpi\in R_{-\alpha}$ and every $x\in R_\alpha$ of nonzero valuation, $[u_{-\alpha,\varpi},u_{-\alpha,x}]$ is contained in ${\mth{H}}_{\alpha,3}$ by {\bf (PT4d'')} applied to $U_{-\alpha,f_0(-\alpha)+1}$ and $U_{\alpha,f_0(\alpha)+1}$ (remember that $f_0(\alpha)+f_0(-\alpha)=1$); we thus obtain, by {\bf (PT4d'')} applied this time to ${\mth{U}}_{-\alpha,f_0(-\alpha)+1}$ and ${\mth{U}}_{-\alpha}$, that every element of ${\mth{H}}_{\alpha,2}$ of valuation $2$ is the image by $b$ of some uniformizer of $R_{-\alpha}$. Since a similar assertion is true for $b'$ as well, for every uniformizer $\varpi'$ of $R_\alpha$, if $b'(\varpi')$ is actually of valuation $2$, the corollary is proved.

It remains to check that $b'(\varpi')$ is of valuation exactly $2$. Let $v$ be the valuation of $b'(\varpi')$; by {\bf (PT4d'')} applied to ${\mth{U}}_{\alpha,f_0(\alpha)+v-1}$  and ${\mth{U}}_{-\alpha}$ and a similar reasoning as above, there exists $x\in R_{-\alpha}$ of valuation at least $v-1$ such that $b'(x)=b'(\varpi')$. By the previous unicity lemma (\ref{unicb}), we must then have $x=\varpi'$, which is of valuation $1$; hence we must have $v\leq 2$. On the other hand, we deduce from ${\bf (PT4d'')}$ that $v\geq 2$; hence $v=2$ and the corollary is proved. $\Box$
\end{proof}

From now on we assume that the uniformizers $\varpi$ and $\varpi'$ we have previously chosen satisfy the condition of corollary \ref{condpi}.

Now we prove that $R_\alpha$ and $R_{-\alpha}$ are isomorphic. Since we already know that $\phi$ is a semigroup isomorphism between $(R_\alpha,*)$ and $(R_{-\alpha},*)$ for any choice of $\varpi$ and $\varpi'$, it remains to prove that, $\varpi$ and $\varpi'$ being chosen as above, the resulting $\phi$ is a group isomorphism between $(R_\alpha,+)$ and $(R_{-\alpha},+)$.

Let $\lambda,\mu$ be two elements of $R_{-\alpha}$. We have:
\[u_{-\alpha,\lambda+\mu}u_{\alpha,1}=u_{-\alpha,\mu}u_{-\alpha,\lambda}u_{\alpha,1}\]
\[=u_{-\alpha,\mu}u_{\alpha,a(\lambda)}h_{\alpha^\vee,b(\lambda)}u_{-\alpha,c(\lambda)}.\]
Assume $v(\lambda)=0$ and let $\nu$ be a square root of $a(\lambda)$ in $R_\alpha$; we have:
\[u_{-\alpha,\mu}u_{\alpha,a(\lambda)}=h_{\alpha^\vee,\nu}u_{-\alpha,\phi(a(\lambda))\mu}u_{\alpha,1}h_{\alpha^\vee,\nu^{-1}}\]
\[=u_{\alpha,a(\lambda)a(\phi(a(\lambda))\mu)}h_{\alpha^\vee,b(\phi(a(\lambda))\mu)}u_{-\alpha,\frac{c(\phi(a(\lambda))\mu)}{\phi(a(\lambda))}},\]
from which we deduce the following expressions:
\begin{equation}\label{alm}a(\lambda+\mu)=a(\lambda)a(\phi(a(\lambda))\mu);\end{equation}
\begin{equation}\label{blm}b(\lambda+\mu)=b(\lambda)b(\phi(a(\lambda))\mu);\end{equation}
\begin{equation}\label{clm}c(\lambda+\mu)=c(\lambda)+\phi(\frac{b(\lambda)^2}{a(\lambda)})c(\phi(a(\lambda))\mu).\end{equation}
We will check that these expressions remain true when $v(\lambda)\neq 0$. Let $\lambda_1,\lambda_2$ be any two elements of $R_{-\alpha}$ such that $v(\lambda_1)=v(\lambda_2)=0$ and $\lambda_1+\lambda_2=\lambda$; we have:
\[b(\lambda+\mu)=b(\lambda_1+(\lambda_2+\mu))\]
\[=b(\lambda_1)b(\phi(a(\lambda_1)))(\lambda_2+\mu))\]
\[=b(\lambda_1)b(\phi(a(\lambda_1))\lambda_2)b(\phi(a(\phi(a(\lambda_1))\lambda_2))\phi(a(\lambda_1))\mu)\]
\[=b(\lambda_1+\lambda_2)b(\phi(a(\lambda_1+\lambda_2))\mu)\]
\[=b(\lambda)b(\phi(a(\lambda))\mu).\]
The proof of (\ref{alm}) is obtained by simply replacing $b$ with $a$ in the above expressions. The proof of (\ref{clm}) is also very similar and is left to the reader.

Similarly, we have, for every $\lambda',\mu'\in R_\alpha$:
\begin{equation}\label{aplm}a'(\lambda'+\mu')=a'(\lambda')a'(\phi^{-1}(a'(\lambda'))\mu');\end{equation}
\begin{equation}\label{bplm}b'(\lambda'+\mu')=b'(\lambda')b'(\phi^{-1}(a'(\lambda'))\mu');\end{equation}
\begin{equation}\label{cplm}c'(\lambda'+\mu')=c'(\lambda')+\frac{b'(\lambda')^2}{\phi^{-1}(a'(\lambda'))}c'(\phi^{-1}(a'(\lambda'))\mu').\end{equation}

For every $\lambda\in R_{-\alpha}$, we can set $\mu=-\lambda$ in (\ref{blm}) to obtain:
\[1=b(0)=b(\lambda)b(-\phi(a(\lambda))\lambda).\]
On the other hand, we have, using (\ref{bplm}):
\[1=b'(0)=b'(\phi^{-1}(\lambda))b'(-\phi^{-1}(a'(\phi^{-1}(\lambda))\lambda)).\]
Assume $v(\lambda)=0$; since we know that $b'(\phi^{-1}(\lambda))=b(\lambda)^{-1}$, we obtain:
\[b(-\phi(a(\lambda))\lambda)=b'(-\phi^{-1}(a'(\phi^{-1}(\lambda))\lambda))^{-1}.\]
On the other hand, we have:
\[b(-\phi(a(\lambda))\lambda)=b'(-a(\lambda)\phi^{-1}(\lambda))^{-1}.\]
By unicity (lemma \ref{unicb} applied to $b'$), we obtain:
\[\phi^{-1}(a'(\phi^{-1}(\lambda))\lambda)=a(\lambda)\phi^{-1}(\lambda),\]
hence:
\[a'(\phi^{-1}(\lambda))=\phi(a(\lambda)).\]
Now let $\mu$ be an element of $R_{-\alpha}$ such that either $v(\mu)=0$ or $\phi(a(\lambda)))\mu=\varpi$. We have, using (\ref{blm}) and (\ref{bplm}):
\[b(\lambda+\mu)=b(\lambda)b(\phi(a(\lambda))\mu)\]
\[=b'(\phi^{-1}(\lambda))^{-1}b'(a(\lambda)\phi^{-1}(\mu))^{-1}\]
\[=b'(\phi^{-1}(\lambda))^{-1}b'(\phi^{-1}(a'(\phi^{-1}(\lambda))\mu))^{-1}\]
\[=b'(\phi^{-1}(\lambda)+\phi^{-1}(\mu))^{-1}.\]
Assume now that either $v(\lambda+\mu)=0$ or $\lambda+\mu=\varpi$, we also have:
\[b(\lambda+\mu)=b'(\phi^{-1}(\lambda+\mu))^{-1}.\]
hence by unicity (lemma \ref{unicb} again):
\[\phi^{-1}(\lambda+\mu)=\phi^{-1}(\lambda)+\phi^{-1}(\mu).\]
Using the facts that $\phi(-x)=-\phi(x)$ for every $x$ (which allows us to replace $(\lambda,\mu,\lambda+\mu)$ by either $(\lambda+\mu,-\lambda,\mu)$ or $(\lambda+\mu,-\mu,\lambda)$) and that $\phi$ is a multiplicative semigroup isomorphism, we see that the above equality is true as soon as two of $\lambda,\mu,\lambda+\mu$ have the same valuation $v$ and that the third one has valuation either $v$ or $v+1$. We will now prove by induction on $v'-v$ that it remains true when the first two have valuation $v$ and the third one has any valuation $v'\geq v$; since by ultrametricity every subset of three elements of $R_\alpha$ such that one of them is the sum of the other two satisfies that condition for some $v$, it will prove that $\phi$ is a ring isomorphism.

We can assume that $\mu$ is the element with the greater valuation. Assume then $v'-v\geq 2$ and let $\nu$ be an element of $R_\alpha$ such that $v<v(\nu)<v'$. We have, applying the induction hypothesis at each step:
\[\phi^{-1}(\lambda+\mu)+\phi^{-1}(\nu)=\phi^{-1}(\lambda+\mu+\nu)\]
\[=\phi^{-1}(\lambda)+\phi^{-1}(\mu+\nu)=\phi^{-1}(\lambda)+\phi^{-1}(\mu)+\phi^{-1}(\nu),\]
which proves the desired identity.

In the sequel we will drop $\phi$ and identify $R_\alpha$ with $R_{-\alpha}$ to simplify the notations. Note though that we still do not know whether $\phi$ is algebraic or not.

For every $\lambda\in R_{-\alpha}$ such that $v(\lambda)=0$, we have:
\[c(\lambda)=\lambda a'(\lambda)=\lambda a(\lambda).\]
This is obviously also true when $\lambda=0$.

Let $\lambda,\mu,\lambda+\mu$ be three elements of $R_\alpha$ which are either zero or of valuation $0$; we have, using (\ref{clm}):
\[(\lambda+\mu)a(\lambda+\mu)=\lambda a(\lambda)+\frac{b(\lambda)^2}{a(\lambda)}a(\lambda)\mu a(a(\lambda)\mu),\]
hence, by (\ref{alm}):
\[(\lambda+\mu-\frac{b(\lambda)^2}{a(\lambda)}\mu)a(\lambda+\mu)=\lambda a(\lambda).\]
When $\lambda+\mu=0$ and $\lambda,\mu\neq 0$, we get:
\[(\frac{b(\lambda)^2}{a(\lambda)}\lambda)=\lambda a(\lambda),\]
hence $a(\lambda)^2$ is the class of $b(\lambda)^2$ mod ${\mth{H}}_{\alpha,d}$. On the other hand, we have the following lemma:

\begin{lemme}
The function $a$ takes its values in $1+\varpi R_\alpha$.
\end{lemme}

\begin{proof}
We have for every $\lambda$:
\[[u_{\alpha,-1},u_{-\alpha,\lambda}]=u_{\alpha,-1}u_{-\alpha,\lambda}u_{\alpha,1}u{-\alpha,-\lambda}\]
\[=u_{\alpha,a(\lambda)-1}h_{\alpha^\vee,b(\lambda)}u_{-\alpha,c(\lambda)-\lambda)}\in{\mth{U}}_{\alpha,1}{\mth{H}}_{\alpha,1}{\mth{U}}_{-\alpha,2}.\]
The lemma follows immediately. $\Box$
\end{proof}

Since this is also true for $b$ by its definition, we obtain that $a(\lambda)$ is the class of $b(\lambda)$ mod ${\mth{H}}_{\alpha,d}$. We thus have:
\[a(\lambda+\mu)=\frac{\lambda a(\lambda)}{\lambda+(1-a(\lambda))\mu}=\frac{\lambda a(\lambda)}{\lambda+\mu-a(\lambda)\mu}.\]
Now we can assume $\varpi$ is the unique uniformizer of $R_{-\alpha}$ such that $b(1)=\frac 1{1+\varpi}$. We then also have $a(1)=\frac 1{1+\varpi}$, and for every $\lambda$ such that $v(\lambda)$ and $v(\lambda-1)$ are both zero, we have:
\[a(\lambda)=a(1+(\lambda-1))=\frac{\frac 1{1+\varpi}}{\lambda+\frac{1-\lambda}{1+\varpi}}\]
\[=\frac 1{(1+\varpi)\lambda+1-\lambda}=\frac 1{1+\varpi\lambda}.\]
Since the set of such $\lambda$ is a dense open subset of $R_{-\alpha}$ and since $a$ and $\lambda\mapsto\frac 1{1+\varpi\lambda}$ are both maps of $k$-varieties, the equality holds for every $\lambda\in R_{-\alpha}$. We obtain similarly that for every $\lambda\in R_{-\alpha}$, we have:
\[b(\lambda)=\frac 1{1+\varpi\lambda},\;\;c(\lambda)=\frac \lambda{1+\varpi\lambda}.\]
By conjugating both sides by a suitable element of ${\mth{H}}_\alpha$ if needed, we deduce from this that for every $x,y\in R_\alpha$, we have:
\begin{equation}\label{scom}u_{-\alpha,x}u_{\alpha,y}=u_{\alpha,\frac y{1+\varpi xy}}h_{\alpha^\vee,\frac 1{1+\varpi xy}}u_{-\alpha,\frac x{1+\varpi xy}}.\end{equation}

Now we finally come to the main result of this subsection:

\begin{prop}\label{soca}
Assume ${\mth{G}}$ is of semisimple parahoric type, solvable and of semisimple rank $1$. Then $R_\alpha$ and $R_{-\alpha}$ are isomorphic, and there exists a henselian nonarchimedean local field $F$ of residual field $k$ and an algebraic group $\underline{G}$ of type $A_1$ defined and split over $F$ such that ${\mth{G}}$ is isomorphic to the quotient of an Iwahori subgroup of the group $G$ of $F$-points of $\underline{G}$ by its $d$-th congruence subgroup.
\end{prop}

\begin{proof}
Let $F$, $\mathcal{O}_F$, $\mathcal{p}_F$ and $\underline{G}$ be defined as in the previous subsection. Since $\underline{\mth{G}}$ is of semisimple parahoric type, $\underline{G}$ is semisimple; it is then isomorphic to either $SL_2$ or $PGL_2$. Let $G$ be the group of $F$-points of $\underline{G}$, let $I$ be the Iwahori subgroup of the elements of $G(\mathcal{O}_F)$ which are upper triangular mod $\mathfrak{p}_F$, and set ${\mth{G}}_0=I/I^d$. Consider now the map:
\[u_{\alpha,\lambda}h_{\alpha^\vee,\mu}u_{-\alpha,\nu}\mapsto \left(\begin{array}{cc}1&\lambda\\&1\end{array}\right) \left(\begin{array}{cc}\mu\\&\mu^{-1}\end{array}\right) \left(\begin{array}{cc}1\\\varpi\nu&1\end{array}\right)\]
from ${\mth{G}}$ to ${\mth{G}}_0$; with the help of easy matrix computations we see that it is an isomorphism. $\Box$
\end{proof}

Now we can also prove the following proposition

\begin{prop}\label{prisec2}
Assume $G$ is any $k$-split group of parahoric type. Then {\bf (PT4d')}) implies ({\bf (PT4d'')}.
\end{prop}

\begin{proof}
Note first that since, in the general case, all subgroups of $\underline{\mth{G}}$ showing up in the statement of {\bf (PT4d')} (resp. {\bf (PT4d'')}) are contained in some given subgroup of semisimple parahoric type of rank $1$ of $\underline{\mth{G}}$, we do not lose any generality by assuming that $\underline{\mth{G}}$ itself is such a group.

Assume $\underline{\mth{G}}$ satisfies {\bf (PT4d')}. Then all results of tist section are valid with $k$ being replaced with $\underline{k}$. In particular, in the nonsolvable case, we still have, with the same definitions as in the previous subsection::
\[u_{-\alpha,\lambda}u_{\alpha,1}=u_{\alpha,\frac 1{1+\lambda}}h_{\alpha^\vee,\frac 1{1+\lambda}}u_{-\alpha,\frac\lambda{1+\lambda}}.\]
Moreover, by choosing the elements $u_{\alpha,1}$ and $u_{-\alpha,1}$ inside ${\mth{G}}$, we see that every element of ${\mth{H}}_\alpha$ shows up as the second term of the right-hand side for some $\lambda\in R_\alpha$, hence for every $i,j$, the projection of $[{\mth{U}}_{\alpha,i},{\mth{U}}_{-\alpha,j)}]$ on ${\mth{H}}$ is ${\mth{H}}_{\alpha,i+j}$ (for the Iwahori decomposition). Moreover, we deduce from proposition \ref{htrans2} that we have, for every suitable $i,j,k$:
\[[\underline{\mth{H}}_{i+j},\underline{\mth{U}}_k]=[\underline{\mth{H}}_{i+j},u],\]
with $u$ being any element of valuation $k$ of $\underline{\mth{U}}_\alpha$, that we can very well choose inside ${\mth{U}}_\alpha$; the identity $[{\mth{H}}_{\alpha,i+j},{\mth{U}}_{\alpha,k}]={\mth{U}}_{\alpha,i+j+k}$ follows. Hence {\bf (PT4d'')} holds.

The solvable case can be proved in a similar way, by also taking the uniformizers $\varpi$ and $\varpi'$ of (\ref{scom}) inside respectively $R_\alpha$ and $R_{-\alpha}$, given that the identity $b(\varpi')=b'(\varpi)^{-1}$ of corollary \ref{condpi} ensures that it is possible. $\Box$
\end{proof}

\section{A more detailed study of the commutator relations}

Now we return to the case where $\Phi$ is of any rank. In \cite{chev}, Chevalley has computed some structural constants for reductive groups, which arise when we study in detail the commutator relations. Our goal in this section is to find similar constants for groups of parahoric type; we will prove that these constants are in fact similar to the Chevalley constants. We will also study the action of the elements of ${\mth{H}}$ on the root subgroups.

\subsection{Some preliminary results}

First we establish some relations between the rings associated to the various elements of $\Phi$.
%Let $u_\varpi$ be a uniformizer of $R_\alpha$; we set $R'_\alpha=R_\alpha/u_\varpi^dR_\alpha$. If $\alpha\not\in\Psi$, since ${\mth{U}}_\alpha$ is of dimension $d$, we must have $u_\varpi^d=0$, hence $R'_\alpha=R_\alpha$.
In the sequel, for every $\alpha\in\Phi$ and for every $x\in R_\alpha$, we will denote by $u_{\alpha,x}$ the element $x$ viewed as an element of the group ${\mth{U}}_\alpha$. Of course $u_{\alpha,x}$ depends on the choice of $u_{\alpha,1}$.

Remember that we are still assuming that $p\neq 2$ and that if $p=3$, $\underline{\mth{G}}$ is not of type $G_2$. To simplify the notations we will also assume that $\Phi$ is irreducible; the general case is not significantly harder.

\begin{prop}\label{anniso3}
Assume either $\underline{\mth{G}}$ is not of type $A_2$ or $p\neq 3$. Let $\beta,\gamma$ be two linearly independent elements of $\Phi$.
\begin{itemize}
\item {\it (i)} Assume either $\underline{\mth{G}}$ is of type other than $A_2$ or $f_0(\beta+\gamma)=f_0(\beta)+f_0(\gamma)$. Then the group $Z_{\mth{H}}({\mth{U}}_\gamma)$ acts transitively on the set of elements of valuation $i$ of ${\mth{U}}_\beta$ for every $i$.
\item {\it (ii)} Assume $\alpha=\beta+\gamma$ is a root. Then all three rings $R'_\alpha$, $R'_\beta$ and $R'_\gamma$ are isomorphic to each other. Moreover, when all three of $\alpha,\beta,\gamma$ belong to $\Psi$, the rings $R_\alpha$, $R_\beta$ and $R_\gamma$ are isomorphic to each other as well.
\end{itemize}
\end{prop}

\begin{proof}
We first prove the following lemma:

\begin{lemme}\label{1ind2}
Let $\alpha,\beta,\gamma$ be three elements of $\Phi$ such that $\alpha=\beta+\gamma$. Assume that {\it (i)} holds for $\beta$ and $\gamma$ and that $f_0(\alpha)=f_0(\beta)+f_0(\gamma)$; then $R'_\alpha$ and $R'_\beta$ (resp. $R_\alpha$ and $R_\beta$ when all three roots belong to $\Psi$) are isomorphic.
\end{lemme}

\begin{proof}
According to the commutator relations, for every $h\in Z_{\mth{H}}({\mth{U}}_\gamma)$, every $u\in{\mth{U}}_\beta$ and every $u'\in{\mth{U}}_\gamma$, we then have:
\begin{equation}\label{adhu}[Ad(h)u,u']=Ad(h)([u,u']).\end{equation}
%where the $\alpha'$ are the roots of the form $i\beta+j\gamma$, with $i,j\geq 1$ and $\alpha'\neq\alpha$.
Moreover, since $f_0(\alpha)=f_0(\beta)+f_0(\gamma)$, it is possible to chose $u_{\alpha,1}$, $u_{\beta,1}$ and $u_{\gamma,1}$ in such a way that the projection of $[u_{\beta,1},u_{\gamma,1}]$ on ${\mth{U}}_\alpha$ is $u_{\alpha,1}$.

Let $\phi$ be the application from ${\mth{U}}_\beta$ to ${\mth{U}}_\alpha$ defined the following way: for every $u\in{\mth{U}}_\beta$, $\phi(u)$ is the projection on ${\mth{U}}_\alpha$ of $[u,u_{\gamma,1}]$; We deduce from {\it (i)} that the action of $Z_{\mth{H}}({\mth{U}}_\gamma)$ on the set of elements of valuation $i$ of ${\mth{U}}_\beta$ is transitive for every $i$.

We deduce from {\bf (PT4c'')} that $\phi$ factors into a group morphism between $R'_\beta$ and $R'_\alpha$, that by a slight abuse of notation we will also denote by $\phi$, and by the assumptions we have made, $\phi(u_{\beta,1})=u_{\alpha,1}$. To prove that $\phi$ is a ring morphism, it remains to check that $\phi$ preserves the product of two elements.

Let $u,u'$ be in ${\mth{U}}_\beta$. Assume first $v(u)=0$, and let $h_u$ be an element of $Z_{\mth{H}}({\mth{U}}_\gamma)$ such that $u=Ad(h_u)u_{\beta,1}$; from (\ref{adhu}) we deduce that we also have $\phi(u)=Ad(h_u)\phi(u_{\beta,1})$, hence:
\[\phi(u.u')=\phi(Ad(h_u)u')=Ad(h_u)\phi(u')=\phi(u).\phi(u').\]
If now $v(u)>0$, we have:
\[\phi(u.u')=\phi((u+u_{\beta,1}).u'-u')=(\phi(u)+\phi(u_{\beta,1})).\phi(u')-\phi(u')=\phi(u).\phi(u').\]

It remains to prove that $\phi$ is an isomorphism. First we prove that $\phi$ is injective. Consider the ideal $Ker(\phi)$ of $R'_\beta$; by proposition \ref{ideal}, there exists $i\in\{0,\dots,d\}$ such that $Ker(\phi)=R'_{\beta,i}$. On the other hand, for every $u\in R'_\gamma$ of valuation $0$, by proposition \ref{htrans2} there exists $h\in{\mth{H}}$ such that $hu_{\gamma,1}h^{-1}=u$, and we have:
\[[R'_{\beta,i},u]=[hR'_{\beta,i}h^{-1},hu_{\gamma,1}h^{-1}]=h'[R'_{\beta,i},u_{\gamma,1}]h'^{-1}.\]
Since the set of elements of valuation $0$ of $R'_\gamma$ generates $R'_\gamma$ as a group, we obtain that the projection of $[R'_{\beta,1},u]$  on $R'_\alpha$ is zero for every $u\in R'_\gamma$; by {\bf (PT4c'')} we must then have $i=d$, hence $Ker(\phi)=0$ and $\phi$ is injective.

Now we prove surjectivity. We deduce easily from {\bf (PT4c'')} that $\phi$ factors into a ring morphism between $R'_\beta/R'_{\beta,1}$ and $R'_\alpha/R'_{\alpha,1}$; since both rings are isomorphic as algebraic $k$-groups over the field $k$ and $k$ is perfect, that morphism is an isomorphism. On the other hand, by {\bf (PT4c'')} $\phi(\varpi)$ is of nonzero valuation; let then $i$ be its valuation. Then $\phi(\varpi^{d-1})$ is of valuation $i(d-1)$; since $\phi$ is injective, this is possible only if $i=1$. Hence $\phi(\varpi)$ is an uniformizer of $R_\alpha$.

Since $\varpi$ (resp. $\phi(\varpi)$) and any set of representatives of $R'_\beta/R'_{\beta,1}$ (resp. $R'_\alpha/R'_{\alpha,1}$) generate $R'_\beta$ (resp. $R'_\alpha$) as a ring, we obtain that $\phi(R'_\beta)=R'_\alpha$. Hence $\phi$ is a ring isomorphism between $R'_\beta$ and $R'_\alpha$ and lemma \ref{1ind2} is proved. The proof that when all three of $\alpha$, $\beta$ and $\gamma$ belong to $\Psi$, $R_\alpha$ and $R_\beta$ are isomorphic is exactly similar. $\Box$
\end{proof}

\begin{cor}
Assume {\it (i)} holds. Then {\it (ii)} holds as well.
\end{cor}

\begin{proof}
Let $\alpha,\beta,\gamma$ be three elements of $\Phi$ such that $\alpha=\beta+\gamma$. When $f_0(\alpha)=f_0(\beta)+f_0(\gamma)$, the fact that the rings are isomorphic is simply lemma \ref{1ind2}. Assume now $f_0(\alpha)<f_0(\beta)+f_0(\gamma)$. Then by proposition \ref{f0ab3}, $f_0(\beta)=f_0(\alpha)+f_0(-\gamma)$; since we know by the case of rank $1$ that $R'_\gamma$ and $R'_{-\gamma}$ are isomorphic, we are reduced to the previous case. $\Box$
\end{proof}

Now we prove {\it (i)}. We start by the following lemma:

\begin{lemme}
Assume $\beta$ and $\gamma$ are strongly orthogonal. Then ${\mth{H}}_\beta\subset Z_{\mth{H}}({\mth{U}}_\gamma)$.
\end{lemme}

\begin{proof}
Let $s_\beta$ be the reflection associated to $\beta$ in the Weyl group of $\Phi$. Since $\gamma+\beta$ is not a root, $\gamma-\beta=s_\beta(\gamma+\beta)$ is not a root either, and if $y_\beta,y_{-\beta},y_\gamma$ are elements of respectively $Lie({\mth{U}}_\beta)$, $Lie({\mth{U}}_{-\beta})$ and $Lie({\mth{U}}_\gamma)$, in the Jacobi identity:
\[[[y_\beta,y_{-\beta}],y_\gamma]+[[y_{-\beta},y_\gamma], y_\beta]+[[y_\gamma,y_\beta],y_{-\beta}]=0,\]
the last two terms of the left-hand side are zero, hence the first one must be zero too. Hence $Lie({\mth{H}}_{\beta,1})$ and $Lie({\mth{U}}_\gamma)$ commute, which implies that ${\mth{H}}_{\beta,1}$ is contained in $Z_{\mth{H}}({\mth{U}}_\gamma)$. Moreover, since $\beta$ and $\gamma$ are orthogonal, $\gamma$ is trivial on ${\mth{T}}_\beta$, hence ${\mth{T}}_\beta\subset Z_{\mth{H}}({\mth{U}}_\gamma)$, which proves the lemma. $\Box$
\end{proof}

\begin{cor}
Assume $\Phi$ is not of type $A_r$, and if $\Phi$ is of type $B_r$, $C_r$ or $F_4$, $\beta$ and $\gamma$ are both long. Then {\it (i)} holds.
\end{cor}

\begin{proof}
Let $\beta^\perp$ (resp. $\gamma^\perp$) be the subsystem of the elements of $\Phi$ which are orthogonal to $\beta$ (resp. $\gamma$). Since $\Phi$ is not of type $A_r$, we see by examining its extended Dynkin diagram that one of the long roots contained in some given extended set of simple roots of $\Phi$ is strongly orthogonal to $r-1$ other elements of that set (for example, any special element of that extended set of simple roots satisfies that condition); since all long roots of $\Phi$ are conjugates by elements of the Weyl group of $\Phi$, it implies that $\beta^\perp$ and $\gamma^\perp$ are of rank $r-1$, which implies in particular that they are nonempty; since $\beta$ and $\gamma$ are linearly intependent, $\gamma^\perp$ cannot be contained in $\beta^\perp$, hence contains at least one element $\alpha'$ which is not orthogonal to $\beta$. On the other hand, it is easy to check (details are left to the reader) that with the hypotheses we have made, every element of $\beta^\perp$ (resp. $\gamma^\perp$) is strongly orthogonal to $\beta$ (resp. $\gamma$); we deduce then from the previous lemma that the dimension of $Z_{\mth{H}}({\mth{U}}_\beta)$ is at least $(d+1)(r-1)$; its codimension in ${\mth{H}}$ is then at most $d+1$, hence exactly $d+1$, which proves that ${\mth{H}}_{\alpha'}Z_{\mth{H}}({\mth{U}}_\beta)={\mth{H}}$. Since ${\mth{H}}_{\alpha'}\subset Z_{\mth{H}}({\mth{U}}_\gamma)$, the corollary now follows from proposition \ref{htrans2}. $\Box$
\end{proof}

Now we assume that $\Phi$ is of type $B_r$, $C_r$ or $F_4$ and either $\beta$ or $\gamma$ is short. If $\beta$ is strongly orthogonal to $\gamma$, then ${\mth{H}}_\beta\subset Z_{\mth{H}}({\mth{U}}_\gamma)$ and {\it (i)} is an immediate consequence of proposition \ref{htrans2}. If $\beta$ and $\gamma$ are not strongly orthogonal, we have three cases to consider:
\begin{enumerate}
\item either $\beta$ or $\gamma$ is long (which makes one long root and one short root);
\item both $\beta$ and $\gamma$ are short, and $\beta$ and $\gamma$ generate a root subsystem of type $B_2$ of $\Phi$ (in which case they are orthogonal but not strongly orthogonal);
\item both $\beta$ and $\gamma$ are short, and $\beta$ and $\gamma$ generate a root subsystem of type $A_2$ of $\Phi$.
\end{enumerate}
We first observe that in the first case, $\beta$ and $\gamma$ always generate a root subsystem of type $B_2$ of $\Phi$. We will thus prove the first and second case simultaneously.

Assume we are in the second case: $\beta$ is a short root of $\Phi$ which is orthogonal but not strongly orthogonal to $\gamma$. Then $\beta\pm\gamma$ is a long root, and the pairs($ \beta-\gamma,\gamma)$ and $(\gamma-\beta,\gamma)$ correspond to the first case.

\begin{lemme}
The groups ${\mth{H}}_{\beta+\gamma}$ and ${\mth{H}}_{\beta-\gamma}$ are isomorphic to both ${\mth{H}}_\beta$ and ${\mth{H}}_\gamma$.
\end{lemme}

\begin{proof}
We will prove the result for ${\mth{H}}_{\beta+\gamma}$, the proof for ${\mth{H}}_{\beta-\gamma}$ being similar. By the proof of proposition \ref{basprod}, we have ${\mth{H}}_\beta{\mth{H}}_\gamma={\mth{H}}_{\beta+\gamma}{\mth{H}}_\gamma$. On the other hand, we deduce from {\bf (PT23} that ${\mth{H}}_\beta\cap{\mth{H}}_\gamma=\{1\}$; for every $x\in{\mth{H}}_{\beta-\gamma}$, there exists then a unique pair $(y,z)\in{\mth{H}}_\beta\times{\mth{H}}_\gamma$ such that $x=yz$; checking that $x\mapsto y$ and $x\mapsto z$ are group isomorphisms is then straightforward. $\Box$
\end{proof}

Let then $R^*=R_\beta^*=R_\gamma^*=R_{\beta+\gamma}^*=R_{\beta-\gamma}^*$ be that group. For every $x\in R^*$, let $h_{(\beta+\gamma) ^\vee,x}$ be an element of ${\mth{H}}_{\beta+\gamma}$ acting on $R_{\beta+\gamma}$ by multiplication by $x^2$; it then acts on $R_\beta$ (resp. $R_\gamma$) by multiplication by a constant that we will call $c(x)$ (resp. $d(x)$), hence on $R_{\beta+\gamma}$ by multiplication by $c(x)d(x)$, which implies that we have $c(x)d(x)=x^2$.  The element $h_{(\beta+\gamma)^\vee,x}$ also acts on $R_{-\beta}$ (resp. $R_{-\gamma}$) by multiplication by $c(x)^{-1}$ (resp. $d(x)^{-1}$), hence on $R_{\beta-\gamma}$ by multiplication by $c(x)d(x)^{-1}$. Since $\beta+\gamma$ and $\beta-\gamma$ are strongly orthogonal we must then have $c(x)d(x)^{-1}=1$, hence $c(x)=d(x)$, for every $x$; we thus obtain $c(x)^2=x^2$, and by continuity of the map $x\mapsto c(x)$ we must have $c(x)=x$ for every $x$. In particular, ${\mth{H}}_{\beta+\gamma}$ is contained in $Z_{\mth{H}}({\mth{U}}_{\beta-\gamma})$ and acts transitively on the set of elements of any given valuation of ${\mth{U}}_\beta$ (resp. ${\mth{U}}_\gamma$), which proves (with $\beta+\gamma$ playing here the role of the long $\gamma$) the "$\gamma$ long" part of the first case.

We can make the exact same reasoning with $\gamma$ being replaced by $-\gamma$, and we obtain that $h_{(\beta-\gamma)^\vee,x}$ acts on $R_{\beta}$ by multiplication by a constant $c'(x)$ and on $R_\gamma$ by multiplication by $c'(x)^{-1}$, and we also obtain $c'(x)=x$ for every $x$ by the same reasoning as above. We thus obtain the following result:

\begin{lemme}\label{myz}
For every $y,z\in R^*$, the elements $h_{(\beta+\gamma)^\vee,y}h_{(\beta-\gamma)^\vee,z}$ of ${\mth{H}}$ acts of ${\mth{U}}_\beta$ (resp. ${\mth{U}}_\gamma$, ${\mth{U}}_{\beta+\gamma}$, ${\mth{U}}_{\beta-\gamma}$) by multiplication by $yz$ (resp. $yz^{-1}$, $y^2$, $z^2$). In particular, for every $x$, $h_{\beta+\gamma)^\vee,x}h_{(\beta-\gamma)^\vee,a,x}$ belongs to $Z_{\mth{H}}({\mth{U}}_\gamma)$ and acts on ${\mth{U}}_\beta$ (resp. ${\mth{U}}_{\pm\beta\pm\gamma}$ by multiplication by $x^2$ (resp. either $x^2$ or $x^{-2}$).
\end{lemme}

Since every element of $R^*$ admits two opposite square roots by lemma \ref{hsq}, the first two cases follow.

(Remark: We will prove later that, as we can expect, the set of elements of ${\mth{H}}$ of the form $h_{(\beta+\gamma)^\vee,x}h_{(\beta-\gamma)^\vee,x}$ is actually the group ${\mth{H}}_\beta$; we do not need this fact for the moment.)

It remains to consider the case where $\beta$ is short and $\beta$ and $\gamma$ generate a root subsystem of type $A_2$ of $\Phi$. This leaves out the case of a system of type $B_r$, in which two short roots are always orthogonal to each other, and leaves us with the cases $C_r$, $r\geq 3$, and $F_4$; in both cases, the Dynkin diagram of $\Phi$ contains two short roots which are next to each other and such that one of them is neighboring a long root that we will call $\alpha'$. We may thus assume that $\beta$ is the one of these two short root neighboring $\alpha'$ and, by eventually replacing it with $-\beta-\gamma$, that $\gamma$ is the other one. This is compatible with the action of the Weyl group, by the folowing lemma:

\begin{lemme}
Assume $\beta$ and $\gamma$ are two short roots of a root system $\Phi$ of type either $C_r$, $r\geq 3$, or $F_4$ which generate a subsystem of type $A_2$ of $\Phi$. There exists then an element of the Weyl group $W_\Phi$ of $\Phi$ which fixes $\beta$ and sends $\gamma$ to $-\beta-\gamma$.
\end{lemme}

\begin{proof}
Let $s_\beta$ be the reflection associated to $\beta$; we have $s_\beta(\beta)=-\beta$ and, since $\beta$ and $\gamma$ generate a subsystem of type $A_2$ of $\Phi$, $s_\beta(\gamma)=\beta+\gamma$. On the other hand, by \cite[plates  III and VIII]{bou}, when $\Phi$ is of type either $C_d$ of $F_4$, its Weyl group contains an element $w_0$ which sends every root to its opposite. The element $w_0s_\beta$ of $W_\Phi$ then satisfies the conditions of the lemma. $\Box$
\end{proof}

Consider now the group ${\mth{H}}_{\alpha'}$; since $\alpha'$ and $\gamma$ are strongly orthogonal, it is contained in $Z_{\mth{H}}({\mth{U}}_\gamma)$, and since $\alpha'$ and $\beta$ generate a root subsystem of type $B_2$ of $\Phi$, we can apply the previous case to $\beta$ and $\beta+\alpha'$ to obtain that ${\mth{UH}}_{\alpha'}$ acts transitively on the set of elements of any given valuation of ${\mth{U}}_{\beta}$. Hence the first condition of the theorem is satisfied for $\beta$ and $\gamma$, which completes the proof for $\Phi$ of type different from $A_r$.

Now we consider the case where $\Phi$ is of type $A_r$, $r\geq 3$. Let $\alpha'$ be an element of $\Phi$ such that $\alpha'$, $\beta$ and $\gamma$ are linearly independent and that $\alpha'+\beta$ is a root; then $\alpha'$ and $\gamma$ are orthogonal and, since $\Phi$ is simply-laced, $\beta-\alpha'$ cannot be a root, hence if $y_{\alpha'}$, $y_{-\alpha'}$ and $y_\beta$ are elements of respectively $Lie({\mth{U}}_{\alpha'})$, $Lie({\mth{U}}_{-\alpha'})$ and $Lie({\mth{U}}_\beta)$, the last two commute and the Jacobi identity between then reduces to:
\[[[y_{\alpha'},y_{-\alpha'}],y_\beta]=[[y_{\alpha'},y_\beta],y_{-\alpha'}].\]
We thus deduce from {\bf (PT4c'')} that $[{\mth{H}}_{\alpha,1},{\mth{U}}_\beta]$ contains ${\mth{U}}_{\gamma,f_0(\alpha)+1}$; the proof that $Z_{\mth{H}}({\mth{U}}_\gamma)\supset{\mth{H}}_{\alpha'}$ acts transitively on the set of elements of ${\mth{U}}_\beta$ of any given valuation is then identical to the proof of proposition \ref{htrans2}.

We finally consider the case where $\Phi$ is of type $A_2$ and $p\neq 3$; in that case, no convenient $\alpha'$ exists and we have to proceed differently. Remember that we are now assuming that $f_0(\beta+\gamma)=f_0(\beta)+f_0(\gamma)$.

\begin{lemme}\label{grobeta}
We can assume that $f_0((-\gamma)=f_0(\beta)+f_0(-\gamma-\beta)$ as well.
\end{lemme}

\begin{proof}
This is a consequence of proposition \ref{f0ab} if $\Psi=\Phi$, and of proposition \ref{f0ab2} if $\Psi$ is empty. Now we consider the case where $\Psi$ is of rank $1$. Ist is easy to check that $\Psi$ must be either $\{\pm\beta\}$ or $\{\pm\gamma\}$ for the identity $f_0(\beta+\gamma)=f_0(\beta)+f_0(\gamma)$ to hold; we will thus assume that $\Psi=\{\pm\beta\}$, the other case being treated by simply switching $\beta$ and $\gamma$. The result follows then from proposition \ref{f0ab} once again. $\Box$
\end{proof}

No we go back to the proof of proposition \ref{anniso3}. By the Jacobi identity, for every $y_\beta\in Lie({\mth{U}}_\beta)$, every $y_\gamma\in Lie({\mth{U}}_\gamma)$ and every $y_{-\beta-\gamma}\in Lie({\mth{U}}_{-\beta-\gamma})$, we have, using the fact that $\gamma-\beta$ and $2\gamma+\beta$ are not elements of $\Phi$:
\[[[y_\gamma,y_{-\beta-\gamma}],y_\beta],y_\gamma]=[[y_\gamma,y_{-\beta-\gamma}],[y_\beta,y_\gamma]];\]
\[[[y_\beta,y_\gamma],y_{-\beta-\gamma}],y_\gamma]=[[y_\beta,y_\gamma],[y_{-\beta-\gamma},y_\gamma]].\]
The right-hand sides of the above equalities being equal to each other, we obtain that $[[y_\gamma,y_{-\beta-\gamma}],y_\beta]-[[y_\beta,y_\gamma],y_{-\beta-\gamma}]$ commutes with $y_\gamma$. Since this is true for every $y_\beta$ and every $y_{-\beta-\gamma}$, we obtain that for every $y_\gamma$, the subspace of $Lie({\mth{H}})$ generated by these elements when $y_\beta$ (resp. $y_{-\beta-\gamma}$) runs over $Lie({\mth{U}}_\beta)$ (resp. $Lie({\mth{U}}_{-\beta-\gamma})$ commutes with $y_\gamma$. On the other hand, we deduce from proposition \ref{htrans2} that ${\mth{H}}$ acts transitively on the sets of elements of a given valuation of both $Lie({\mth{U}}_\beta)$ and $Lie({\mth{U}}_{-\beta-\gamma})$, which implies that the aforementioned subspace depends only on the valuation of $y_\gamma$;, we will now assume that $y_\gamma$ is of valuation $0$. Our subspace then commutes witn $y_\gamma$ for every $y_\gamma$ of valuation $0$, hence for every $y_\gamma$ by density.

With a similar reasoning as above with $(\beta,\gamma,-\beta-\gamma)$ being replaced with $(-\beta-\gamma,\beta,-\gamma)$ (which satisfy $f_0(-\gamma)=f_0(-\beta-\gamma)+f_0(\beta))$ by lemma \ref{grobeta}), we obtain that $[[y_\beta,y_\gamma],y_{-\beta-\gamma}]-[[y_{-\beta-\gamma},y_\beta],y_\gamma]$ commutes with $y_\beta$ for every $y_\beta,y_\gamma,y_{-\beta-\gamma}$. On the other hand, by Jacobi once again, we have:
\[\big[[[y_\beta,y_\gamma],y_{-\beta-\gamma}\big]+\big[[y_\gamma,y_{-\beta-\gamma}],y_\beta\big]+\big[[y_{-\beta-\gamma},y_\beta],y_\gamma\big]=0,\]
which implies:
\[\big[[y_\gamma,y_{-\beta-\gamma}],y_\beta\big]-\big[[y_\beta,y_\gamma],y_{-\beta-\gamma}\big]=\frac 32\big[[y_\gamma,y_{-\beta-\gamma}],y_\beta\big]-\frac 12(\big[[y_\beta,y_\gamma],y_{-\beta-\gamma}\big]-\big[[y_{-\beta-\gamma},y_\beta],y_\gamma\big]),\]
hence:
\[\big[\big[[y_\gamma,y_{-\beta-\gamma}\big],y_\beta\big]-\big[[y_\beta,y_\gamma],y_{-\beta-\gamma}\big],y_\beta\big]=\frac 32\big[\big[[y_\gamma,y_{-\beta-\gamma}],y_\beta\big],y_\beta\big];\]
since $[[y_\gamma,y_{-\beta-\gamma}],y_\beta]\in Lie({\mth{H}}_\beta)$, and since $p\neq 3$, we conclude with (PT4d) and proposition \ref{htrans2}. $\Box$
\end{proof}

\begin{lemme}\label{sumroots}
Let $\Delta$ be an extended system of simple roots of $\Phi$. Then every $\alpha\in\Phi$ is a sum of elements of $\Delta$.
\end{lemme}

\begin{proof}
Since $\Delta$ is the union of a set of simple roots of $\Phi$ and of the negative $\alpha_0$ of the highest root relative to that set of simple roots, every positive (resp. negative) element of $\Phi$ can be written as a sum of simple roots (resp. as a sum of simple roots and $\alpha_0$). $\Box$
\end{proof}

\begin{cor}\label{isoconn}
Assume $\Phi$ is irreducible. Then all rings $R'_\alpha$, $\alpha\in\Phi$, are isomorphic to each other.
\end{cor}

\begin{proof}
Since the case of rank $1$ groups has already been dealt with in the previous section, we can assume $\Phi$ is not of type $A_1$. Let $\Delta$ be any extended set of simple roots of $\Phi$ (i.e. the union of a set of simple roots with the singleton containing the inverse of the highest root of $\Phi$ with respect to that set); since by (lemma \ref{sumroots}, every element $\alpha$ of $\Phi$ can be written as a sum of elements of $\Delta$, assuming $\alpha$ itself does not belong to $\Delta$, we deduce from \cite[\S 1, proposition 19]{bou} that there exists $\delta\in\Delta$ such that $\alpha-\delta$ is a root. Applying proposition \ref{anniso3} we obtain that $R'_\alpha$ and $R'_\delta$ are isomorphic. It is then enough to show that all $R'_\delta$, $\delta\in\Delta$, are isomorphic to each other.

Since $\Phi$ is irreducible and not of type $A_1$, its extended Dynkin diagram is connected and for every $\delta,\delta'\in\Delta$ such that there is an edge between them, $\delta+\delta'$ is a root, hence $R'_\delta$ and $R'_{\delta'}$ are isomorphic by proposition \ref{anniso3}. Hence all $R'_\delta$ are isomorphic and the corollary is proved. $\Box$
\end{proof}

Since the root system of the quotient ${\mth{G}}/R_u({\mth{G}})$ is $\Psi$, the Weyl group $N_{\mth{G}}({\mth{H}})/{\mth{H}}$ can be identified to the Weyl group $W_\Psi$ of $\Psi$, and we have:

\begin{prop}
Let $w$ be an element of $W_\Psi$ and let $n$ be a representative of $w$ in $N_{\mth{G}}({\mth{H}})$. Then $Ad(n)$ acts on ${\mth{H}}$, and, for every $\alpha\in\Phi$, induces an isomorphism between $R_\alpha$ and $R_{w\alpha}$ (resp. between $R'_\alpha$ and $R'_{w\alpha}$).
\end{prop}

\begin{proof}
Let $\alpha$ be any element of $\Phi$. We already know that $Ad(n)$ is a group isomorphism between ${\mth{U}}_\alpha$ and ${\mth{U}}_{w\alpha}$. The fact that it induces ring isomorphisms between $R_\alpha$ and $R_{w\alpha}$ and between $R'_\alpha$ and $R'_{w\alpha}$  is proved the same way as in proposition \ref{anniso3}. $\Box$
\end{proof}

Now consider the rings $R_\alpha$ when they are different from $R'_\alpha$, which is true if and only if $\alpha\in\Psi$. We can prove, always the same way, the following proposition:

\begin{prop}\label{isoconn2}
Let $\alpha,\beta$ be two elements of $\Psi$. Assume $\alpha$ and $\beta$ belong to the same irreducible component of $\Psi$; then $R_\alpha$ and $R_\beta$ are isomorphic. In particular, when $\Psi$ is irreducible, all the $R_\alpha$, $\alpha\in\Psi$, are isomorphic to each other.
\end{prop}

Unfortunately this may not be true anymore when $\Psi$ is reducible. We will give an explicit counterexample at the end of the paper.

\subsection{Definition and first properties of the constants}

Now we have to define our structural constants. Since we make use of proposition \ref{anniso3}, we also have to assume, in addition to the hypotheses we have already made, that when $p=3$, $\Phi$ is not of type $A_2$.

We will also assume that all the $R_\alpha$, $\alpha\in\Psi$, are isomorphic to each other. This allows us to introduce a nonarchimedean henselian local field $F$ such that if $\mathcal{O}_F$ is the ring of integers of $F$ and $\mathfrak{p}_F$ the maximal ideal of $\mathcal{O}_F$, every $R_\alpha$ is isomorphic to the quotient $\mathcal{O}_F/\mathfrak{p}_F^i$, $i$ being either $d+1$ or $d$ depending on whether $\alpha$ lies in $\Psi$ or not. We will also set $R=\mathcal{O}_F/\mathfrak{p}_F^{d+1}$; for every $\alpha\in\Phi$, ${\mth{H}}_\alpha$ is then isomorphic to the unit ring $R^*$ of $R$.

Let $\varpi$ be a uniformizer of $F$. For every $\alpha\in\Phi$, we define the element $u_{\alpha,x}$, $x\in\varpi^{f_0(\alpha)}\mathcal{O}_F/\varpi^{d+1-f_0(-\alpha)}\mathcal{O}_F$ the following way: choose $u_{\alpha,\varpi^{f_0(\alpha)}}$ arbitrarily (for the moment) among the elements of ${\mth{U}}_\alpha$ which are invertible as elements of the ring $R_\alpha$, and for every $x\in\varpi^{f_0(\alpha)}\mathcal{O}_F$, set $u_{\alpha,x}$ as the element of ${\mth{U}}_\alpha$ corresponding to the class of $x$ mod $\varpi^{d+1-f_0(-\alpha)}\mathcal{O}_F$.

Note that this definition is not consistent with the one we have previously used in the case of solvable groups of rank one, but it has the advantage of greatly simplifying the notations when dealing with the Chevalley constants . When $\alpha\in\Psi$, with a suitable choice of $u_{\pm\alpha,\varpi^{f_0(\pm\alpha)}}$, we can always manage to make the proposition \ref{ualphac} still work. When $\alpha\not\in\Psi$, it is not hard to see that the expressions we have found for the functions $a,b,c$ also yield the result of proposition \ref{ualphac} for a suitable choice of $u_{\alpha,\varpi^{f_0(\alpha)}}$. More precisely, (\ref{scom}) becomes, for every $\alpha\in\Phi$ and for suitable $x$ and $y$:
\begin{equation}\label{scom1}u_{-\alpha,x}u_{\alpha,y}=u_{\alpha,\frac y{1+xy}}h_{\alpha^\vee,\frac 1{1+xy}}u_{-\alpha,\frac x{1+xy}},\end{equation}
just as in the nonsolvable case.

First we check the existence of our constants.

\begin{prop}\label{exct}
Let $\alpha,\beta$ be two elements of $\Phi$ such that $\alpha+\beta\in\Phi$ There exist elements $c_{\alpha,\beta,i,j}$ of $R$, with $i,j$ being positive integers such that $i\alpha+j\beta\in\Phi$, such that we have, for every $\lambda,\mu\in R$:
\[[u_{\alpha,\lambda},u_{\beta,\mu}]=\prod_{i,j}u_{i\alpha+j\beta,c_{\alpha,\beta,i,j}\lambda^i\mu^j}.\]
\end{prop}

\begin{proof}
Note first that since the groups ${\mth{U}}_{i\alpha+j\beta}$ do not necessarily commute, the $c_{\alpha,\beta,i,j}$ may depend on the order in which we take the product, and we have to choose such an order. In the sequel, if $(i,j)\neq(i',j')$, the term in $i\alpha+j\beta$ will then be placed left of the term in $i'\alpha+j'\beta$ if either $j<j'$ or $j=j'$ and $i<i'$.

To simplify the notations and because the remark preceding lemma \ref{fazero} allows us to do so, we will assume that $f_0(\alpha)=f_0(\beta)=0$. We define the $c_{\alpha,\beta,i,j}$ as the ones satisfying the above equality for $\lambda=\mu=1$; we first prove that the above equality holds for any $\lambda$ (still with $\mu=1$). First assume $v(\lambda)=0$. Since $\alpha+\beta$ is a root and $\Phi$ is reduced, $\alpha$ and $\beta$ are linearly independent; by proposition \ref{anniso3}, there exists then an element $h$ of ${\mth{H}}$ commuting with ${\mth{U}}_\beta$ and such that $Ad(h)u_{\alpha,1}=u_{\alpha,\lambda}$. We then have:
\[[u_{\alpha,\lambda},u_{\beta,1}]=Ad(h)[u_{\alpha,1},u_{\beta,1}]=\prod_{i,j}Ad(h)u_{i\alpha+j\beta,c_{\alpha,\beta,i,j}}.\]
By proposition \ref{anniso3} again, with the help of an obvious induction, we see that $Ad(h)u_{\alpha+j\beta,c_{\alpha,\beta,1,j}}=u_{\alpha+j\beta,c_{\alpha,\beta,1,j}\lambda}$ for every $j$ such that $\alpha+j\beta$ is a root. Moreover, for every $j$ such that $\alpha+j\beta$ and $2\alpha+j\beta$ are both roots, applying the previous equality to $(\alpha,\alpha+j\beta)$, we see that $Ad(h)u_{2\alpha+j\beta,c_{\alpha,\beta,2,j}}$ is the projection on ${\mth{U}}_{2\alpha+j\beta}$ of:
\[[Ad(h)u_{\alpha,c_{\alpha,\beta,2,j}c_{\alpha,\alpha+j\beta,1,1}^{-1}},Ad(h)u_{\alpha+j\beta,1}]=[u_{\alpha,c_{\alpha,\beta,2,j}c_{\alpha,\alpha+j\beta,1,1}^{-1}\lambda},u_{\alpha+j\beta,\lambda}],\]
hence $Ad(h)u_{2\alpha+j\beta,c_{\alpha,\beta,2,j}}=u_{2\alpha+j\beta,c_{\alpha,\beta,2,j}\lambda^2}$. We also obtain in a similar way as above that $Ad(h)u_{2\alpha+j\beta,c_{\alpha,\beta,2,j'}}=u_{2\alpha+j'\beta,c_{\alpha,\beta,2,j'}\lambda^2}$ for every $j'>j$ such that $2\alpha+j'\beta$ is a root.

Since for every positive $i,j$ such that $i\alpha+j\beta$ is a root, by \cite[\S 1, proposiition 19]{bou}, either $(i-1)\alpha+j\beta$ or $i\alpha+(j-1)\beta$ is also a root, the case of the remaining pairs $(i,j)$ is proved similarly with the help of an easy induction. Since the set of $\lambda$ of valuation $0$ is dense in $R$, the relations remain true when $v(\lambda)>0$. Finally, the case of any $\mu$ with $\lambda$ arbitrarily fixed is treated symmetrically. $\Box$
\end{proof}

The following corollary will be useful in the next subsections.

\begin{cor}\label{cswap}
We have $u_{\alpha,\lambda}u_{\beta,\mu}=u_{\beta,\mu}u_{\alpha,\lambda}\prod_{i,j}u_{i\alpha+j\beta,(-1)^{i+j}c_{\alpha,\beta,i,j}\lambda^i\mu^j}$.
\end{cor}

\begin{proof}
For every $h,h'$ belonging to any group $H$, we have $hh'=h'h[h^{-1},h'^{-1}]$. The corollary then follows immediately from proposition \ref{exct} applied to $u_{\alpha,-\lambda}$ and $u_{\beta,-\mu}$. $\Box$
\end{proof}

To simplify notations, when $i=j=1$, we can write $c_{\alpha,\beta}$ instead of $c_{\alpha,\beta,1,1}$. Similarly, when $j=1$ and there is no ambiguity, we can drop it (i.e. write $c_{\alpha,\beta,i}$ instead of $c_{\alpha,\beta,i,1}$); we can do the same when $i=1$ and there is no ambiguity. (Note that, as we will see later, such an ambiguity only arises when $\Phi$ is of type $G_2$ and $\alpha$ and $\beta$ are two short roots whose sum is also a short root.)

Note that although the $c_{\alpha,\beta,i,j}$ are elements of $R_{\alpha+\beta}$, they are in general only defined up to some ideal of $R_{\alpha+\beta}$. More precisely, set $\delta=d+1-f_0(-i\alpha-j\beta)-if_0(\alpha)-jf_0(\beta)$; the constant $c_{\alpha,\beta,i,j}$ is an element of the quotient of $R_\alpha$ by its ideal of elements of valuation $\geq \delta$. (If $\delta\leq 0$, that simply means that the term belonging to ${\mth{U}}_{i\alpha+j\beta}$ in the decomposition of $[u_{\alpha,\lambda},u_{\beta,\mu}]$ is trivial.)

Also to simplify notations, by a slight abuse, for such constants we will write $c=c'$ instead of "if $\delta$ (resp. $\delta'$) is the integer associated to $c$ (resp. $c'$) as above, $c$ and $c'$ are equal up to an element of valuation at least $Inf(\delta,\delta')$ of $R$". We will make a similar abuse with expressions of the form $cc'=x$, with $x$ being an element of $R$. Note that this kind of "equality" is not transitive; we thus have to be particularly careful when manipulating it.

We first prove a few preliminary results:

\begin{prop}\label{cswitch}
For every $\alpha,\beta\in\Phi$ such that $\alpha+\beta\in\Phi$ and every $i,j$ such that either $i$ or $j$ is $1$, we have $c_{\beta,\alpha,j,i}=-c_{\alpha,\beta,i,j}$.
\end{prop}

\begin{proof}
Let ${\mth{U}}'$ (resp. ${\mth{U}}''$ be the product of the ${\mth{U}}_{i\alpha+j\beta}$ with $i,j\geq 1$ (resp. $i,j>1$); we see easily using the commutator relations ({\bf (PT4c'')}) that ${\mth{U}}'$ normalizes ${\mth{U}}''$ and that the quotient ${\mth{U}}'/{\mth{U}}''$ is abelian. Moreover, we have $[u_{\alpha,1},u_{\beta,1}]\in \prod_{i,j}u_{i\alpha+j\beta,c_{\alpha,\beta,i,j}}{\mth{U}}''$ and $[u_{\beta,1},u_{\alpha,1}]\in\prod_{i,j}u_{i\alpha+j\beta,c_{\beta,\alpha,j,i}}{\mth{U}}''$, when in both cases $(i,j)$ runs over the set of pairs of positive integers such that at least one of them is $1$; since $[u_{\alpha,1},u_{\beta,1}]$ and $[u_{\beta,1},u_{\alpha,1}]$ are inverse of each other, we must have $u_{i\alpha+j\beta,c_{\alpha,\beta,i,j}}=u_{i\alpha+j\beta,c_{\beta,\alpha,j,i}}^{-1}$ for every $(i,j)$ satisfying that condition, and the result follows. $\Box$
\end{proof}

Note that this is not true when both $i$ and $j$ are strictly greater than $1$; we will give an explicit counterexample later on (proposition \ref{cabij5}).

\begin{lemme}\label{hcomm}
Let $h,h',h''$ be three elements of a group $H$, such that $h'$ and $h''$ commute; we have:
\[[[h,h'],h'']=h[[h^{-1},h''],h']^{-1}h^{-1}.\]
\end{lemme}

\begin{proof}
By rewriting the right-hand side as $h[[h',[h^{-1},h'']]h^{-1}$ and developing both sides we immediately obtain the result. $\Box$
\end{proof}

We now prove a technical lemma:

\begin{lemme}\label{totor}
Let $n$ be a positive integer. The relation on the subsets of $\{1,\dots,n\}$ defined by $I\leq J$ if either $I=J$ or the smallest element of their symmetrical difference belongs to $I$ is a total order.
\end{lemme}

\begin{proof}
It is obvious that the relation is reflexive, antisymmetrical, and that for every $I,J$ we have either $I\leq J$ or $J\leq I$; we only have to check transitivity. Let $I\leq J\leq K$ be three subsets of $\{1,\dots,n\}$; since the result is trivial when two of them are identical, we may assume $I<J<K$. Let $x$ (resp. $y$) be the smallest element of the symmetrical difference between $I$ and $J$ (resp. $J$ and $K$). Then $y$ belongs to $J$ and $x$ does not, hence $x\neq y$. Moreover, we have:
\begin{itemize}
\item if $x<y$, then $x$ does not belong to the symmetrical difference between $J$ and $K$, hence $x\not\in K$;
\item if $x>y$, then $y$ does not belong to the symmetrical difference between $I$ and $J$, hence $y\in I$.
\end{itemize}
In both cases, $Inf(x,y)$ belongs to $I$ and not to $K$. On the other hand, for every $z<Inf(x,y)$, since $z$ does not belong th any of the above two symmetrical differences, it must belong either to all three of $I,J,K$ or to none of them. Hence $I<K$ and the lemma is proved. $\Box$
\end{proof}

Now we use this lemma to prove the following proposition, which will be useful in the sequel.

\begin{prop}\label{comimbr}
Let $a_1,\dots,a_n,b$ be elements of a group $G$ such that all the $a_i$ commute. We have:
\[[a_1,[a_2,\dots,[a_n,b]\dots]]=\prod_{I\subset\{1,\dots,n\}}[\prod_{i\in I}a_i,b]^{(-1)^{n-\#(I)}}\]
where the product is taken with the following order on the subsets of $\{1,\dots,n\}$: for every $I\neq J$, $I$ precedes $J$ if the smallest element of their symmetrical difference belongs to $I$.
\end{prop}

\begin{proof}
We prove the proposition by induction on $n$, the case $n=1$ being trivial. Assume first $n=2$; we have:
\[[a_1,[a_2,b]]=a_1a_2ba_2^{-1}b^{-1}a_1^{-1}ba_2b^{-1}a_2^{-1}\]
\[=(a_1a_2ba_2^{-1}a_1^{-1}b)(ba_1b^{-1}a_1^{-1})(ba_2b^{-1}a_2^{-1})\]
\[=[a_1a_2,b][a_1,b]^{-1}[a_2,b]^{-1}.\]

Now assume $n>2$. Applying the case $n=2$ to $a_1$, $a_2$ and $[a_3,\dots,[a_n,b]\dots,]$, we obtain:
\[[a_1,[a_2,\dots,[a_n,b]\dots]]\]
\begin{equation}\label{tterm}=[a_1a_2,[a_3,\dots,a_n,b]\dots]][a_1,[a_3,\dots[a_n,b]\dots]]^{-1}[a_2,[a_3,\dots,[a_n,b]\dots]]^{-1};\end{equation}
Applying the induction hypothesis to each one of the three terms of the right-hand side, we obtain:
\[[a_1a_2,\dots,[a_n,b]\dots]]=\prod_{I\subset\{2,\dots,n\}}[\prod_{i\in I}c_i,b]^{(-1)^{n-1-\#(I)}}\]
where $c_2=a_1a_2$ and $c_i=a_i$ for every $i>2$;
\[[a_1,[a_3,\dots,[a_n,b]\dots]]=\prod_{I\subset\{1,3,\dots,n\}}[\prod_{i\in I}a_i,b]^{(-1)^{n-1-\#(I)}};\]
\[[a_2,[a_3,\dots,[a_n,b]\dots]]=\prod_{I\subset\{2,\dots,n\}}[\prod_{i\in I}a_i,b]^{(-1)^{n-1-\#(I)}}\]

We remark that in the first two products, the terms coming from $[a_3,\dots,[a_n,b]\dots]$ care identical, all on the right side of the expression and in the same order, which implies that they cancel out in (\ref{tterm}); since by lemma \ref{totor} the remaining terms can posiblly be in only one order, the result follows. $\Box$
\end{proof}

\begin{prop}\label{com3}
Let $\alpha,\beta,\gamma$ be three elements of $\Phi$ such that $\alpha+\beta$, $\alpha+\gamma$ and $\alpha+\beta+\gamma$ are roots. Assume there exists a subset of positive roots of $\Phi$ containing all three of them, and assume also that one of the following conditions is satisfied:
\begin{itemize}
\item $f_0(\alpha+\beta+\gamma)=f_0(\alpha)+f_0(\beta)+f_0(\gamma)$;
\item $f_0(\alpha+\beta+\gamma)=f_0(\alpha)+f_0(\beta)+f_0(\gamma)-1$ and $\alpha+\beta+\gamma\in\Psi$.
\end{itemize}
Then we have:
\[c_{\alpha,\beta}c_{\gamma,\alpha+\beta}=c_{\gamma,\alpha}c_{\alpha+\gamma,\beta}.\]
\end{prop}

\begin{proof}
We will first check that the only pair $(j,k)$ of positive integers such that $j\alpha+k\beta=\alpha+\beta$ is $(1,1)$. Assume $(j,k)\neq(1,1)$; since both roots are contained in some given set of positive roots and $(j-1)\alpha)+(k-1)\beta=0$, Since we assumes $j$ and $k$ are both positive integers, we obtain a contradiction. Hence we must have $j=k=1$.

By applying lemma \ref{hcomm} to $u_{\alpha,\lambda},u_{\beta,1},u_{\gamma,1}$ for every $\lambda\in R_a$ and considering only the terms in $u_{\alpha+\beta+\gamma}$ in both sides of the expression, we obtain an equality of the form:
\[c_{\alpha,\beta}c_{\gamma,\alpha+\beta}\lambda+\lambda^2P(\lambda)=c_{\gamma,\alpha}c_{\alpha+\gamma,\beta}\lambda+\lambda^2Q(\lambda),\]
for every $\lambda$, where $P$ and $Q$ are polynomials. The result of the proposition follows immediately, the condition on $f_0$ ensuring that the result is actually true and not true only up to some nontrivial ideal. $\Box$
\end{proof}

The following lemma, whose proof is straightforward, will also be useful in the sequel:

\begin{lemme}\label{prodcom}
Let $h,h',h''$ be three elements of a group $H$, such that $h'$ and $h''$ commute; then $[hh',h'']=[h,h'']$.
\end{lemme}

The next three subsections will be devoted to establishing relations between our structural constants by examining what happens in subsystems of rank $2$ of $\Phi$. Let $\alpha,\beta$ be two linearly independent elements of $\Phi$ such that $\alpha+\beta$ is a root and $\alpha-\beta$ is not; $\alpha$ and $\beta$ then generate a subsystem of rank $2$ of $\Phi$. By proposition \ref{crs} and lemma \ref{fazero}, we can always assume that $f_0$ is such that $f_0(\alpha)=f_0(\beta)=0$. We will also assume that when $\alpha$ and $\beta$ are of different lengths, $\alpha$ is the shortest one.

\subsection{The case $A_2$}

Assume first that $\alpha$ and $\beta$ generate a subsystem of type $A_2$ of $\Phi$. Then the roots it contains are $\pm\alpha$, $\pm\beta$ and $\pm(\alpha+\beta)$, which are all of the same length.

\begin{prop}\label{com2}
We have:
\[c_{\alpha,\beta}c_{\alpha+\beta,-\beta}=1.\]
Moreover, for every $\lambda\in R^*$ and every $\mu\in R_\alpha$, we have $Ad(h_{\beta^\vee,\lambda})(u_{\alpha^\vee,\mu})=u_{\alpha^\vee,\frac\mu\lambda}$.
\end{prop}

\begin{proof}
For every $\lambda\neq 1$ such that $u_{-\beta,\lambda}$ exists, we have:
\[[[u_{\alpha,1},u_{\beta,1}],u_{-\beta,\lambda}]=[u_{\alpha+\beta,c_{\alpha,\beta}},u_{-\beta,\lambda}]=u_{\alpha,c_{\alpha,\beta}c_{\alpha+\beta,-\beta}\lambda},\]
hence:
\[[u_{\alpha,1},u_{\beta,1}]u_{-\beta,\lambda}=u_{\alpha,c_{\alpha,\beta}c_{\alpha+\beta,-\beta}\lambda}u_{-\beta,\lambda}u_{\alpha+\beta,c_{\alpha,\beta}}.\]
\[=u_{-\beta,\lambda}u_{\alpha,c_{\alpha,\beta}c_{\alpha+\beta,-\beta}\lambda}u_{\alpha+\beta,c_{\alpha,\beta}}.\]
On the other hand, we have:
\[[u_{\alpha,1},u_{\beta,1}]u_{-\beta,\lambda}=u_{\alpha;1}u_{\beta,1}u_{\alpha,-1}u_{\beta,-1}u_{-\beta,\lambda}\]
\[=u_{\alpha;1}u_{\beta,1}u_{\alpha,-1}u_{-\beta,\frac\lambda{1-\lambda}}h_{\beta^\vee,1-\lambda}u_{\beta,\frac 1{\lambda-1}}\]
\[=u_{\alpha;1}u_{\beta,1}u_{-\beta,\frac\lambda{1-\lambda}}u_{\alpha,-1}h_{\beta^\vee,1-\lambda}u_{\beta,\frac 1{\lambda-1}}\]
\begin{equation}\label{ca2bl}=u_{\alpha;1}u_{-\beta,\lambda}h_{\beta^\vee,\frac 1{1-\lambda}}u_{\beta,1-\lambda}u_{\alpha,-1}h_{\beta^\vee,1-\lambda}u_{\beta,\frac 1{\lambda-1}}\end{equation}
\[=u_{-\beta,\lambda}u_{\alpha;1}u_{\beta,\frac 1{1-\lambda}}u_{\alpha,\mu}u_{\beta,\frac 1{\lambda-1}},\]
where $u_{\alpha,\mu}=Ad(h_{\beta^\vee,\frac 1{1-\lambda}})u_{\alpha,-1}$. From the previous equality we see that we must have for every $\lambda$:
\[\mu=c_{\alpha,\beta}c_{\alpha+\beta,-\beta}\lambda-1.\]
Since the right-hand side must be a monomial in $\lambda-1$, we obtain $c_{\alpha,\beta}c_{\alpha+\beta,-\beta}=1$, as desired.

Now we prove the second assertion. We see from above that for every $\lambda\in {\mth{H}}_\beta$, we have $Ad(h_{\beta^\vee,\lambda})u_{\alpha,-1}=u_{\alpha,-\frac 1\lambda}$ for every $\lambda\in R^*$ of valuation at least $f_0(-\beta)$. Since $f_0(-\beta)$ is either $0$ or $1$, tha previous equality holds as soon as $h_{\beta^\vee,\lambda}$ belongs to ${\mth{H}}_{\beta,1}$. By definition of a root subgroup and since $<\alpha,\beta^\vee>=-1$ it also holds when $h_{\beta^\vee,\lambda}\in{\mth{T}}_\beta$, hence it holds for any $\lambda$. Since by the properties of $R_\alpha$, for every element $h$ of ${\mth{H}}$, $Ad(h)$ acts on $R_\alpha$ by multiplication by some given element of $R_\alpha$, the end of the proof is straightforward. $\Box$
\end{proof}

\begin{prop}\label{prodh}
For every $\lambda\in R^*$, we have $h_{\alpha^\vee+\beta^\vee,\lambda}=h_{\alpha^\vee,\lambda}h_{\beta^\vee,\lambda}$.
\end{prop}

\begin{proof}
Since $\alpha$ and $\beta$ generate a subsystem of type $A_2$ of $\Phi$, we have $(\alpha+\beta)^\vee=\alpha^\vee+\beta^\vee$. By proposition \ref{basprod}, we know that there exist $\mu,\nu\in R^*$ such that $h_{\alpha^\vee+\beta^\vee,\lambda}=h_{\alpha^\vee,\mu}h_{\beta^\vee,\nu}$. By considering the actions of these elements respectively on ${\mth{U}}_\alpha$, ${\mth{U}}_\beta$ and ${\mth{U}}_{\alpha+\beta}$, with the help of the previous proposition, we obtain:
\[\lambda=\mu^2\nu^{-1};\;\;\lambda=\mu^{-1}\nu^2;\;\;\lambda^2=\mu\nu.,\]
from which we deduce:
\[\lambda\mu^{-1}=\mu\nu^{-1}=\nu\lambda^{-1}.\]
Since the product of these three terms is $1$, they are all equal to some given third root of $1$. On the other hand, $\lambda\mapsto\lambda\mu^{-1}$ is a continuous group automorphism of $R^*$, which must then be trivial. The result follows. $\Box$
\end{proof}

Now we come to the point where we actually have to be careful about the rings in which the constants live. The following proposition, for example, can be easily deduced from proposition \ref{com2} when the equalities are compatible: we obtain $c_{\alpha,\beta}c_{\alpha+\beta,-\beta}=1$ and (by proposition \ref{com2} applied to $-\beta$ and $\alpha+\beta$) $c_{\alpha+\beta,-\beta}c_{-\alpha-\beta,\alpha}=1$, hence $c_{\alpha,\beta}=c_{-\alpha-\beta,\alpha}$, and the equality with $c_{\beta,-\alpha-\beta}$ is obtained the same way. It can be easily checked (details are left to the reader) that the compatibility condition is satisfied when at least one of $\alpha,\beta,\alpha+\beta$ belongs to $\Psi$, or when none does and $f_0(\alpha+\beta)<f_0(\alpha)+f_0(\beta)$; on the other hand, in the remaining case, that is when none of those three roots lies in $\Psi$ and $f_0(\alpha+\beta)=f_0(\alpha)+f_0(\beta)$, $c_{\alpha,\beta}$, $c_{\beta,-\alpha-\beta}$ and $c_{-\alpha-\beta,\alpha}$ are defined up to an element of valuation $d$ of $R$, but for example, $c_{\alpha+\beta,-\beta}$ os defined only up to an element of valuation $d-1$ of $R$. and the above reasoning only yields the relations below up to an element of valuation $d-1$ of $R$. We thus need to prove them directly.

Note though that for every $\alpha,\beta$, $c_{\alpha,\beta}$ and $c_{\beta,\alpha}$ always live in the same ring. Hence we are at least allowed to apply lemma \ref{cswitch} freely.

\begin{prop}\label{cycla2}
We have:
\[c_{\alpha,\beta}=c_{\beta,-\alpha-\beta}=c_{-\alpha-\beta,\alpha}.\]
\end{prop}

\begin{proof}
For every $\lambda\neq 1$ such that $u_{\-\alpha-\beta,\lambda}$ exists, we have:
\[[u_{\alpha,1},u_{\beta,1}]u_{-\alpha-\beta,\lambda}=u_{\alpha+\beta,c_{\alpha,\beta}}u_{-\alpha-\beta,\lambda}\]
\[=u_{-\alpha-\beta,\frac\lambda{1+c_{\alpha,\beta}\lambda}}h_{\alpha^\vee+\beta^\vee,1+c_{\alpha,\beta}\lambda}u_{\alpha+\beta,\frac{c_{\alpha,\beta}}{1+c_{\alpha,\beta}\lambda}}.\]
On the other hand, we have, using the fact that $c_{-\alpha-\beta,\alpha}=-c_{\alpha,-\alpha-\beta}$:
\[[u_{\alpha,1},u_{\beta,1}]u_{-\alpha-\beta,\lambda}=u_{\alpha,1}u_{\beta,1}u_{\alpha,-1}u_{\beta,-1}u_{-\alpha-\beta,\lambda}\]
\[=u_{-\alpha-\beta,\lambda}u_{-\beta,c_{\alpha,-\alpha-\beta}\lambda}u_{\alpha,1}u_{-\alpha,c_{\beta,-\alpha-\beta}\lambda}u_{\beta,1}u_{-\beta,-c_{\alpha,-\alpha-\beta}\lambda}u_{\alpha,-1}u_{-\alpha,-c_{\beta,-\alpha-\beta}\lambda}u_{\beta,1}\]
\[=u_{-\alpha-\beta,\lambda}u_{-\beta,c_{\alpha,-\alpha-\beta}\lambda}u_{-\alpha,\frac{c_{\beta,-\alpha-\beta}\lambda}{1+c_{\beta,-\alpha-\beta}\lambda}}h_{\alpha^\vee,1+c_{\beta,-\alpha-\beta}\lambda}u_{\alpha,\frac 1{1+c_{\beta,-\alpha-\beta}\lambda}}u_{-\beta,\frac{c_{-\alpha-\beta,\alpha}\lambda}{1+c_{-\alpha-\beta,\alpha}\lambda}}\]
\[*h_{\beta^\vee, 1+c_{-\alpha-\beta,\alpha}\lambda}]u_{\beta,\frac 1{1+c_{-\alpha-\beta,\alpha}\lambda}}u_{-\alpha,-\frac{c_{\beta,-\alpha-\beta}\lambda}{1+c_{\beta,-\alpha-\beta}\lambda}}h_{\alpha^\vee,1+c_{\beta,-\alpha-\beta}\lambda}u_{\alpha,-\frac 1{1+c_{\beta,-\alpha-\beta}\lambda}}u_{\beta,1}.\]
Since ${\mth{U}}_{-\alpha}$ and ${\mth{U}}_\beta$ (resp. ${\mth{U}}_{-\beta}$ and ${\mth{U}}_\alpha$) commute, and since we have:
\[u_{\alpha,\frac 1{1+c_{\beta,-\alpha-\beta}\lambda}}h_{\beta^\vee, 1+c_{-\alpha-\beta,\alpha}}u_{-\alpha,-\frac{c_{\beta,-\alpha-\beta}\lambda}{1+c_{\beta,-\alpha-\beta}\lambda}}=u_{\alpha,\frac 1{1+c_{\beta,-\alpha-\beta}\lambda}}u_{-\alpha,-c_{\beta,-\alpha-\beta}\lambda}h_{\beta^\vee, 1+c_{-\alpha-\beta,\alpha}}\]
and:
\[u_{\alpha,\frac 1{1+c_{\beta,-\alpha-\beta}\lambda}}u_{-\alpha,-c_{\beta,-\alpha-\beta}\lambda}=u_{-\alpha,\frac{-c_{\beta,-\alpha-\beta}\lambda}{1+c_{\beta,-\alpha-\beta}\lambda}}h_{\alpha^\vee,\frac 1{1+c_{\beta,-\alpha-\beta}\lambda}}u_{\alpha,1},\]
we obtain that $[u_{\alpha,1},u_{\beta,1}]u_{-\alpha-\beta,\lambda}$ is equal to:
\[u_{-\alpha-\beta,\lambda}u_{-\beta,c_{\alpha,-\alpha-\beta}\lambda}u_{-\alpha,\frac{c_{\beta,-\alpha-\beta}\lambda}{1+c_{\beta,-\alpha-\beta}\lambda}}h_{\alpha^\vee,1+c_{\beta,-\alpha-\beta}\lambda}u_{-\beta,\frac{c_{-\alpha-\beta,\alpha}\lambda}{1+c_{-\alpha-\beta,\alpha}\lambda}}u_{-\alpha,\frac{-c_{\beta,-\alpha-\beta}\lambda}{1+c_{\beta,-\alpha-\beta}\lambda}}\]
\[*h_{\alpha^\vee,\frac 1{1+c_{\beta,-\alpha-\beta}\lambda}}u_{\alpha,1}h_{\beta^\vee, 1+c_{-\alpha-\beta,\alpha}}u_{\beta,\frac 1{1+c_{-\alpha-\beta,\alpha}\lambda}}h_{\alpha^\vee,1+c_{\beta,-\alpha-\beta}\lambda}u_{\alpha,-\frac 1{1+c_{\beta,-\alpha-\beta}\lambda}}u_{\beta,1}.\]
Using the fact that ${\mth{H}}$ normalizes every ${\mth{U}}_\alpha$, we are finally left with an expression of the form:
\[[u_{\alpha,1},u_{\beta,1}]u_{-\alpha-\beta,\lambda}\in{\mth{U}}^-h_{\alpha^\vee,\frac{(1+c_{\beta,-\alpha-\beta}\lambda)^2}{1+c_{-\alpha-\beta,-\alpha}\lambda}}h_{\beta^\vee, 1+c_{-\alpha-\beta,\alpha}\lambda}{\mth{U}}^+\]
where ${\mth{U}}^-$ (resp. ${\mth{U}}^+$) is the product of the ${\mth{U}}_\gamma$, where $\gamma$ runs over the set of roots which are linear combinations with negative or zero (resp. positive or zero) integer coefficients of $\alpha$ and $\beta$; proposition \ref{prodh} now yields for every $\lambda$:
\[1+c_{\alpha,\beta}\lambda=\frac{(1+c_{\beta,-\alpha-\beta}\lambda)^2}{1+c_{-\alpha-\beta,\alpha}\lambda}=1+c_{-\alpha-\beta,\alpha}\lambda,\]
from which we deduce, by comparing the terms in $\lambda$ in the above four polynomials:
\begin{equation}\label{blu}c_{\alpha,\beta}=2c_{\beta,-\alpha-\beta}-c_{-\alpha-\beta,\alpha}=c_{-\alpha-\beta,\alpha}.\end{equation}
The proposition follows immediately. $\Box$
\end{proof}

\begin{cor}\label{paba2}
We have:
\[c_{\alpha,\beta}c_{-\alpha,-\beta}=-1.\]
\end{cor}

\begin{proof}
We will in fact prove that $c_{\alpha,\beta}c_{-\beta,-\alpha}=1$, which is equivalent by lemma \ref{cswitch}. By propositions \ref{com2} and \ref{cycla2} we have $c_{\alpha,\beta}c_{\alpha+\beta,-\beta}=1$ and $c_{\alpha+\beta,-\beta}=c_{-\beta,-\alpha}$. On the other hand, by proposition \ref{f0ab3} we have either $f_0(\alpha+\beta)=f_0(\alpha)+f_0(\beta)$ or $f_0(\alpha)=f_0(\alpha+\beta)+f_0(-\beta)$, and we can easily deduce from this  that $c_{\alpha+\beta,-\beta}$ lives in the same ring as either $c_{\alpha,\beta}$ or $c_{-\beta,-\alpha}$. The result follows. $\Box$
\end{proof}

\subsection{The case $B_2$}

Assume now $\alpha$ and $\beta$ generate a subsystem of $\Phi$ of type $B_2$. The elements of $\Phi$ are then $\pm\alpha$ (short), $\pm\beta$ (long), $\pm(\alpha+\beta)$ (short) and $\pm(2\alpha+\beta)$ (long). Moreover we have $\alpha^\vee+\beta^\vee=(2\alpha+\beta)^\vee$ and $2\alpha^\vee+\beta^\vee=(\alpha+\beta)^\vee$.

\begin{prop}\label{cabij2}
We have:
\[c_{\alpha,\beta,2}=\frac 12c_{\alpha,\beta}c_{\alpha,\alpha+\beta}.\]
\end{prop}

\begin{proof}
According to proposition \ref{comimbr}, we have:
\[[u_{\alpha,1},[u_{\alpha,1},u_{\beta;1}]]=[u_{\alpha,2},u_{\beta,1}][u_{\alpha_1},u_{\beta,1}]^{-2}.\]
By evaluating both sides, we obtain:
\[u_{2\alpha+\beta,c_{\alpha,\beta}c_{\alpha,\alpha+\beta}}=u_{\alpha+\beta,2c_{\alpha,\beta}}u_{2\alpha+\beta,4c_{\alpha,\beta,2}}(u_{\alpha+\beta,c_{\alpha,\beta}}u_{2\alpha+\beta,c_{\alpha,\beta,2}})^{-2},\]
Since ${\mth{U}}_{\alpha+\beta}$ and ${\mth{U}}_{2\alpha+\beta}$ commute,  considering the terms in ${\mth{U}}_{2\alpha+\beta}$, we obtain:
\[c_{\alpha,\beta}c_{\alpha,\alpha+\beta}=2c_{\alpha,\beta,2},\]
which proves the proposition. $\Box$
\end{proof}

\begin{prop}\label{com1b2}
We have:
\[c_{\alpha,\beta}c_{\alpha+\beta,-\beta}=1.\]
Moreover, for every $\lambda\in R^*$ and every $\mu\in R_\alpha$, we have $Ad(h_{\beta^\vee,\lambda})u_{\alpha,\mu}=u_{\alpha,\frac\mu\lambda}$.
\end{prop}

\begin{proof}
For every suitable $\lambda$, we have:
\[[[u_{\alpha,1},u_{\beta,1}],u_{-\beta,\lambda}]=[u_{\alpha+\beta,c_{\alpha,\beta}}u_{2\alpha+\beta,c_{\alpha,\beta,2}},u_{-\beta,\lambda}].\]
Since $(2\alpha+\beta)-\beta=2\alpha$ is not a root, ${\mth{U}}_{2\alpha+\beta}$ and ${\mth{U}}_{-\beta}$ commute and we obtain, using lemma \ref{prodcom}:
\[[[u_{\alpha,1},u_{\beta,1}],u_{-\beta,\lambda}]=[u_{\alpha+\beta,c_{\alpha,\beta}},u_{-\beta,\lambda}]=u_{\alpha,c_{\alpha,\beta}c_{\alpha+\beta,-\beta}\lambda}.\]
The rest of the proof is exactly similar to the proof of proposition \ref{com2}. (Note also that the second assertion has already been proved in the course of the proof of proposition \ref{anniso3}.) $\Box$
\end{proof}

\begin{prop}\label{com2b2}
We have:
\[c_{\beta,\alpha}c_{\alpha+\beta,-\alpha}=2.\]
Moreover, for every $\lambda\in R^*$ and every $\mu\in R_\beta$, we have $Ad(h_{\alpha^\vee,\lambda})u_{\beta,\mu}=u_{\beta,\frac\mu{\lambda^2}}$.
\end{prop}

\begin{proof}
For every suitable $\lambda$, we have:
\[[[u_{\beta,1},u_{\alpha,1}],u_{-\alpha,\lambda}]=[u_{\alpha+\beta,c_{\beta,\alpha}}u_{2\alpha+\beta,c_{\beta,\alpha,2}},u_{-\alpha,\lambda}].\]
Using proposition \ref{comimbr}, we see that the right-hand side is equal to:
\[[u_{\alpha+\beta,c_{\beta,\alpha}},[u_{2\alpha+\beta,c_{\beta,\alpha,2}},u_{-\alpha,\lambda}]][u_{2\alpha+\beta,c_{\beta,\alpha,2}},u_{-\alpha,\lambda}][u_{\alpha+\beta,c_{\beta,\alpha}},u_{-\alpha,\lambda}].\]
It is easy to see by developing it that the first of these three commutators is trivial. We thus have:
\[[[u_{\beta,1},u_{\alpha,1}],u_{-\alpha,\lambda}]=u_{\alpha+\beta,c_{\beta,\alpha,2}c_{2\alpha+\beta,-\alpha}\lambda}u_{\beta,c_{\beta,\alpha,2}c_{2\alpha+\beta,-\alpha,2}\lambda^2+c_{\beta,\alpha}c_{\alpha+\beta,-\alpha}\lambda}.\]
On the other hand, we have:
\[u_{\beta,1}u_{\alpha,1}u_{\beta,-1}u_{\alpha,-1}u_{-\alpha,\lambda}=u_{\beta,1}u_{\alpha,1}u_{\beta,-1}u_{-\alpha,\frac\lambda{1-\lambda}}h_{\alpha^\vee,1-\lambda}u_{\alpha,\frac 1{\lambda-1}}\]
\[=u_{\beta,1}u_{-\alpha,\lambda}h_{\alpha^\vee,\frac 1{1-\lambda}}u_{\alpha,1-\lambda}u_{\beta,-1}h_{\alpha^\vee,1-\lambda}u_{\alpha,\frac 1{\lambda-1}}\]
By a similar reasoning as in the proof of proposition \ref{com2}, we must have:
\[Ad(h_{\alpha^\vee,\frac 1{1-\lambda}})u_{\beta,1}=u_{\beta,-1+c_{\beta,\alpha,2}c_{2\alpha+\beta,-\alpha,2}\lambda^2+c_{\beta,\alpha}c_{\alpha+\beta,-\alpha}\lambda}.\]
Since the coefficient of the right-hand side is a monomial in $1-\lambda$, it must be $(1-\lambda)^2$ and the first assertion of the proposition follows immediately. The second one is proved the same way as in proposition \ref{com2}. $\Box$
\end{proof}

\begin{prop}\label{prodh2}
For every $\lambda\in R^*$, we have $h_{\alpha^\vee+\beta^\vee,\lambda}=h_{\alpha(^\vee,\lambda}h_{\beta^\vee,\lambda}$.
\end{prop}

\begin{proof}
SInce $\alpha$ and $\beta$ generate a subsystem of $\Phi$ of type $B_2$, we have $\alpha^\vee+\beta^\vee=(2\alpha+\beta)^\vee$. For every $\lambda\in R^*$, let $\mu,\nu\in R^*$ be such that $h_{\alpha^\vee+\beta^\vee,\lambda}=h_{\alpha^\vee,\mu}h_{\beta^\vee,\nu}$; by examining the actions of both members on respectively ${\mth{U}}_\alpha$, ${\mth{U}}_\beta$ and ${\mth{U}}_{2\alpha+\beta}$, we obtain, using propositions \ref{com1b2} and \ref{com2b2}:
\[\lambda=\mu^2\nu^{-1};\;\;1=\mu^{-2}\nu^2;\;\;\lambda^2=\mu^2.\]
We deduce from the first two relations that $\lambda=\nu$. By the last one, $\lambda=\pm\mu$, and by continuity of $\lambda\mapsto\mu$ we must have $\lambda=\mu$. Hence the result. $\Box$
\end{proof}

\begin{prop}\label{act1b2}
The group ${\mth{H}}_\alpha$ acts trivially on ${\mth{U}}_{\alpha+\beta}$.
\end{prop}

\begin{proof}
Since $\alpha$ and $\alpha+\beta$ are orthogonal short roots, we deduce from lemma \ref{myz} that every element of ${\mth{H}}$ of the form $h_{-\beta^\vee,\lambda}h_{(2\alpha+\beta)^\vee,\lambda}$ acts trivially on ${\mth{U}}_{\alpha+\beta}$. We thus only have to prove that every element of ${\mth{H}}_\alpha$ is of that form.

By the previous proposition, we have for every $\lambda\in R^*$:
\[h_{-\beta^\vee,\lambda}=h_{\alpha^\vee,\lambda}h_{(-\alpha-\beta)^\vee,\lambda};\]
\[h_{(2\alpha+\beta)^\vee,\lambda}=h_{\alpha^\vee,\lambda}h_{(\alpha+\beta)^\vee,\lambda},\]
hence:
\[h_{-\beta^\vee,\lambda}h_{(2\alpha+\beta)^\vee,\lambda}=h_{\alpha^\vee,\lambda^2}.\]
Since by lemma \ref{hsq}, every element of ${\mth{H}}_\alpha$ (resp. ${\mth{H}}_{2\alpha+\beta}$ is a square in ${\mth{H}}_\alpha$ (resp. ${\mth{H}}_{2\alpha+\beta}$, the corollary follows. $\Box$
\end{proof}

\begin{prop}\label{cyclb2}
We have:
\[c_{\alpha,\beta}=c_{\beta,-\alpha-\beta}=\frac 12c_{-\alpha-\beta,\alpha}.\]
\end{prop}

\begin{proof}
First we observe that if $f_0(-\alpha-\beta)=0$, we have:
\[0=-f_0(-\alpha-\beta)\leq f_0(\alpha+\beta)\leq f_0(\alpha)+f_0(\beta)=0,\]
hence $\alpha+\beta\in\Psi$. By a similar reasoning, $\beta$ and $-\alpha-\beta$ also lie in $\Psi$ as well as their opposites, hence $c_{\alpha,\beta}$, $c_{\beta,-\alpha-\beta}$ and $c_{-\alpha-\beta,\alpha}$ are all elements of $R$, and we then only have to apply the previous propositions to get the result: we have, using successively propositions \ref{com1b2}, \ref{com2b2}, \ref{com2b2} again and \ref{com1b2} again:
\[c_{\alpha,\beta}=c_{\alpha+\beta,-\beta}^{-1}=\frac 12c_{-\alpha-\beta,\alpha}=c_{-\beta,-\alpha}^{-1}=c_{\beta,-\alpha-\beta}.\]

Assume now $f_0(-\alpha-\beta)=1$; $R_{-\alpha-\beta}$ will then be assimilated to the maximal ideal of $R$. In particular, for every $\lambda\in R_{-\alpha-\beta}$, we have $\lambda^{d+1}=0$, and every rational function on $R_{-\alpha-\beta}$ of the form $\frac 1{1+\lambda P(\lambda)}$ is then also identically equal to a polynomial in $\lambda$, since we have for every $\lambda$:
\[\frac 1{1+\lambda P(\lambda)}=1-\lambda P(\lambda)+\dots+(-1)^d\lambda^dP(\lambda)^d.\]

For every suitable $\lambda$, we have:
\[[u_{\alpha,1},u_{\beta,1}]u_{-\alpha-\beta,\lambda}=u_{\alpha+\beta,c_{\alpha,\beta}}u_{2\alpha+\beta,-c_{\alpha,\beta,2}}u_{-\alpha-\beta,\lambda}\]
\[=u_{-\alpha-\beta,\frac\lambda{1+c_{\alpha,\beta}\lambda}}h_{(\alpha+\beta)^\vee,1+c_{\alpha+\beta}\lambda}u_{\alpha+\beta,\frac{c_{\alpha,\beta}}{1+c_{\alpha,\beta}\lambda}}\]
\[*u_{2\alpha+\beta,c_{-\alpha,\beta,2}}u_{\alpha,-c_{2\alpha+\beta,-\alpha-\beta}c_{\alpha,\beta,2}\lambda}u_{-\beta,-c_{2\alpha+\beta,-\alpha-\beta}c_{\alpha,\beta,2}\lambda^2}.\]
On the other hand, we have:
\[[u_{\alpha,1},u_{\beta,1}]u_{-\alpha-\beta,\lambda}=u_{\alpha,1}u_{\beta,1}u_{\alpha,-1}u_{\beta,-1}u_{-\alpha-\beta,\lambda}\]
\[=u_{-\alpha-\beta,\lambda}u_{-\beta,c_{\alpha,-\alpha-\beta}\lambda}u_{\alpha,1}u_{-\alpha,c_{\beta,-\alpha-\beta}\lambda}u_{\beta,1}u_{-\beta,-c_{\alpha,-\alpha-\beta}\lambda}u_{\alpha,-1}u_{-\alpha,-c_{\beta,-\alpha-\beta}\lambda}u_{\beta,1}u',\]
where $u'$ is a product of terms of the form $u_{\gamma,\lambda^2P(\lambda)}$, $P(\lambda)$ being a polynomial in $\lambda$. Moreover, it is easy to check that for every root $\delta$ and every $\mu$, $[u',u_{\delta,\mu}]$ is also a product of terms of the same form; by a similar reasoning as in proposition \ref{cycla2}, we are finally left with an expression of the form:
\[[u_{\alpha,1},u_{\beta,1}]u_{-\alpha-\beta,\lambda}\in{\mth{U}}^-h_{\alpha^\vee,\frac{(1+c_{\beta,-\alpha-\beta}\lambda)^3}{1+c_{-\alpha-\beta,\alpha}\lambda}}h_{\beta^\vee, 1+c_{-\alpha-\beta,\alpha}\lambda}u''{\mth{U}}^+,\]
with ${\mth{U}}^-$ and ${\mth{U}}^+$ being defined in a similar way as in that proposition and $u''$ being of the same form as $u'$. Using proposition \ref{prodh2} and the fact that $(\alpha+\beta)^\vee=\alpha^\vee+2\beta^\vee$, we obtain:
\[c_{\alpha,\beta}=3c_{\beta,-\alpha-\beta}-c_{-\alpha-\beta,\alpha}=\frac 12c_{-\alpha-\beta,\alpha},\]
The result follows immediately. $\Box$
\end{proof}

\begin{cor}\label{pabb2}
We have:
\[c_{\alpha,\beta}c_{-\alpha,-\beta}=-1,\]
\[c_{\alpha,\alpha+\beta}c_{-\alpha,-\alpha-\beta}=-4,\]
\end{cor}

\begin{proof}
Similarly to the case $A_2$, we deduce those equalities from propositions \ref{com1b2}, \ref{com2b2} and \ref{cyclb2}. $\Box$
\end{proof}

\subsection{The case $G_2$}

Assume now $\alpha$ and $\beta$ generate a subsystem of $\Phi$ of type $G_2$ (which implies in fact that $\Phi$ itself is of type $G_2$). The elements of $\Phi$ are then $\pm\alpha$ (short), $\pm\beta$ (long), $\pm(\alpha+\beta)$ (short) $\pm(2\alpha+\beta)$ (short), $\pm(3\alpha+\beta)$ (long) and $\pm(3\alpha+2\beta)$ (long); moreover, we have:
\[(\alpha+\beta)^\vee=\alpha^\vee+3\beta^\vee;\]
\[(2\alpha+\beta)^\vee=2\alpha^\vee+3\beta^\vee;\]
\[(3\alpha+\beta)^\vee=\alpha^\vee+\beta^\vee;\]
\[(3\alpha+2\beta)^\vee=\alpha^\vee+2\beta^\vee.\]

Remember that in the case $G_2$, we are assuming that the characteristic of $k$ is neither $2$ nor $3$.

\begin{prop}\label{cabij3}
We have:
\[c_{\alpha,\alpha+\beta,2,1}=\frac 12c_{\alpha,\alpha+\beta}c_{\alpha,2\alpha+\beta};\]
\[c_{\alpha,\alpha+\beta,1,2}=\frac 12c_{\alpha,\alpha+\beta}c_{2\alpha+\beta,\alpha+\beta}.\]
\end{prop}

\begin{proof}
According to proposition \ref{comimbr}, we have:
\[[u_{\alpha,1},[u_{\alpha,1},u_{\alpha+\beta,1}]]=[u_{\alpha,2},u_{\alpha+\beta,1}][u_{\alpha,1},u_{\alpha+\beta,1}]^{-2}.\]
The left-hand side is equal to:
\[[u_{\alpha,1},u_{2\alpha+\beta,c_{\alpha,\alpha+\beta}}u_{3\alpha+\beta,c_{\alpha,\alpha+\beta,2,1}}u_{3\alpha+2\beta,c_{\alpha,\alpha+\beta,1,2}}]\]
\[=[u_{\alpha,1},u_{2\alpha+\beta,c_{\alpha,\alpha+\beta}}]=u_{3\alpha+\beta,c_{\alpha,\alpha+\beta}c_{\alpha,2\alpha+\beta}}.\]
The right-hand side is equal to:
\[u_{2\alpha+\beta,2c_{\alpha,\alpha+\beta}}u_{3\alpha+\beta,4c_{\alpha,\alpha+\beta,2,1}}u_{3\alpha+2\beta,2c_{\alpha,\alpha+\beta,1,2}}(u_{2\alpha+\beta,c_{\alpha,\alpha+\beta}}u_{3\alpha+\beta,c_{\alpha,\alpha+\beta,2,1}}u_{3\alpha+2\beta,c_{\alpha,\alpha+\beta,1,2}})^{-2}\]
\[=u_{3\alpha+\beta,2c_{\alpha,\alpha+\beta,2,1}},\]
hence the first assertion. The second one is obtained symmetrically. $\Box$
\end{proof}

\begin{prop}\label{cabij4}
We have:
\[c_{\alpha,\beta,2}=\frac 12c_{\alpha,\beta}c_{\alpha,\alpha+\beta};\]
\[c_{\alpha,\beta,3}=\frac 16c_{\alpha,\beta}c_{\alpha,\alpha+\beta}c_{\alpha,2\alpha+\beta};\]
\[c_{\alpha,\beta,3,2}=\frac 13c_{\alpha+\beta}^2c_{\alpha,\alpha+\beta}c_{\alpha+\beta,2\alpha+\beta}.\]
\end{prop}

\begin{proof}
According to proposition \ref{comimbr}, we have:
\begin{equation}\label{g2blah}[u_{\alpha,1},[u_{\alpha,1},u_{\beta,1}]]=[u_{\alpha,2},u_{\beta,1}][u_{\alpha_1},u_{\beta,1}]^{-2}.\end{equation}
The left-hand side is equal to:
\[[u_{\alpha,1},u_{\alpha+\beta,c_{\alpha,\beta}}u_{2\alpha+\beta,c_{\alpha,\beta,2}}u_{3\alpha+\beta,c_{\alpha,\beta,3}}u_{3\alpha+2\beta,c_{\alpha,\beta,3,2}}].\]
Since both ${\mth{U}}_{3\alpha+\beta}$ and ${\mth{U}}_{3\alpha+2\beta}$ commute with ${\mth{U}}_\alpha$, we deduce from lemma \ref{prodcom} that we have:
\[[u_{\alpha,1},[u_{\alpha,1},u_{\beta,1}]]=[u_{\alpha,1},u_{\alpha+\beta,c_{\alpha,\beta}}u_{2\alpha+\beta,c_{\alpha,\beta,2}}]\]
\[=[u_{\alpha,1},u_{\alpha+\beta,c_{\alpha,\beta}}][u_{\alpha,1},u_{2\alpha+\beta,c_{\alpha,\beta,2}}]\]
\[=u_{2\alpha+\beta,c_{\alpha,\beta}c_{\alpha,\alpha+\beta}}u_{3\alpha+\beta,c_{\alpha,\beta}c_{\alpha,\alpha+\beta,2,1}}u_{3\alpha+2\beta,c_{\alpha,\beta}^2c_{\alpha,\alpha+\beta,1,2}}u_{3\alpha+\beta,c_{\alpha,\beta,2}c_{\alpha,2\alpha+\beta}}\]
\[=u_{2\alpha+\beta,c_{\alpha,\beta}c_{\alpha,\alpha+\beta}}u_{3\alpha+\beta,c_{\alpha,\beta}c_{\alpha,\alpha+\beta,2,1}+c_{\alpha,\beta,2}c_{\alpha,2\alpha+\beta}}u_{3\alpha+2\beta,c_{\alpha,\beta}^2c_{\alpha,\alpha+\beta,1,2}}.\]
The right-hand side is equal to:
\[u_{\alpha+\beta,2c_{\alpha,\beta}}u_{2\alpha+\beta,4c_{\alpha,\beta,2}}u_{3\alpha+\beta,8c_{\alpha,\beta,3}}u_{3\alpha+2\beta,8c_{\alpha,\beta,3,2}}\]
\[*(u_{\alpha+\beta,c_{\alpha,\beta}}u_{2\alpha+\beta,c_{\alpha,\beta,2}}u_{3\alpha+\beta,c_{\alpha,\beta,3}}u_{3\alpha+2\beta,c_{\alpha,\beta,3,2}})^{-2}.\]
Now we develop the above expression and rearrange its terms; while doing so, we must take into account the fact that ${\mth{U}}_{\alpha+\beta}$ and ${\mth{U}}_{2\alpha+\beta}$ do not commute; more precisely, by proposition \ref{cswap} we have:
\begin{equation}\label{g21121}u_{2\alpha+\beta,-c_{\alpha,\beta,2}}u_{\alpha+\beta,-c_{\alpha,\beta}}=u_{\alpha+\beta,-c_{\alpha,\beta}}u_{2\alpha+\beta,-c_{\alpha,\beta,2}}u_{3\alpha+2\beta,-c_{\alpha,\beta,2}c_{\alpha,\beta}c_{\alpha+\beta,2\alpha+\beta}}.\end{equation}
We finally obtain (details are left to the reader) that the right-hand side of (\ref{g2blah}) is equal to:
\[u_{2\alpha+\beta,2c_{\alpha,\beta,2}}u_{3\alpha+\beta,6c_{\alpha,\beta,3}}u_{3\alpha+2\beta,6c_{\alpha,\beta,3,2}-5c_{\alpha+\beta,2\alpha+\beta}c_{\alpha,\beta}c_{\alpha,\beta,2}},\]
from which we deduce, with the help of proposition \ref{cabij3}:
\[2c_{\alpha,\beta,2}=c_{\alpha,\beta}c_{\alpha,\alpha+\beta};\]
\[6c_{\alpha,\beta,3}=c_{\alpha,\beta}c_{\alpha,\alpha+\beta,2,1}+c_{\alpha,\beta,2}c_{\alpha,2\alpha+\beta}\]
\[=\frac 12(c_{\alpha,\beta}c_{\alpha,\alpha+\beta}c_{\alpha,2\alpha+\beta}+c_{\alpha,\beta}c_{\alpha,\alpha+\beta}c_{\alpha,2\alpha+\beta})=c_{\alpha,\beta}c_{\alpha,\alpha+\beta}c_{\alpha,2\alpha+\beta};\]
\[6c_{\alpha,\beta,3,2}-5c_{\alpha+\beta,2\alpha+\beta}c_{\alpha,\beta}c_{\alpha,\beta,2}=c_{\alpha,\beta}^2c_{\alpha,\alpha+\beta,1,2},\]
hence:
\[6c_{\alpha,\beta,3,2}=5c_{\alpha+\beta,2\alpha+\beta}c_{\alpha,\beta}c_{\alpha,\beta,2}+c_{\alpha,\beta}^2c_{\alpha,\alpha+\beta,1,2}\]
\[=\frac 52c_{\alpha+\beta,2\alpha+\beta}c_{\alpha,\beta}^2c_{\alpha,\alpha+\beta}-\frac 12c_{\alpha,\beta}^2c_{\alpha,\alpha+\beta}c_{\alpha+\beta,2\alpha+\beta}\]
\[=2c_{\alpha+\beta,2\alpha+\beta}c_{\alpha,\beta}^2c_{\alpha,\alpha+\beta}.\]
The proposition follows immediately. $\Box$
\end{proof}

\begin{prop}\label{cabij5}
We have:
\[c_{\beta,\alpha,2,3}=-c_{\alpha,\beta,3,2}-\frac 12c_{\alpha,\beta}^2c_{\alpha,\alpha+\beta}c_{\alpha+\beta,2\alpha+\beta}.\]
\[=-\frac 16c_{\alpha,\beta}^2c_{\alpha,\alpha+\beta}c_{\alpha+\beta,2\alpha+\beta}.\]
\end{prop}

\begin{proof}
We have:
\[[u_{\beta,1},u_{\alpha,1}]=u_{\alpha+\beta,c_{\beta,\alpha}}u_{2\alpha+\beta,c_{\beta,\alpha,2}}u_{3\alpha+\beta,c_{\beta,\alpha,3}}u_{3\alpha+2\beta,c_{\beta,\alpha,2,3}}.\]
On the other hand, we have:
\[[u_{\beta,1},u_{\alpha,1}]=[u_{\alpha,1},u_{\beta,1}]^{-1}\]
\[=u_{3\alpha+2\beta,-c_{\alpha,\beta,3,2}}u_{3\alpha+\beta,-c_{\alpha,\beta,3}}u_{2\alpha+\beta,-c_{\alpha,\beta,2}}u_{\alpha+\beta,-c_{\alpha,\beta}}.\]
The assertion is then an immediate consequence of proposition \ref{cabij4} and (\ref{g21121}). $\Box$
\end{proof}

Next we observe that the long roots $\pm\beta,\pm(3\alpha+\beta),\pm(3\alpha+2\beta)$ form a closed root subsystem of type $A_2$ of $\Phi$; the following result follows then immediately from the case $A_2$:

\begin{prop}\label{com1g2}
We have:
\[c_{\beta,3\alpha+\beta}=c_{3\alpha+\beta,-3\alpha-2\beta}=c_{-3\alpha-2\beta,\beta};\]
\[c_{\beta,3\alpha+\beta}c_{3\alpha+2\beta,-3\alpha-\beta}=1;\]
\[c_{\beta,3\alpha+\beta}c_{-\beta,-3\alpha-\beta}=-1.\]
Moreover, for every $\lambda\in R^*$ and every $\mu\in R_{3\alpha+\beta}$, we have $Ad(h_{\beta^\vee,\lambda})u_{3\alpha+\beta,\mu}=u_{\alpha,\frac\mu\lambda}$.
For every $\lambda\in R\*$, we also have $h_{(3\alpha+2\beta)^\vee,\lambda}=h_{\alpha^\vee,\lambda}h_{(3\alpha+\beta)^\vee,\lambda}$.
\end{prop}

Note that, by proposition \ref{crs}, we are still allowed to assume $f_0(\alpha)=f_0(\beta)=0$.

\begin{prop}\label{com2g2}
We have:
\[c_{\alpha,\beta}c_{\alpha+\beta,-\beta}=1.\]
Moreover, for every $\lambda\in R^*$ and every $\mu\in R_\alpha$, we have $Ad(h_{\beta^\vee,\lambda})u_{\alpha,\mu}=u_{\alpha,\frac\mu\lambda}$.
\end{prop}

\begin{proof}
For every suitable $\lambda$, we have:
\[[u_{\alpha,1},u_{\beta,1}]u_{-\beta,\lambda}=u_{\alpha+\beta,c_{\alpha,\beta}}u_{2\alpha+\beta,-c_{\alpha,\beta,2}}u_{3\alpha+\beta,c_{\alpha,\beta,3}}u_{3\alpha+2\beta,-c_{\alpha,\beta,3,2}}u_{-\beta,\lambda}\]
\[=u_{\alpha+\beta,c_{\alpha,\beta}}u_{-\beta,\lambda}u_{2\alpha+\beta,-c_{\alpha,\beta,2}}u_{3\alpha+\beta,c_{\alpha,\beta,3}-c_{\alpha,\beta,3,2}c_{3\alpha+2\beta,-\beta}\lambda}u_{3\alpha+2\beta,-c_{\alpha,\beta,3,2}}\]
\[=u_{-\beta,\lambda}u_{\alpha+\beta,c_{\alpha,\beta}}u_{\alpha,c_{\alpha,\beta}c_{\alpha+\beta,-\beta}\lambda}u',\]
where $u'$ is an element of ${\mth{U}}_{2\alpha+\beta}{\mth{U}}_{3\alpha+\beta}{\mth{U}}_{3\alpha+2\beta}$. On the other hand, we have:
\[[u_{\alpha,1},u_{\beta,1}]u_{-\beta,\lambda}=u_{\alpha,1}u_{\beta,1}u_{\alpha,-1}u_{\beta,-1}u_{-\beta,\lambda}\]
\[=u_{-\beta,\lambda}u_{\alpha,1}h_{\beta^\vee,\frac 1{1-\lambda}}u_{\beta,1-\lambda}u_{\alpha,-1}h_{\beta^\vee,1-\lambda}u_{\beta,-\frac 1{1-\lambda}}.\]
The above expression being identical to (\ref{ca2bl}), we can once again finish the proof the same way as for proposition \ref{com2}. $\Box$
\end{proof}

\begin{prop}\label{com3g2}
We have:
\[c_{\beta,\alpha}c_{\alpha+\beta,-\alpha}=3;\]
Moreover, for every $\lambda\in R^*$ and every $\mu\in R_\beta$, we have $Ad(h_{\alpha^\vee,\lambda})u_{\beta,\mu}=u_{\beta,\frac\mu{\lambda^3}}$.
\end{prop}

\begin{proof}
For every suitable $\lambda$, we have:
\[[u_{\beta,1},u_{\alpha,1}]u_{-\alpha,\lambda}=u_{\alpha+\beta,-c_{\alpha,\beta}}u_{2\alpha+\beta,c_{\alpha,\beta,2}}u_{3\alpha+\beta,-c_{\alpha,\beta,3}}u_{3\alpha+2\beta,c_{\beta,\alpha,3,2}}u_{-\alpha,\lambda}.\]
Obviously, for every $\gamma$ of the form $i\alpha+j\beta$, $i,j>0$ and for every $i',j'\geq 0$ such that $i'(-\alpha)+j'\gamma$ is a root, that root is either $-\alpha$ or of the form $U_{i''\alpha+j''\beta}$, with $i''\geq 0$ and $j''>0$. We thus obtain:
\[[u_{\beta,1},u_{\alpha,1}]u_{-\alpha,\lambda}=u_{-\alpha,\lambda}u_{\beta,c_{\beta,\alpha}c_{\alpha+\beta,-\alpha}\lambda-c_{\beta,\alpha,2}c_{2\alpha+\beta,-\alpha,2,1}\lambda^2+c_{\beta,\alpha,3}c_{3\alpha+\beta,\-\alpha,3}\lambda^3}u',\]
where $u'$ is a product of elements of ${\mth{U}}_{i\alpha+j\beta}$, $i,j>0$. On the other hand, we have:
\[[u_{\beta,1},u_{\alpha,1}]u_{-\alpha,\lambda}=u_{\beta,1}u_{\alpha,1}u_{-\beta,1}u_{-\alpha_1}u_{-\alpha,\lambda}\]
\[=u_{-\alpha,\lambda}u_{\beta,1}h_{\alpha^\vee,\frac 1{1-\lambda}}u_{\alpha,1-\lambda}u_{\beta,-1}h_{\alpha^\vee,1-\lambda}u_{\alpha,\frac 1{\lambda-1}}\]
By a similar reasoning as in propositions \ref{com2} and \ref{com2b2}, we obtain:
\[Ad(h_{\alpha^\vee,\frac 1{1-\lambda}})u_{\beta,1}=u_{\beta,1+c_{\beta,\alpha}c_{\alpha+\beta,-\alpha}\lambda-c_{\beta,\alpha,2}c_{2\alpha+\beta,-\alpha,2,1}\lambda^2+c_{\beta,\alpha,3}c_{3\alpha+\beta,-\alpha,3}\lambda^3}.\]
Since the coefficient of the right-hand side is a monomial in $1-\lambda$, it must be $(1-\lambda)^3$. We finish the proof in a similar way as for proposition \ref{com2b2}. $\Box$
\end{proof}

\begin{prop}\label{prodh3}
For every $\lambda\in R^*$, we have $h_{\alpha^\vee+\beta^\vee,\lambda}=h_{\alpha(^\vee,\lambda}h_{\beta^\vee,\lambda}$.
\end{prop}

\begin{proof}
Since $\alpha$ and $\beta$ generate a root system of type $G_2$, we have $\alpha^\vee+\beta^\vee=(3\alpha+\beta)^\vee$. For every $\lambda\in R^*$, let $\mu,\nu\in R^*$ be such that $h_{\alpha^\vee+\beta^\vee,\lambda}=h_{\alpha^\vee,\mu}h_{\beta^\vee,\nu}$; by examining the actions of both members on respectively ${\mth{U}}_\alpha$, ${\mth{U}}_\beta$ and ${\mth{U}}_{3\alpha+\beta}$, we obtain, using propositions \ref{com2g2} and \ref{com3g2}:
\[\lambda=\mu^2\nu^{-1};\;\;\lambda^{-1}=\mu^{-3}\nu^2\;\;\lambda^2=\mu^3\nu^{-1}.\]
We deduce from the first and third relations that $\lambda=\mu$, and from the second and third ones that $\lambda=\nu$, which proves the proposition. $\Box$
\end{proof}

\begin{prop}\label{com4g2}
We have:
\[c_{\alpha+\beta,\alpha}c_{2\alpha+\beta,-\alpha}=4;\]
\end{prop}

\begin{proof}
As in proposition \ref{com3g2}, for every root $\gamma$ of the form $i\alpha+j\beta$, $i,j>0$ and for every $i',j'\geq 0$ such that $i'(-\alpha)+j'\gamma$ is a root, that root is either $-\alpha$ or of the form $U_{i''\alpha+j''\beta}$, with $i''\geq 0$ and $j''>0$. For every suitable $\lambda$, with the help of equation \ref{blu} (see the proof of proposition \ref{com3g2}), we obtain:
\[[u_{\alpha+\beta,1},u_{\alpha,1}]u_{-\alpha,\lambda}=u_{2\alpha+\beta,c_{\alpha+\beta,\alpha}}u_{3\alpha+\beta,-c_{\alpha+\beta,\alpha,1,2}}u_{3\alpha+2\beta,-c_{\alpha+\beta,\alpha,2,1}}u_{-\alpha,\lambda}\]
\[=u_{-\alpha,\lambda}u_{\alpha+\beta,c_{\alpha+\beta,\alpha}c_{2\alpha+\beta,-\alpha}\lambda-c_{3\alpha+\beta,\alpha,1,2}\lambda^2+\lambda^2x}u',\]
where $x$ is a polynomial in $\lambda$ and $u'$ is a product of elements of ${\mth{U}}_\gamma$ with $\gamma\in\{\beta,2\alpha+\beta,3\alpha+\beta,3\alpha+2\beta\}$. On the other hand, we have:
\[[u_{\alpha+\beta,1},u_{\alpha,1}]u_{-\alpha,\lambda}=u_{\alpha+\beta,1}u_{\alpha,1}u_{\alpha+\beta,-1}u_{\alpha,-1}u_{-\alpha,\lambda}\]
\[=u_{-\alpha,\lambda}u_{\alpha+\beta,1}u_{\beta,c_{\alpha+\beta,-\alpha}\lambda}h_{\alpha^\vee,\frac 1{1-\lambda}}u_{\alpha,1-\lambda}u_{\alpha+\beta,-1}u_{\beta,c_{\alpha+\beta,-\alpha}\frac\lambda{\lambda-1}}h_{\alpha^\vee,1-\lambda}u_{\alpha,\frac 1{\lambda-1}}\]
\[=u_{-\alpha,\lambda}u_{\alpha+\beta,1}u_{\beta,c_{\alpha+\beta,-\alpha}\lambda}u_{\alpha,\frac 1{1-\lambda}}u_{\alpha+\beta,\lambda-1}u_{\beta,(\lambda-1)^2c_{\alpha+\beta,-\alpha}\lambda}u_{\alpha,\frac 1{\lambda-1}}\]
\[=u_{-\alpha,\lambda}u_{\alpha+\beta,1}u_{\alpha,\frac 1{1-\lambda}}u_{\alpha+\beta,\lambda-1+\frac{c_{\alpha,\beta}c_{\alpha+\beta,-\alpha}\lambda}{1-\lambda}}u_{\beta,(1+(\lambda-1)^2)c_{\alpha+\beta,-\alpha})\lambda}u_{\alpha,\frac 1{\lambda-1}}u'\]
\[=u_{-\alpha,\lambda}u_{\alpha+\beta,1}u_{\alpha,\frac 1{1-\lambda}}u_{\alpha+\beta,\lambda-1-\frac {3\lambda}{1-\lambda}}u_{\alpha,\frac 1{\lambda-1}}u_{\alpha+\beta,\frac{(1+(\lambda-1)^2)c_{\alpha+\beta,-\alpha}c_{\beta,\alpha}\lambda}{\lambda-1}}\]
\[*u_{\beta,(1+(\lambda-1)^2)c_{\alpha+\beta,-\alpha}\lambda}u'',\]
where $u'$ and $u''$ are elements of ${\mth{U}}_{2\alpha+\beta}{\mth{U}}_{3\alpha+\beta}{\mth{U}}_{3\alpha+2\beta}$; we have used proposition \ref{com3g2} in the process. Applying it once again, we see that the term in $u_{\alpha+\beta}$ we finally obtain is $u_{\alpha+\beta,x}$, with:
\[x=\lambda-\frac{3\lambda}{1-\lambda}-\frac{3(-1-(\lambda-1)^2)\lambda}{1-\lambda}\]
\[=1+3\lambda(1-\lambda)=4\lambda-3\lambda^2\]
Hence we have $c_{\alpha+\beta,\alpha}c_{2\alpha+\beta,-\alpha}=4$ and the proposition is proved. $\Box$
\end{proof}

\begin{prop}\label{prodh4}
For every $\lambda\in R^*$, we have $h_{(2\alpha+\beta)^\vee,\lambda}=h_{\alpha(^\vee,\lambda}h_{(\alpha+\beta)^\vee,\lambda}$.
\end{prop}

\begin{proof}
The three roots involved are all short, and generate a (nonclosed) subsystem of type $A_2$ of $\Phi$; the proof is then similar to the proof of proposition \ref{prodh}. $\Box$
\end{proof}

\begin{prop}\label{cyclg2}
We have:
\[c_{\alpha,\beta}=c_{\beta,-\alpha-\beta}=\frac 13c_{-\alpha-\beta,\alpha};\]
\[c_{\alpha,\alpha+\beta}=c_{\alpha+\beta,-2\alpha-\beta}=c_{-2\alpha-\beta,\alpha}.\]
\end{prop}

\begin{proof}
The first assertion can be proved in a very similat way as the corresponding one in the case $B_2$. Details are left to the reader.

For the second one, since the roots we are considering generate a subsystem of type $A_2$ of $\Phi$, we will proceed in a similar way as for proposition \ref{cycla2}; however, that subsystem is not closed, which leads to additional terms showing up while developing our commutators. For that reason, instead of $[u_{\alpha,1},u_{\beta,1}]u_{-\alpha-\beta,\lambda}$, we will work with $[u_{\alpha,\lambda},u_{\beta,\lambda}]u_{-\alpha-\beta,\mu}$ for every suitable $\lambda,\mu$. We now have to compare the terms in $\lambda^2\mu$ in the expressions we obtain, and since all additional terms, which come from the long roots showing up in the development of $[u_{\alpha,\lambda},u_{\beta,\lambda}]$ and similar commutators, turn out to contribute only to terms of higher degrees in either $\lambda$ or $\mu$, the equality we finally obtain is the same as in the proof of proposition \ref{cycla2}, and the result follows. $\Box$
\end{proof}

\begin{cor}\label{pabg2}
We have:
\[c_{\alpha,\beta}c_{-\alpha,-\beta}=-1;\]
\[c_{\alpha,\alpha+\beta}c_{-\alpha,-\alpha-\beta}=-4;\]
\[c_{\alpha,2\alpha+\beta}c_{-\alpha,-2\alpha-\beta}=-9.\]
\end{cor}

\begin{proof}
Using propositions \ref{com2g2} and \ref{cyclg2} and a similar remark as in the case $A_2$ (see the proof of corollary \ref{paba2}), we obtain:
\[c_{\alpha,\beta}c_{-\alpha,-\beta}=-\frac 13c_{\beta,\alpha}c_{\alpha+\beta,-\alpha}=-1.\]
The proof of the other two assertions is similar, using proposition \ref{com4g2} as well when needed. $\Box$
\end{proof}

For practical use, we will now write some of the relations of the last three sections in a more condensed form:

\begin{prop}\label{pab}
For every $\alpha,\beta\in\Phi$ such that $\alpha+\beta$ is a root, there exists an integer $p_{\alpha,\beta}$ such that $c_{\alpha,\beta}c_{-\alpha,-\beta}=-p_{\alpha,\beta}^2$. Moreover, we have:
\[\frac{c_{\alpha,\beta}}{p_{\alpha,\beta}}=\frac{c_{\beta,-\alpha-\beta}}{p_{\beta,-\alpha-\beta}}=\frac{c_{-\alpha-\beta,\alpha}}{p_{-\alpha-\beta,\alpha}}.\]
\end{prop}

\begin{proof}
The proposition is an immediate consequence of corollaries \ref{paba2}, \ref{pabb2} and \ref{pabg2} for the first assertion and of propositions \ref{cycla2}, \ref{cyclb2} and \ref{cyclg2} for the second assertion. $\Box$
\end{proof}

It can be checked that, as in \cite{chev}, $p_{\alpha,\beta}$ is the smallest integer $c$ such that $\beta-c\alpha$ is not a root (or rather that we can choose $p_{\alpha,\beta}$ to be such an integer, since the $c_{\alpha,\beta}$ do not necessary live in rings of characteristic zero), but we will not use this fact in the sequel. We can still assume $p_{\alpha,\beta}\in\{1,2,3\}$ for every $\alpha,\beta$, and thus depends only on the relative position of $\alpha$ and $\beta$ since we have assumed $p\neq 2$.

\subsection{The unicity result}

Now we prove that the $c_{\alpha,\beta}$ correspond up to normalization (and up to the fact that some of them are elements of a quotient of $R$ and not of $R$ itself) to the Chevalley constants for reductive groups. The fact that those constants satisfy the previous propositions is proved in \cite{chev} for most of them, and can be easily deduced from the results of that paper for the remaining ones; the converse follows immediately from the following unicity result. Note that since the $c_{\alpha,\beta,i,j}$, $i+j\geq 3$, are entirely determined by the $c_{\alpha,\beta}$ (propositions \ref{cabij2}, \ref{cabij3}, \ref{cabij4} and \ref{cabij5}), we do not have to worry about them.

\begin{prop}\label{comuniq}
Let $(c'_{\alpha,\beta})_{\alpha,\beta}$ be another family of nonzero constants such that:
\begin{itemize}
\item for every $\alpha,\beta\in\Phi$ such that $\alpha+\beta\in\Phi$, $c'_{\alpha,\beta}$ lives in the same ring as $c_{\alpha,\beta}$;
 \item the $c'_{\alpha,\beta}$ also satisfy all the results we have proved about the $c_{\alpha,\beta}$ in the previous four subsections (namely propositions \ref{cswitch}, \ref{com3}, \ref{com2}, \ref{cycla2}, \ref{com1b2}, \ref{com2b2}, \ref{cyclb2}, \ref{com1g2}, \ref{com2g2}, \ref{com3g2}, \ref{com4g2}, \ref{cyclg2} and \ref{pab}).
\end{itemize}
There exists then a family of nonzero constants $(N_\alpha)_{\alpha\in\Phi}$ satisfying the following conditions:
\begin{itemize}
\item for every $\alpha\in\Phi$, $N_{-\alpha}=N_\alpha^{-1}$;
\item for every $\alpha,\beta$, $c'_{\alpha,\beta}=\frac{N_\alpha N_\beta}{N_{\alpha+\beta}}c_{\alpha,\beta}$.
\end{itemize}
\end{prop}

\begin{proof}
Consider the system of affine roots $\Phi'=\Phi\times{\mth{Z}}$; let ${\mth{B}}$ be any Borel subgroup of ${\mth{G}}$ containing ${\mth{T}}$ and let $\Phi'_{\mth{B}}$ be the subset of the elements $(\alpha,n)$ of $\Phi'$ such that $n\geq f_{\mth{B}}(\alpha)$, where $f_{\mth{B}}$ is the concave function associated to ${\mth{B}}$; since $f_{\mth{B}}$ is concave, $\Phi'_{\mth{B}}$ is a closed subset of $\Phi'$. Consider the partial order on $\Phi'_{\mth{B}}$ defined by $(\alpha,n)\leq(\beta,m)$ if there exists a pair $(\gamma,l)$ which is a sum of elements of $\Phi'_{\mth{B}}$ and such that $\beta=\alpha+\gamma$ and $m=n+l$, and let $\Delta'$ be the set of the minimal elements of $\Phi'_{\mth{B}}$ for that order; the set $\Delta_0=\{\alpha|(\alpha,n)\in\Delta'\;\mathrm{for}\;\mathrm{some}\;n\}$ is then an extended set of simple roots of $\Phi$.

Let $\Delta$ be any (standard) set of simple roots of $\Phi$ contained in $\Delta_0$ and let $\Phi^+$ be the set of positive roots of $\Phi$ generated by $\Delta$. By minimality of the $(\alpha,f(\alpha))$, $\alpha\in\Delta$, and an easy induction, for every $\beta\in\Phi^+$, writing $\beta=\sum_{\alpha\in\Delta}\lambda_\alpha\alpha$, we have:
\[f_{\mth{B}}(\beta)=\sum_{\alpha\in\Delta}\lambda_\alpha f_{\mth{B}}(\alpha).\]
This implies in particular that for every $\alpha,\beta,\gamma\in\Phi^+$ such that $\alpha+\beta+\gamma\in\Phi^+$, $f_{\mth{B}}(\alpha+\beta+\gamma)=f_{\mth{B}}(\alpha)+f_{\mth{B}}(\beta)+f_{\mth{B}}(\gamma)$.

On the other hand, for every $\alpha\in\Phi$, we have either $f_0(\alpha)=f_{\mth{B}}(\alpha)$ or $f_0(\alpha)=f_{\mth{B}}(\alpha)-1$ and $\alpha\in\Psi$. We deduce from this, using the concavity of $f_0$, that for every $\alpha,\beta,\gamma$ as above, either $f_0(\alpha+\beta+\gamma)=f_0(\alpha)+f_0(\beta)+f_0(\gamma)$ or $f_0(\alpha+\beta+\gamma)=f_0(\alpha)+f_0(\beta)+f_0(\gamma)-1$ and $\alpha+\beta+\gamma\in\Psi$. Hence we are allowed to apply proposition \ref{com3} to every $\alpha,\beta,\gamma\in\Phi^+$ satisfying the other required conditions.

We first prove the existence of constants $N_\alpha$, $\alpha\in\Phi^+$, such that the second condition is satisfied for every $\alpha,\beta\in\Phi^+$ such that $\alpha+\beta\in\Phi$; when either $\alpha$ or $\beta$ (say $\alpha$) is not in $\Psi$, the constant $c_{\alpha,\beta}$ then lives in $R_{\alpha+\beta}$. (If both $\alpha$ and $\beta$ belong to $\Psi$, $c_{\alpha,\beta}$ always lives in $R_{\alpha+\beta}$ anyway). For every $\alpha\in\Phi^+$, let $h_\alpha$ be the height of $\alpha$, that is the number of elements of $\Delta$ (counted with multiplicities) $\alpha$ is the sum of; we define the $N_\alpha$ by induction on $h_\alpha$. For every $\alpha\in\Delta$, we can choose $N_\alpha$ arbitrarily. Now assume $\alpha$ is an element of $\Delta$ of height $2$; there exists then a unique pair $\beta,\gamma$ of elements of $\Delta$ such that $\beta+\gamma=\alpha$, and we can set:
\[N_\alpha=\frac{c_{\beta,\gamma}N_\beta N_\gamma}{c'_{\beta,\gamma}}.\]
(By proposition \ref{cswitch}, the order in which we take $\beta$ and $\gamma$ does not matter.) Now assume $\alpha$ is of height $>2$, and let $\beta\in\Phi$ and $\gamma\in\Delta$ be such that $\alpha=\beta+\gamma$; if $\alpha\in\Psi$, we also assume $\gamma\in\Psi$ whenever possible. Define $N_\alpha$ as above; if $\alpha$ belongs to $\Psi$ and $\beta$ and $\gamma$ do not, we can choose $N_{\alpha}$ arbitrarily among the elements of $R$ whose image in $R_\beta$ is the right-hand side of the above equality. (This case only occurs when $\alpha$ is the smallest positive root in the irreducible component of $\Psi$ containing the highest root of $\Phi$ relative to $\Delta$.)

Let now $\beta',\gamma'$ be any two elements of $\Phi^+$ such that $\alpha=\beta'+\gamma'$; we will check that $c'_{\beta',\gamma'}$ satisfies the required condition. According to \cite[\S 1, proposition 19]{bou}, switching $\beta'$ and $\gamma'$ if necessary, we can assume $\delta=\beta'-\gamma$ is either a root or $0$. In the second case there is nothing to prove ;.in the first case, since $\gamma$ is a simple root, $\delta$ must be positive. We then have $\beta=\delta+\gamma'$ and, by proposition \ref{com3}:
\[c'_{\beta',\gamma'}c'_{\delta,\gamma}=c'_{\beta,\gamma}c'_{\delta,\gamma'}\]
from which we deduce, using the induction hypothesis:
\[c'_{\beta',\gamma'}\frac{N_\gamma N_\delta}{N_{\beta'}}c_{\delta,\gamma}=c'_{\beta,\gamma}\frac{N_\delta N_{\gamma'}}{N_\beta}c_{\delta,\gamma'}\]
\[=\frac{N_\delta N_{\gamma'}N_\gamma}{N_\alpha}c_{\beta,\gamma}c_{\delta,\gamma'}\]
hence:
\[c'_{\beta',\gamma'}c'_{\delta,\gamma}=\frac{N_{\beta'}N_{\gamma'}}{N_\alpha}c_{\beta,\gamma}c_{\delta,\gamma'}.\]
The assertion then follows from proposition \ref{com3}.

Now consider the case of negative roots. For any $\alpha\in-\Phi^+$, set $N_\alpha=N_{-\alpha}^{-1}$. Let $\alpha,\beta$ be both negative and such that $\alpha+\beta\in\Phi$; we then have, using the first assertion of proposition\ \ref{pab}:
\[c'_{\alpha,\beta}=\frac{- p_{\alpha,\beta}^2}{c'_{-\alpha,-\beta}}\]
\[=\frac{p_{\alpha,\beta}^2N_{-\alpha-\beta}}{N_{-\alpha}N_{-\beta}c_{-\alpha,-\beta}}\]
\[=\frac{N_\alpha N_\beta}{N_{\alpha+\beta}}c_{\alpha,\beta},\]
as required.

Finally assume $\alpha$ is positive and $\beta$ negative. Assume $\alpha+\beta>0$, hence $-\alpha-\beta<0$; we have, using the second assertion of proposition \ref{pab}:
\[c'_{\alpha,\beta}=\frac{p_{\alpha,\beta}^2}{p_{\beta,-\alpha-\beta}^2}c'_{\beta,-\alpha-\beta}=\frac{N_\beta N_{-\alpha-\beta}p_{\alpha,\beta}^2}{N_{-\alpha}p_{\beta,-\alpha-\beta}^2}c_{\beta,-\alpha-\beta}=\frac{N_\alpha N_\beta}{N_{\alpha+\beta}}c_{\alpha,\beta},\]
as required. The case $\alpha+\beta<0$ being similar, the proof is now complete. $\Box$
\end{proof}

\section{The main theorem}
\
Now we come to the main result of the paper. Let $k$ be any perfect field of characteristic $p\neq 2$ (eventually zero), let $d$ be a nonnegative integer and let $\underline{\mth{G}}$ be a connected group of parahoric type of depth $d$ defined over $k$. Let $\underline{\Phi}$ (resp. $\Phi$ be the absolute (resp. relative) root system of $\underline{\mth{G}}$ relative to some given maximal $k$-torus $\underline{\mth{T}}$ and let ${\mth{G}}$  be the group of $k$-points of $\underline{\mth{G}}$; we assume that $\Phi$ is irreducible, and that $\underline{\Phi}$ has no components of type either $A_2$ or $G_2$ when $p=3$. Note that we do not require $\underline{\Phi}$ to be irreducible. We will also denote by $\underline{\Psi}$ (resp. $\Psi$) the absolute (resp. relative) root system of the reductive part of $\underline{\mth{G}}$.

Remember that $\underline{\mth{H}}$ is the Cartan subgroup of $\underline{\mth{G}}$ containing $\underline{\mth{T}}$, or in other words the centralizer of $\underline{\mth{T}}$ in $\underline{\mth{G}}$. For every root $\alpha\in\underline{\Phi}$, we denote by $\underline{\mth{U}}_\alpha$ the root subgroup of $\underline{\mth{G}}$ associated to $\alpha$, and by $\underline{R}_\alpha$ the truncated valuation ring associated to $\underline{\mth{U}}_\alpha$ as in subsection $3.2$. If $\varpi$ is any uniformizer of $\underline{R}_\alpha$, we also set $\underline{R'}_\alpha=\underline{R}_\alpha/\varpi^d\underline{R}_\alpha$. We also set $\underline{\mth{T}}_\alpha=Im(\alpha^\vee)$, and set $\underline{\mth{H}}_\alpha$ to be the subgroup of $\underline{\mth{H}}$ generated by $\underline{\mth{T}}_\alpha$ and the intersection $\underline{\mth{H}}_{\alpha_1}$ of $\underline{\mth{H}}$ with the subgroup of $\underline{\mth{G}}$ generated by $\underline{\mth{U}}_\alpha$ and $\underline{\mth{U}}_{-\alpha}$.

We also denote respectively by ${\mth{G}}$, ${\mth{T}}$, ${\mth{H}}$, ${\mth{U}}_\alpha$, $R_\alpha$, $R'_\alpha$, ${\mth{T}}_\alpha$ and ${\mth{H}}_\alpha$ the groups or rings of $k$-points of $\underline{\mth{G}}$, $\underline{\mth{T}}$, $\underline{\mth{H}}$, $\underline{\mth{U}}_\alpha$, $\underline{R}_\alpha$, $\underline{R'}_\alpha$, $\underline{\mth{T}}_\alpha$ and $\underline{\mth{H}}_\alpha$ when these varous groups and rings happen to be defined over $k$.

Let $\underline{\mth{S}}$ be a maximal $k$-split torus of $\underline{\mth{G}}$; we can assume $\underline{\mth{T}}$ contains $\underline{\mth{S}}$. As in \cite{cou}, a group of parahoric type is said to be {\em quasisplit} over $k$ if it satisfies the following equivalent conditions:

\begin{itemize}
\item the centralizer of $\underline{\mth{S}}$ is a Cartan subgroup $\underline{\mth{H}}$ of $\underline{\mth{G}}$;
\item let $k'$ be a Galois extension of $k$ on which $\underline{\mth{G}}$ splits, and let $\Gamma=Gal(k'/k)$. There exists a $\Gamma$-stable set of positive roots in $\underline{\Phi}$;
\item with $\Gamma$ defined as above, there existe a $\Gamma$-stable pseudo-Borel subgroup of $\underline{\mth{G}}$ containing $\underline{\mth{T}}$.
\end{itemize}

In the sequel, we will make the additional assumptions that  $\underline{\mth{G}}$ is quasisplit and satisfies the following condition: let $\underline{\Psi}$ be the root system of its reductive part; then all the rings $\underline{R}_\alpha$, $\alpha\in\underline{\Psi}$ (resp. $\underline{R}'_\alpha$, $\alpha\in\underline{\Phi}-\underline{\Psi}$) are isomorphic. By proposition \ref{isoconn2} (resp. corollary \ref{isoconn}), the first condition (resp. the second condition) is true as soon as $\underline{\Psi}$ (resp. $\underline{\Phi}$) is connected, and we will see in the next section that it is even enough to assume that $\Psi$ (resp; $\Phi$) is connected.

 To avoid problems with the center $Z(\underline{\mth{G}})$ of $\underline{\mth{G}}$ we have to add the following condition: let $r$ (resp. $r'$) be the semisimple rank (resp. the rank) of $\underline{\mth{G}}$. Then the unipotent radical of $Z(\underline{\mth{G}})$ is the direct product of its intersection with the product of the $\underline{\mth{H}}_\alpha$ and of $r'-r$ copies of the group $1+\varpi\underline{R}_\alpha$, where $\varpi$ is a uniformizer of $R_\alpha$. This is in particular always true when $\underline{\mth{G}}$ is of semisimple parahoric type, since $r'-r=0$ in that case.

\begin{theo}\label{main}
Assume $\underline{\mth{G}}$ is connected and satisfies the above conditions. There exists a nonarchimedean henselian local field $F$ whose residual field is $k$ and a connected reductive group $\underline{G}$ defined and quasisplit over $k$ such that , if $F_{nr}$ is the maximal unramified algebraic extension of $F$ and $G_{nr}$ is the group of $F_{nr}$-points of $\underline{G}$, then $\underline{\mth{G}}$ is $k$-isomorphic to the quotient of a parahoric subgroup of $G_{nr}$ defined over $F$ by the product of its $d$-th congruence subgroup and a finite group.
\end{theo}

\begin{proof}
First assume that $\underline{\mth{G}}$ is split over $k$ (hence $\Phi=\underline{\Phi}$, $\Psi=\underline{\Psi}$ and for every $\alpha\in\Phi$, $\underline{R}_\alpha$ is defined over $k$) and that every element of $k$ admits two opposite square roots (which is in particular the case when $k$ is algebraically closed). By proposition \ref{exfield}, there exists a henselian local field $F$ such that if $\mathcal{O}_F$ is the integer ring of $F$ and $\mathfrak{p}_F$ is the maximal ideal of $\mathcal{O}_F$, for some given $\alpha\in\Phi$, $R'_\alpha$ is isomorphic to $\mathcal{O}_F/\mathfrak{p}_F^d$. We deduce from the hypotheses that for every $\alpha\in\Phi$, the ring $R'_\alpha$ is isomorphic to $\mathcal{O}_F/\mathfrak{p}_F^d$. When $\Psi$ is nonempty, this is also true for the rings $R_\alpha$ of $k$-points of the $\underline{R}_\alpha$, $\alpha\in\Psi$, with $d$ replaced by $d+1$.

Let $(X,\Phi,X^\vee,\Phi^\vee)$ be the root datum of $\underline{\mth{G}}$. According to the theorem 10.1.1 of \cite{spr}, there exists a connected reductive algebraic group $\underline{G}$, defined and split over $F$, whose root datum is also $(X,\Phi,X^\vee,\Phi^\vee)$. Let $\underline{T}$ be a $k$-split maximal torus of $\underline{G}$; for every $\alpha\in\underline{\Phi}=\Phi$, let $\underline{U_\alpha}$ be the root subgroup associated to $\alpha$ and let $v=(v_\alpha)$ be a valuation (in the sense of \cite[Section I.6.2]{bt}) on $\underline{T}$ and the $\underline{U_\alpha}$; since $\underline{G}$ is $k$-split we can assume that the image of the group $T$ (resp. $U_\alpha$ for every $\alpha$) of $F$-points of $\underline{T}$ (resp. $\underline{U_\alpha}$ for every $\alpha$) by $v_\alpha$ is ${\mth{Z}}\cup\{\infty\}$.

Let $F_{nr}$ be the maximal unramified algebraic extension of $F$; $F_{nr}$ is an henselian local field with discrete valuation whose residual field is $\overline{k}$. Let $G_{nr}$ (resp. $T_{nr}$, $U_{\alpha,nr}$ for every $\alpha$) be the group of $F_{nr}$-points of $\underline{G}$ (resp. $\underline{T}$, $\underline{U}_\alpha$ for every $\alpha$). For every $\alpha$ and every $i\in{\mth{Z}}$, let $U_{\alpha,nr,i}$ be the subgroup of $U_{\alpha,nr}$ defined the usual way; according to \cite[proposition 2.6]{cou2}, the group $K_{nr}$ generated by the maximal bounded subgroup $K_{T_{nr}}$ of $T_{nr}$ and the $U_{\alpha,nr,f_0(\alpha)}$, with $f_0$ being the concave function defined in $({\bf PT4})$, is then a parahoric subgroup of $G_{nr}$.

Let $K_{nr}^d$ be the $d$-th congruence subgroup of $K_{nr}$; we prove now that $\underline{\mth{G}}$ and $K_{nr}/K_{nr}^d$ are $k$-isomorphic. We already know that we have $k$-isomorphisms between $\underline{\mth{U}}_\alpha$ and $U_{\alpha,nr,f_0(\alpha)}/U_{\alpha,nr,d+1-f_0(-\alpha)}$ for every $\alpha$; moreover, we see by the unicity result (proposition \ref{comuniq}) that these isomorphisms can be chosen in such a way that the constants $c_{\alpha,\beta}$ are the images of the corresponding constants for $G_{nr}$ in the quotient rings in which they live; since these constants are elements of subquotient rings of $F$ anyway, it is clearly compatible with the fact that the isomorphisms are defined over $k$.

We can also deduce from the propositions of the previous sections (propositions \ref{prodh}, \ref{prodh2}, \ref{com1g2}, \ref{prodh3} and \ref{prodh4}) a $k$-isogeny from $K_{T_{nr}}/K_{T_{nr}}^d$ to ${\mth{H}}$; let $K'$ be its finite kernel.

Set $K=K_{nr}\cap G$ and $K^d=K_{nr}^d\cap G$; $K$ is then a parahoric subgroup of $G$ and $K^d$ is its $d$-th congruence subgroup. We only have to find an isomorphism between ${\mth{G}}$ and ${\mth{G}}_0=K/K^d$; the isomorphism between $\underline{\mth{G}}$ and $\underline{\mth{G}}_0$ will then be given by the particular case $k=\underline{k}$.

Assume first ${\mth{G}}$ is solvable.
 We can then write ${\mth{G}}={\mth{H}}\prod_{\alpha\in\Phi}{\mth{U}}_\alpha$, with the $\alpha$ being taken in some fixed (arbitrarily chosen) order. Moreover, the reductive group ${\mth{G}}/R_u({\mth{G}})$ is solvable too, hence $\Psi$ must be empty, which implies by proposition \ref{sfun} that $f_0(\alpha)+f_0(-\alpha)=1$ for every $\alpha$; hence the reductive group $K/K^0$ is a torus, which implies that $K$ is an Iwahori subgroup of $G$. As an Iwahori subgroup, $K$ admits an Iwahori decomposition, which is similar to the above decomposition of ${\mth{G}}$. Using those decompositions, the above isomorphisms then yield a bijection $\phi$ from ${\mth{G}}$ to $K/K'K^d$, which is also an isomorphism of $k$-varieties; it thus only remains to prove that $\phi$ preserves the multiplication. By an obvious induction we can see that it is enough to prove that $\phi(gg')=\phi(g)\phi(g')$ when $g$ is any element of ${\mth{G}}$ and $g'$ belongs either to ${\mth{H}}$ or to one of the ${\mth{U}}_\alpha$. Note also that we can replace ${\mth{G}}$ by any one of its subgroups which is a product of root subgroups and eventually of ${\mth{H}}$, and ${\mth{G}}_0$ by a suitable subgroup, and the following reasoning still works. 

First we assume that $g'\in{\mth{H}}$. Since we then have ${\mth{U}}_\alpha g'=g'{\mth{U}}_\alpha$ for every $\alpha$ it is enough to prove the assertion when $g\in{\mth{U}}_\alpha$ for some $\alpha$. But since, as a consequence of propositions \ref{com2}, \ref{com1b2}, \ref{com2b2}, \ref{act1b2}, \ref{com1g2}, \ref{com2g2} and \ref{com3g2}, ${\mth{H}}$ and $K_T/K'K_T^d$ act on respectively ${\mth{U}}_\alpha$ and its image by the respective characters $\alpha$ and $\phi\circ\alpha\circ\phi^{-1}$, we obtain $\phi(g'^{-1}gg')=\phi(g'^{-1})\phi(g)\phi(g')$. On the other hand, it is immediate from the definition of $\phi$ that $\phi(g'^{-1}gg')=\phi(g'^{-1})\phi(gg')$, hence the desired assertion.

Now we assume that $g'\in{\mth{U}}_\alpha$ for some $\alpha=\alpha_{g'}$. For every nonnegative integer $i$, let $\mathcal{C}_i$ be the normal subgroup of $R_u({\mth{G}})$ defined inductively by $\mathcal{C}_0=R_u({\mth{G}})$ and $\mathcal{C}_i=[R_u({\mth{G}}),\mathcal{C}_{i-1}]$ for every $i\geq 1$. Since ${\mth{G}}$ is solvable, by \cite[corollary 6.3.3]{spr}, $R_u({\mth{G}})$ is nilpotent, or in other words the groups $\mathcal{C}_i$ are trivial for $i$ large enough. Moreover, since ${\mth{H}}$ normalizes $R_u({\mth{G}})$, it also normalizes every $\mathcal{C}_i$; since ${\mth{H}}$ contains ${\mth{T}}$, for every $i$, the group $\mathcal{C}_i$ is then the product of some subgroup of ${\mth{H}}$ and of unipotent groups of the form ${\mth{U}}_{\alpha,j}$, which means that if $g=h\prod_\alpha u_\alpha$ belongs to $\mathcal{C}_i$, then so do $h$ and every $u_\alpha$. For every $g\in G$, we denote by $c_g$ the largest integer $i$ such that $g\in\mathcal{C}_i$; by convention $c_1=+\infty$.

Of course, $\phi(gg')=\phi(g)\phi(g')$ holds if either $g$ or $g'$ is the identity. Moreover, let $\beta=\beta_g$ be the greatest root (for the order we have used in the definition of $\phi$) such that the component of $g$ in ${\mth{U}}_\beta$ is nontrivial; when $g\in{\mth{H}}$, we can set $\beta=0$ and consider that $0<\alpha$ for every root $\alpha$. When $\beta\leq\alpha$, the fact that $\phi(gg')=\phi(g)\phi(g')$ is immediate from the definition of $\phi$. We now prove the general case by descending induction on $i=Sup(c_g,c_{g'})$.

If $i$ is large enough, then either $g$ or $g'$ is trivial. Assume this is not the case, and write $g=g_1u_\beta$, where $u_\beta$ is the component of $g$ in ${\mth{U}}_\beta$; we will also assume $\beta>\alpha$. We then have:
\[\phi(g)\phi(g')=\phi(g_1)\phi(u_\beta)\phi(g')\]
\[=\phi(g_1)\phi(g')\phi(u_\beta)\phi([u_\beta^{-1},g'^{-1}]).\]
The fact that $\phi(u_\beta)\phi(g')=\phi(g')\phi(u_\beta)\phi([u_\beta^{-1},g'^{-1}])$ comes from the case of rank 1 solvable groups (proposition \ref{soca}) when $\beta=-\alpha$ and from the commutator relations of the previous section when $\alpha$ and $\beta$ are linearly independent.  Now if we assume that $\phi(g_1)\phi(g')=\phi(g_1g')$, we obtain:
\[\phi(g)\phi(g')=\phi(g_1g')\phi(u_\beta)\phi([u_\beta^{-1},g'^{-1}]).\]
Since every root $\gamma$ such that the component of $g_1g'$ in ${\mth{U}}_\gamma$ is nontrivial is $\leq\beta$, we have $\phi(g_1g')\phi(u_\beta)=\phi(g_1g'u_\beta)$. Moreover, since either $u_\beta$ or $g'$ lies in $\mathcal{C}_i$, $[u_\beta^{-1},g'^{-1}]$ belongs to $\mathcal{C}_{i+1}$, and we can apply the induction hypothesis to see that $\phi(g_1g'u_\beta)\phi([u_\beta^{-1},g'^{-1}])=\phi(g_1g'u_\beta[u_\beta^{-1},g'^{-1}])=\phi(gg')$, as desired.

We are thus reduced to proving the assertion for $g_1$ and $g'$, where $g_1$ belongs to $\mathcal{C}_i$ whenever $g$ does, and the largest root $\beta_1$ such that $g_1$ has a component in ${\mth{U}}_{\beta_1}$ is strictly smaller than $\beta$. By an obvious induction, after repeating the process a finite number of times we reach a $g_1$ such that $\beta_1$ is smaller than $\alpha$, which proves the assertion. Hence the theorem holds when ${\mth{G}}$ is solvable.

Now we assume ${\mth{G}}$ is not solvable. Although ${\mth{G}}$ does not admit an Iwahori decomposition anymore, the subset ${\mth{H}}\prod_{\alpha\in\Phi}{\mth{U}}_\alpha$, which of course depends on the order on $\Phi$ that we have chosen, is a dense open subset of ${\mth{G}}$ for any choice of that order. We will thus be able to use proposition \ref{isoo} during the proof of the theorem.

Assume that, with some arbitrary choice of a set of positive roots $\Phi^+$ in $\Phi$, every negative root is smaller than every positive root for our order; we then have, using the fact that ${\mth{H}}$ normalizes every root subgroup, ${\mth{H}}\prod_\alpha{\mth{U}}_\alpha={\mth{U}}^-{\mth{B}}$, where ${\mth{B}}$ is the pseudo-Borel subgroup of ${\mth{G}}$ generated by ${\mth{H}}$ and the ${\mth{U}}_\alpha$, $\alpha>0$, and ${\mth{U}}^-$ is the unipotent subgroup of ${\mth{G}}$ generated by the ${\mth{U}}_\alpha$, $\alpha<0$. Let $\phi$ be the bijection between ${\mth{U}}^-{\mth{B}}$ and the corresponding subset of $K/K'K^d$ defined as in the solvable case. Since ${\mth{U}}^-$ and ${\mth{B}}$ are solvable, the restriction of $\phi$ to ${\mth{U}}^-$ (resp. ${\mth{B}}$) is an isomorphism, which does not depend on the restriction of our order to $\Phi^-$ (resp. $\Phi^+$) since its restriction to the generating subset $\bigcup_{\alpha\in\Phi^-}{\mth{U}}_\alpha$ of ${\mth{U}}^-$ (resp. the generating subset ${\mth{H}}\cup\bigcup_{\alpha\in\Phi^+}{\mth{U}}_\alpha$ of ${\mth{B}}$) does not either. Hence for every $g=ub\in{\mth{U}}^-{\mth{B}}$, $\phi(ub)=\phi(u)\phi(b)$ does not depend on the choice of the order on $\Phi$ as long as it satisfies the required condition. Moreover, let ${\mth{B}}_0={\mth{B}}{\mth{U}}_1^-$ be the unique Borel subgroup of ${\mth{G}}$ containing ${\mth{B}}$ (recall that ${\mth{B}}$ is a pseudo-Borel subgroup, not a true Borel subgroup); since ${\mth{B}}_0$ is the group of $k$-points of a solvable group of parahoric type defined over $k$, we know from the previous case that for every $g,g'\in{\mth{B}}_0$, $\phi(gg')=\phi(g)\phi(g')$. Now assume $g=ub$ is any element of ${\mth{U}}^-{\mth{B}}$, with $g'$ still in ${\mth{B}}_0$. We then have:
\[\phi(g)\phi(g')=\phi(u)\phi(b)\phi(g')=\phi(u)\phi(bg').\]
Since $bg'\in{\mth{B}}_0$, we can write $bg'=u'b'$, with $u'\in{\mth{U}}^-$ and $b'\in{\mth{B}}$. Since ${\mth{U}}^-$ is solvable, we obtain:
\[\phi(u)\phi(u'b')=\phi(u)\phi(u')\phi(b')=\phi(uu')\phi(b')=\phi(uu'b')=\phi(gg').\]
By a similar reasoning, we obtain that we also have $\phi(gg')=\phi(g)\phi(g')$ when $g$ is an element of the only Borel subgroup ${\mth{B}}_0^-$ of ${\mth{G}}$ containing ${\mth{H}}$ and opposite to ${\mth{B}}$, and $g'$ is any element of ${\mth{U}}^-{\mth{B}}$.

Let $\phi'$ be the bijection between ${\mth{B}}{\mth{U}}^-$ and the corresponding subset of $K/K'K^d$ defined the same way as $\phi$, but using the opposite order on $\Phi$. Since $\phi$ and $\phi'$ coincide on ${\mth{H}}$ and on every ${\mth{U}}_\alpha$, they also coincide on every solvable subgroup of ${\mth{G}}$ containing ${\mth{H}}$. We prove now that they coincide in fact on a dense open subset of ${\mth{G}}$.

Let $\alpha$ be a simple root of $\Phi^+$, and set $\Phi^+_\alpha=\Phi^+\cup\{-\alpha\}-\{\alpha\}$.

\begin{lemme}
The set $\Phi^+_\alpha$ is a set of positive roots of $\Phi$.
\end{lemme}

\begin{proof}
Let $W_\Phi$ be the Weyl group of $\Phi$ and let $s_\alpha\in W_\Phi$ be the reflection associated to $\alpha$. Since $\alpha$ is a simple root of $\Phi^+$, the length of $s_\alpha$ (relative to $\Phi^+$) is $1$, hence $\alpha$ is the only positive root such that $s_\alpha(\alpha)$ is negative. Hence $\Phi^+_\alpha=s_\alpha(\Phi^+)$ and the lemma follows. $\Box$
\end{proof}

Now we go back to the proof of theorem \ref{main}. We will assume that the order on $\Phi$ has been chosen in such a way that $\alpha$ (resp. $-\alpha$)  is the smallest positive (resp. largest negative) root for that order. Let ${\mth{G}}_\alpha$ be the group generated by ${\mth{U}}_\alpha$, ${\mth{H}}$ and ${\mth{U}}_{-\alpha}$; for every $g\in{\mth{U}}^-{\mth{B}}$, if we write $g=ug_\alpha u'$, with $g_\alpha\in{\mth{G}}_\alpha$ and $u$ (resp. $u'$) belonging to the product of the root subgroups associated to negative (resp. positive) roots distinct from $\pm\alpha$, we deduce from the previous remarks that we have:
\[\phi(g)=\phi(u)\phi(g_\alpha)\phi(u').\]

Now consider an order on $\Phi$ satisfying similar properties with $\Phi^+$ replaced by $\Phi^+_\alpha$ and $\alpha$ and $-\alpha$ switched, and let $\phi_\alpha$ be the corresponding bijection between the suitable open dense subset of ${\mth{G}}$ and the corresponding subset of $K/K'K^d$; if $g$ belongs to the domain of definition of $\phi_\alpha$, we also have:
\[\phi_\alpha(g)=\phi_\alpha(u)\phi_\alpha(g_\alpha)\phi_\alpha(u').\]
Assume first that $u=u'=1$, or in other words that $g\in {\mth{G}}_\alpha$; by either proposition \ref{soca} or proposition \ref{a1nr} depending on whether ${\mth{G}}_\alpha$ is solvable or not, both $\phi$ and $\phi_\alpha$ extend to an isomorphism between ${\mth{G}}_\alpha$ and the corresponding subgroup of $K/K'K^d$; on the other hand, they coincide on ${\mth{H}}$, ${\mth{U}}_\alpha$ and ${\mth{U}}_{-\alpha}$, which generate ${\mth{G}}_\alpha$, hence they must coincide on the whole intersection of their domains of definition. The assertion for ${\mth{G}}$ id an immediate consequence of the assertion for ${\mth{G}}_\alpha$, proposition \ref{isoo} and the above two equalities.

We can iterate the process, replacing positive roots of $\Psi$ by negative roots one at a time,, and after a finite number of steps we reach a situation where $\Phi^+$ has been replaced by $-\Phi^+$ and $\phi$ by $\phi'$; by an obvious induction, we obtain that $\phi$ and $\phi'$ coincide on a dense open subset $S$ of ${\mth{G}}$, which is the intersection of the domains of definition of all the $\phi$-like maps  we have defined and used in the process. In particular, we have $S\subset{\mth{U}}^-{\mth{B}}$.

Let now $g=ub$, $g'=u'b'$ be two elements of ${\mth{U}}^-{\mth{B}}$; such that $bu\in S$. Then in particular $bu\in{\mth{U}}^-{\mth{B}}$, hence $gg'\in{\mth{U}}^-{\mth{B}}$; on the other hand, sonce $bu\in{\mth{B}}{\mth{U}}^-$, and $\phi'$ uses the order on $\Phi$ which is opposite to the one used by $\Phi$ , we have $\phi'(bu)=\phi'(b)\phi'(u)$, and we obtain:
\[\phi(gg')=\phi(ubu'b')=\phi(u)\phi(bu')\phi(b')\]
\[=\phi(u)\phi'(bu')\phi(b')\]
\[=\phi(u)\phi'(b)\phi'(u')\phi(b')\]
\[=\phi(u)\phi(b)\phi(u')\phi(b')=\phi(g)\phi(g').\]
Setting ${\mth{L}}={\mth{G}}$, ${\mth{L}}'=K/K'K^d$, $\Omega={\mth{U}}^-{\mth{B}}$ and $\Omega''=\{(ub,u'b')\in\Omega\times\Omega|bu'\in S\}$, we finally use proposition \ref{isoo} to see that $\phi$ extends to an isomorphism between ${\mth{G}}$ and $K/K'K^d$; the fact that $\Omega''$ satisfies the required condition can be easily checked by taking $b=u'=1$. The assertion of theorem \ref{main} is now proved.

Assume now $k$ is any perfect field (up to the conditions on its characteristic), with $\underline{\mth{G}}$ still split over $k$. As in the previous case, there exists a nonarchimedean local field $F_{nr}$ whose residual field is $\overline{k}$ and a connected reductive group $\underline{G}$ defined over $F_{nr}$ such that, if $G_{nr}$ is the group of $F_{nr}$-points of $\underline{G}$, $\underline{\mth{G}}$ is isomorphic to the quotient of a parahoric subgroup $K_{nr}$ of $G_{nr}$ by its $d$-th congruence subgroup $K_{nr}^d$, via an isomorphism $\phi$. Moreover, if $F$ is any henselian local field associated to the rings $R_\alpha$ by proposition \ref{exfield}, we can assume that $F_{nr}$ is the maximal unramified extension of $F$.

If $\underline{G}$ is defined over $F$ and $\phi$ is defined over $k$, then $\phi({\mth{G}})$ is simply the group of $k$-points of $K_{nr}/K'_{nr}K_{nr}^d$, with $K'_{nr}$ being a suitable finite group, and the result follows immediately; we thus only have to check that these hypotheses actually hold. The isomorphism $\phi$ is defined by its restrictions to $\underline{\mth{H}}$ and the $\underline{\mth{U}}_\alpha$ and all these groups as well as their images are defined over $k$; moreover, consider the extended root datum $(X^*(\underline{T}),\Phi,X_*(\underline{T}),\Phi^\vee,\Delta_0,\tau)$, where $\Delta_0$ is a set of simple roots of the anisotropic part of ${\mth{G}}$ (see \cite[section 2.2]{tits}) and $\tau$ is the index of ${\mth{G}}$ (see \cite[section 2.3]{tits}): since we are dealing here with a split group, we have $\Delta_0=\emptyset$ and $\tau=Id$, hence according to \cite[chapter 17]{spr}, we can always choose $\underline{G}$ in such a way that it is defined and split over $F$, and in that case, since $F_{nr}/F$ is unramified, the Bruhat-Tits building of the group $G$ of $F$-points of $\underline{G}$ is a subcomplex of the building of $G_{nr}$ containing at least one apartment of that building, which implies in particular that every parahoric subgroup of $G_{nr}$ is conjugate to a $Gal(F_{nr}/F)$-stable one; we can then always assume that $K_{nr}$ satisfies that property. As an immediate consequence, we can thus always construct $\phi$ in such a way that it is defined over $k$, which proves the assertion.

Now we go to the general case. We can assume that $\underline{\mth{T}}$ contains a maximal $k$-split torus $\underline{\mth{S}}$ of $\underline{\mth{G}}$; moreover, let $k_0$ be the smallest Galois extension of $k$ such that $\underline{\mth{G}}$ splits over $k$; as for reductive groups, $k_0$ is the Galois extension corresponding to the largest normal subgroup of $Gal(\overline{k}/k)$ fixing $X(\underline{\mth{T}})$ pointwise, hence it is unique. Let $F_0$ be the finite unramified extension of $F$ whose residual field is $k_0$, let $G_0$ be the group of $F_0$-points of $\underline{G}$, let $K_{nr}$ be defined as in the previous case, set $K_0=G_0\cap K_{nr}$ and define $K'_0$ in a similar way as the groups $K'$ and $K'_{nr}$ of the previous cases; $K_0/K'_0K_0^d$ is then the group of $k_)$-points of $K_{nr}/K'_{nr}K_{nr}^d$. Set $\Gamma=Gal(k_0/k)$; $\Gamma$ is a finite group acting on the root datum $(X^*(\underline{T}),\underline{\Phi},X_*(\underline{T})\underline{,\Phi}^\vee)$.

Let $\alpha$ be any element of $\underline{\Phi}$ and let $\Gamma_\alpha$ be the subgroup of the elements of $\Gamma$ which stabilize $\alpha$; we can always set the unit element $1$ in $R_\alpha$ in such a way that it belongs to the group of $\Gamma_\alpha$-fixed points of ${\mth{U}}_\alpha$, and choose the unit element of every ${\mth{U}}_{\gamma(\alpha)}$ to be $\gamma(1)$. By making similar choices in $K_0/K'_0K_0^d$, we obtain:
\[\phi(\prod_{\gamma\in\Gamma}\gamma(1))=\prod_{\gamma\in\Gamma}\gamma(\phi(u)_1).\]
Let $\Phi$ be the relative root system of ${\mth{G}}$, or in other words the root system of ${\mth{G}}$ relative to the group of $k$-points ${\mth{S}}$ of $\underline{\mth{S}}$. Since $\underline{\mth{G}}$ is quasisplit, we deduce from \cite[Corollary 3.7]{bot} that the restriction to ${\mth{S}}$ of every $\alpha\in\underline{\Phi}$ is nontrivial, hence belongs to $\Phi$. By the remark after \cite[7.1]{bot}, for every $\beta\in\Phi$, there are two possible cases: either $\beta$ is nonmultipliable, which means that no $n\beta$; with $n$ being a positive integer different from $1$, belongs to $\Phi$, or $\beta$ is multipliable, in which case the only such $n$ is $n=2$. The following properties of respectively nonmultipliable and multipliable roots are respectively the cases I and II of \cite[II. 4.1.4]{bt}.

If $\beta$ is a nonmultipliable element of $\Phi$, the elements of the relative root subgroup ${\mth{U}}_\beta$ of ${\mth{G}}$ are precisely the $\prod_{\gamma\in\Gamma}\gamma(u)$, $u\in{\mth{U}}_\alpha$. Hence $\phi$ sends ${\mth{U}}_\beta$ to the corresponding relative root subgroup of $K/K'K^d$. Assume now $\beta$ is multipliable, and let $\alpha$ be an element of $\underline{\Phi}$ whose restriction to ${\mth{S}}$ is $\beta$. Then $\Gamma_\alpha$ is a subgroup of $\Gamma$ of index $2n$, where $n$ is the number of irreducible components of $\underline{\Phi}$ containing elements of the inverse image of $\beta$. Let $\Gamma'_\alpha$ be the subgroup of $\Gamma$ stabilizing each one of these components; $\Gamma_\alpha$ is then a subgroup of $\Gamma'_\alpha$ of index $2$, hence normal in $\Gamma'_\alpha$, and if $\gamma$ is the nontrivial element of $\Gamma'_\alpha/\Gamma_\alpha$, $\alpha+\gamma(\alpha)$ is a $\Gamma'_\alpha$-stable element of $\underline{\Phi}$ and we have:
\[u_{\alpha,1}u_{\gamma(\alpha),1}=u_{\alpha+\gamma(\alpha),c_{\alpha,\gamma(\alpha)}}u_{\gamma(\alpha),1}u_{\alpha,1},\]
with $c_{\alpha,\gamma(\alpha)}$ being an element of the extension $k_\alpha$ of $k$ such that $Gal(k_0/k_\alpha)=\Gamma_\alpha$. By applying $\gamma$ on both sides and comparing the new equality with the previous one, we obtain:
\[c_{\alpha,\gamma(\alpha)}+\gamma(c_{\gamma(\alpha),\alpha})=0.\]
Since $c_{\gamma(\alpha),\alpha}=-c_{\alpha,\gamma(\alpha)}$, we see that it simply means that $c_{\alpha,\gamma(\alpha)}$ is an element of the fixed field $k'_\alpha$ of $\Gamma_\alpha$, hence $u_{\alpha+\gamma(\alpha),c_{\alpha,\gamma(\alpha)}}\in\underline{\mth{G}}(k'_\alpha)$. On the other hand, we have for every $x,y\in R_\alpha$:
\[\prod_{\gamma'\in\Gamma_\alpha}\gamma'(\gamma(u_{\alpha+\gamma(\alpha),x}u_{\alpha,y}u_{\gamma(\alpha),\gamma(y)}))=\prod_{\gamma'\in\Gamma_\alpha}\gamma'(u_{\alpha+\gamma(\alpha),\gamma(x)-c_{\alpha,\gamma(\alpha)}y\gamma(y)}u_{\alpha,y}u_{\gamma(\alpha),\gamma(y)}),\]
and the elements of ${\mth{U}}_\beta$ are precisely the elements of that form such that $\prod_{\gamma'\in\Gamma_\alpha}\gamma'\gamma(x)-x=\prod_{\gamma'\in\Gamma_\alpha}\gamma'(c_{\alpha,\gamma(\alpha)}y\gamma(y))$. Since this is true for the corresponding subgroup of $K/K'K^d$ as well, that subgroup must be the image of ${\mth{U}}_\beta$ by $\phi$.

Now we have to find a $k$-automorphism between $\underline{\mth{H}}$ and the corresponding Cartan subgroup $\underline{\mth{H}}_0$ of the reductive group of which $K_0/K'_0K_0^d$ is the group of $k_0$-points. We have $\underline{\mth{H}}=\underline{\mth{T}}R_u(\underline{\mth{H}})$; write also $\underline{\mth{H}}_0=\underline{\mth{T}}_0R_u(\underline{\mth{H}}_0)$, where $\underline{\mth{T}}_0$ is the unique maximal torus of $\underline{\mth{H}}_0$. Since $\underline{\mth{G}}$ and $\underline{\mth{G}}_0$ have isomorphic root data and $\Gamma$ acts the same way on both, $\underline{\mth{T}}$ and $\underline{\mth{T}}_0$ must be $k$-isomorphic. On the other hand, since ${\mth{H}}$ and ${\mth{H}}_0$ are the centralizers of $k$-isomorphic maximal tori of respectively ${\mth{G}}$ and ${\mth{G}}_0$, they must be $k$-isomorphic as well; hence we see using proposition \ref{basprod} that $R_u(\underline{\mth{H}})$ and $R_u(\underline{\mth{H}}_0)$ are abelian unipotent $k$-groups both $k$-isomorphic to the product of $r'$ copies of the group $1+\underline{R}_{\alpha,1}$, hence their respective groups of $k_0$-points
%can be viewed as representations of $\Gamma$ in finite-dimensional $k$-vector spaces. Since both these representations are isomorphic to the direct sum of $d$ copies of the representation of $\Gamma$ in $Lie(\underline{\mth{T}})$ yielded by its action on $\Phi$, they
are isomorphic to each other, which implies tnat $R_u(\underline{\mth{H}})$ and $R_u(\underline{\mth{H}}_0)$ are $k$-isomorphic. We thus obtain that $\underline{\mth{H}}$ and $\underline{\mth{H}}_0$ are $k$-isomorphic as well.

Since ${\mth{G}}$ is generated by ${\mth{H}}$ and the ${\mth{U}}_\beta$, $\beta\in\Phi$, we see that with the choice we have made $\phi$ must be defined over $k$ as well, which completes the proof of the theorem.
\end{proof}

\section{Some particular cases}

In this section, we give a few examples of cases in which the result of the previous section always holds, as well as an explicit example of a group of parahoric type which does not satisfy the required conditions.

\begin{prop}\label{car0}
Assume $k$ is of characteristic $0$. Then all the $\underline{R}_\alpha$, $\alpha\in\underline{\Psi}$ (resp. $\alpha\not\in\underline{\Psi}$), are $k$-isomorphic.
\end{prop}

\begin{proof}
According to \cite[II, theorem 2]{ser}, the only complete local field with discrete valuation admitting $k$ as its residual field is $k((X))$. Since every henselian local subfield of $k((X))$ yields the same quotient rings as $k((X))$, there is only one possible choice for $\underline{R}_\alpha$ when $\alpha\in\underline{\Psi}$ (resp. when $\alpha\not\in\underline{\Psi}$), and the result follows immediately. $\Box$
\end{proof}

\begin{prop}
Assume $\Phi$ and $\Psi$ are irreducible or empty. Then all the $\underline{R}_\alpha$, $\alpha\in\underline{\Psi}$ (resp. $\alpha\not\in\underline{\Psi}$), are $k$-isomorphic.
\end{prop}

\begin{proof}
Let $k'$ be a Galois extension of $k$ over which $\underline{\mth{G}}$ splits, and set $\Gamma=Gal(k'/k)$; we deduce from lemma 16.2.8 of \cite{spr} that the fact that $\Phi$ (resp. $\Psi$) is irreducible means that the action of $\Gamma$ on the irreducible components of $\underline{\Phi}$ (resp. $\underline{\Psi}$) is transitive, and since for every $\alpha\in\Phi$ and every $\gamma\in\Gamma$, $\underline{R}_{\gamma(\alpha)}=\gamma(\underline{R}_\alpha)$ is $k$-isomorphic to $\underline{R}_\alpha$, we only have to check that for some fixed irreducible component $\underline{\Phi}_1$ of $\underline{\Phi}$, all the $\underline{R}_\alpha$, with $\alpha$ being an element of $\underline{\Phi}_1$ contained (resp. not contained) in $\underline{\Psi}$, are $k$-isomorphic The result for the roots outside $\underline{\Psi}$ is now an immediate consequence of corollary \ref{isoconn}; for roots inside $\underline{\Psi}$ we can use a similar reasoning. $\Box$
\end{proof}

We say that the truncated ring $\underline{R}_\alpha$ is {\em absolutely unramified} if it is of residual characteristic $p>0$ and the valuation of $p$ in $\underline{R}_\alpha$ is $1$. Note that it implies $d\geq 1$.

For every $\alpha\in\underline{\Phi}$, let $\underline{R}'_\alpha$ be the quotient of $\underline{R}_\alpha$ by the ideal of its elements of valuation at least $d$.

\begin{prop}
Assume that $\underline{\Phi}$ is irreducible, $d>1$ and that for some $\alpha\in\underline{\Psi}$, $\underline{R}_\alpha$ is absolutely unramified. Then for every $\beta\in\underline{\Psi}$, $\underline{R}_\alpha$ and $\underline{R}_\beta$ are $k$-isomorphic to each other.
\end{prop}

\begin{proof}
Since $\underline{R}_\alpha$ is absolutely unramified, its characteristic is $p^{d+1}$ for some prime number $p$, and we have $v(p)=1$ in $\underline{R}_\alpha$. For any choice of the field $F_\alpha$ admitting $\underline{R}_\alpha$ as a quotient of its unit ring, that field must be an absolutely unramified local field whose residual field is $k$. According to \cite[II, theorem 3]{ser}, there exists only one such local complete field (up to isomorphism); let $F$ be that field.

On the other hand, since $d>1$, the image of $p$ in $\underline{R}'_{\alpha}$ is nonzero; we deduce then from proposition \ref{anniso3} and corollary \ref{isoconn} that $p$ is also of valuation $1$ in $\underline{R}_\beta$, hence $\underline{R}_\beta$ is absolutely unramified. The fields $F_\alpha$ and $F_\beta$ associated to respectively $\underline{R}_\alpha$ and $\underline{R}_\beta$ by proposition \ref{exfield} must then be henselian subfields of $F$, which implies that $\underline{R}_\alpha$ and $\underline{R}_\beta$ are $k$-isomorphic too. $\Box$
\end{proof}

Note that the above result does not hold when $d=1$.

We can generalize that proposition to the case of tame ramification as follows. We say that the truncated ring $\underline{R}_\alpha$ is {\em absolutely unramified} if it is of residual characteristic $p>0$ and the valuation $e$ of $p$ in $\underline{R}_\alpha$ is not a multiple of $p$. Similarly as before, these hypotheses imply $d\geq e$.

\begin{prop}
Assume that $\underline{\Phi}$ is irreducible and that for some $\alpha\in\underline{\Psi}$, $\underline{R}_\alpha$ is absolutely tamely ramified of ramification index $e<d$. Then for every $\beta\in\underline{\Psi}$ (resp. $\beta\not\in\underline{\Psi}$), $\underline{R}_\beta$ is $k$-isomorphic to $\underline{R}_\alpha$.
\end{prop}

\begin{proof}
Since $\underline{R}_\alpha$ is absolutely ramified of ramification index $e$, its characteristic is a power of some prime number $p$ and we have $v(p)=e$ in $\underline{R}_\alpha$. Hence every field $F_\alpha$ associated to $\underline{R}_\alpha$ by proposition \ref{exfield} is also absolutely ramified of index $e$. Since the ramification is tame, we deduce easily from \cite[II, theorem 3]{ser} and Hensel's lemma applied to a polynomial of the form $X^e-\varpi_\alpha^{-e}\varpi$, where $\varpi_\alpha$ is an uniformizer of $F_\alpha$, that $F_\alpha$ is of the form $F_{nr}[\sqrt[e]{\varpi}]$, where $F_{nr}$ is a henselian subfield of the unique absolutely unramified local field with residual field $\overline{k}$ and $\varpi$ is some uniformizer of $F_{nr}$; moreover, $F_{nr}[\sqrt[e]{\varpi}]$ only depends on the class of $\varpi$ mod $\mathfrak{p}_{nr}^2$, where $\mathfrak{p}_{nr}$ is the maximal ideal of the ring of integers of $F_{nr}$.

On the other hand, since $d>e$, $\underline{R}'_{\alpha}$ is also tamely ramified of index $e$, hence $\underline{R}_\beta$ is too, and there exists a uniformizer $\varpi '$ of $F_{nr}$ which is a $e$-th power of some element of $F_\beta$ and which belongs to $\varpi+\mathfrak{p}_{nr}^2$; since $d>e$, $F_\alpha$ and $F_\beta$ must be henselian subfields of the same complete local field, hence $\underline{R}_\alpha$ and $\underline{R}_\beta$ are $k$-isomorphic. $\Box$
\end{proof}

Note that the above result does not hold when the ramification is wild (that is when the valuation of $p$ is a multiple of $p$).

In these two cases, assuming $\underline{\mth{G}}$ is of semisimple rank at least $2$, the fact that $\underline{\mth{G}}$ matches the conditions of theorem \ref{main} is an immediate consequence of proposition \ref{isoconn}.

Now we give an example of (the group of $k$-points ${\mth{G}}$ of) a group of parahoric type $\underline{\mth{G}}$ with $\Psi$ reducible that does not satisfy the required conditions. Consider the fields $F_1={\mth{Q}}_3[\sqrt{3}]$ and $F_2={\mth{Q}}_3[\sqrt{6}]$; both are totally ramified extensions of ${\mth{Q}}_3$ of ramification index $2$, but they are not isomorphic to each other. For $i=1,2$, let $\mathcal{O}_i$ be the ring of integers of $F_i$, let $\mathfrak{p}_i$ be its maximal ideal, and let $R_i$ be the ring $\mathcal{O}_i/\mathfrak{p}_i^3$. The $R_i$ are not isomorphic to each other, but the quotient rings $R'_i=\mathcal{O}_i/\mathfrak{p}_i^2$ are (via the map sending $\sqrt{3}$ to $\sqrt{6}$); by a slight abuse of notation we can set $R'=R'_1=R'_2$. Moreover, the maximal ideals $I_i$ of the $R_i$ are isomorphic as $R'$-modules (via the map sending $\sqrt{3}$ to $\sqrt{6}$ and $3$ to $6$); by a similar abuse of notation we can also set $I=I_1=I_2$.

Let ${\mth{G}}$ be the set of matrices of the form: $\left(\begin{array}{cc}A&B\\C&D\end{array}\right)$, with $A\in GL_n(R_1)$, $B\in M_n(R')$, $C\in M_n(I)$ and $D\in GL_n(R_2)$. By the remarks above the usual matrix multiplication is well-defined, and it is easy to check that it defines a group law on ${\mth{G}}$. Since that map is obviously algebraic, ${\mth{G}}$ is the group of ${\mth{F}}_3$-points of a split algebraic group $\underline{\mth{G}}$ defined on ${\mth{F}}_3$; checking that $\underline{\mth{G}}$ is of parahoric type, with $\Phi$ of type $A_{2n-1}$ and $\Psi$ of the form $\Psi_1\times\Psi_2$, with for $i=1,2$, $\Psi_i$ being of type $A_{n-1}$, is also straightforward.

Let $\alpha$ (resp. $\beta$) be an element of $\Psi_1$ (resp. $\Psi_2$). Then $R_\alpha$ is isomorphic to $R_1$ and $R_\beta$ is isomorphic to $R_2$. On the other hand, assume ${\mth{G}}$ is a quotient of a parahoric subgroup of a group $G(F)$ for some nonarchimedean local field $F$ of residual field ${\mth{F}}_3$; then every $R_\alpha$, $\alpha\in\Psi$, must be isomorphic to $\mathcal{O}_F/\mathfrak{p}_F^3$, where $\mathcal{O}_F$ is the ring of integers of $F$ and $\mathfrak{p}_F$ its maximal ideal, and we reach a contradiction.

Note that in this particular case, the maximal unramified extensions of $F_1$ and $F_2$ are isomorphic; which means that if $\alpha$ and $\beta$ are two elements of $\underline{\Psi}=\Psi$ not lying in the same irreducible component, $\underline{R}_\alpha$ and $\underline{R}_\beta$ are $\overline{k}$-isomorphic, but the isomorphism between them cannot be defined over $k$. Replacing our fields $F_1$ and $F_2$ by fields of wild absolute ramification, we can even find examples of groups of parahoric type in which the $\underline{R}_\alpha$, $\alpha\in\Psi$, are not even $\overline{k}$-isomorphic.

\end{document}